%% file: dR.tex
\documentstyle[amssymb,thmsd,a4,11pt]{article}
\input tcilatex

\input diagram
\newtheorem{examples}[theorem]{Examples}
\newtheorem{PropDef}[theorem]{Proposition and Definition}

\QQQ{Language}{
American English
}

\begin{document}

\author{Gabriele Vezzosi \\
%EndAName
{\em Dipartimento di Matematica}\\
{\em Universit\`a di Bologna}\\
{\em Piazza di Porta S. Donato 5, 40127}\\
{\em Bologna - Italy}\\
e-mail: {\tt vezzosi@dm.unibo.it}\\
\smallskip\ \\
Alexandre M. Vinogradov \\
%EndAName
{\em Dipartimento di Matematica e Informatica}\\
{\em Universit\`a di Salerno}\\
{\em Via S. Allende} {\em 84081-Baronissi (SA) Italy }\\
{\em and}\\
{\em the Diffiety Institute 45/6, Polevaya st.,}\\
{\em \ Pereslavl-Zalessky, Russia.}\\
{\em \ }e-mail: {\em \ }{\tt vinograd@ponza.dia.unisa.it}}
\title{{\bf On higher order analogues of de Rham cohomology }}
\date{}
\maketitle

\begin{abstract}
If $K$ is a commutative ring and $A$ is a $K$-algebra, for any sequence $%
\sigma $ of positive integers there exists \cite{Vi2} an higher order
analogue ${\bf dR}_\sigma $ of the standard de Rham complex ${\bf dR}\equiv 
{\bf dR}_{\left( 1,...,1,..\right) }$, which can also be defined starting
from suitable (=differentially closed) subcategories of $A-{\bf Mod}$. The
main result of this paper is that the cohomology of ${\bf dR}_\sigma $ does
not depend on $\sigma $, under some smoothness assumptions on the ambient
category. In \cite{VV}, a weaker result was proved by completely different
methods.

Before proving the main theorem we give a rather detailed exposition of all
relevant (for our present purposes) functors of differential calculus on
commutative algebras. This part can be also of an independent interest.

\smallskip
\ \ 
\end{abstract}

\tableofcontents

\section{Introduction}

Higher analogues of the standard de Rham complex were found by one of the
authors about 20 years ago in course of his study of basic functors of
differential calculus (see \cite{Vi1}). More exactly, such an analogue $%
dR_\sigma $ can be associated with any sequence $\sigma $ of positive
integers. The standard de Rham complex turns out to be associated with the
simplest sequence of this kind, namely, with $\sigma =(1,1,...)$.
Differentials of these complexes are natural differential operators of,
generally, higher orders. Complexes $dR_\sigma $`s as well as related
functors they represent play an important role inside differential calculus
in the sense of \cite{Vi1} and as such are worth to be better investigated.
It is very plausible that they may have important applications to the formal
theory of partial differential equations as well as to Secondary Calculus
(see \cite{Vi3}), to mention the most direct ones. Of course, the first
natural question concerning these complexes is : what are their
cohomologies? At the end of eighties the second author jointly with Yu. N.
Torkhov tested some simple cases (unpublished) in order to check a natural
feeling that all complexes $dR_\sigma $ have the same cohomology. After that
in 1994 the authors proved (see \cite{VV-esi} and \cite{VV}) an analogue of
the ''infinitesimal Stokes formula'' for all $dR_\sigma $`s with
non-decreasing $\sigma $, which allows to develop a homotopy techniques
sufficient to show that these complexes are quasi-isomorphic. One year later
we found a somewhat involved proof of this fact in full generality. Its
essentially simplified version is the main result of this paper.

In this paper we present also for the first time a rather detailed
exposition of all relevant for our purposes functors of differential
calculus (sections 1-3) and as such it could be of an independent interest.

A preliminary version of this paper appeared as a Preprint \cite{VV-sns} of
the Scuola Normale Superiore, Pisa.

Finally we note that recently A.M.Verbovetsky sketched in [Ve] an
alternative proof of the main theorem of this paper based on the theory of
compatibility complexes.

\bigskip\ 

\begin{center}
{\bf Notations and Conventions}

\smallskip\ \ 
\end{center}

$K:$ \qquad a commutative ring with unit ;

$A:\qquad $ a commutative $K$-algebra with unit;

$R-{\bf Mod:\qquad }$ the category of $R$-modules for a commutative ring $R$;

${\bf {\rm Diff}}_A:\qquad $the category whose objects are the $A$-modules
and the morphisms are differential operators (Section 1) between them;

${\bf Ens}:\qquad $the category of sets;

$\left[ {\bf C},{\bf C}\right] :\qquad $the category of functors ${\bf %
C\longrightarrow C}$, ${\bf C}$ being a category;

$Ob\left( {\bf C}\right) :\qquad $the objects of ${\bf C}$, ${\bf C}$ being
a category; $C\in Ob\left( {\bf C}\right) $ means that $C$ is an object of $%
{\bf C}$;

$A-{\bf BiMod:\qquad }$ the category of $A$-bimodules, whose objects are
understood as {\it ordered couples} $(P,P^{+})$ of $A$-modules and whose
morphisms are the usual morphisms of bimodules. Note that $P$ and $P^{+}$
coincide as $K$-modules, hence as sets;

$*A-{\bf BiMod:\qquad }$the category of ''biads'' (see Section 1);

{\bf K}$\left( A-{\bf Mod}\right) $ (resp. {\bf K}$\left( K-{\bf Mod}\right) 
$, resp. {\bf K}$\left( {\bf Diff}_A\right) $)$:$\quad \quad the category of
complexes in $A-{\bf Mod}$ (resp. $K-{\bf Mod}$, resp. ${\bf {\rm Diff}}_A$);

${\frak Ens}:\qquad $ the forgetful functor from a concrete category to $%
{\bf Ens}$;

If ${\frak D}$ is a full subcategory of $A-{\bf Mod}$, a functor $T:{\frak D}%
{\bf \longrightarrow }{\frak D}$ will be said {\em strictly representable}
in ${\frak D}$ if it exists $\tau \in Ob\left( {\frak D}\right) $ and a
functorial isomorphism $T\simeq Hom_A\left( \tau ,\cdot \right) $ in $\left[ 
{\frak D},{\frak D}\right] $.

$A-{\bf BiMod}_{{\frak D}}$ (resp. $*A-{\bf BiMod}_{{\frak D}}$, resp. {\bf K%
}$\left( {\bf Diff}_{A,{\frak D}}\right) $) will be the subcategory of $A-%
{\bf BiMod}$ whose objects are couples of objects in ${\frak D}$ (resp. the
subcategory of $*A-{\bf BiMod}$ whose objects are triples of objects in $%
{\frak D}$, resp. the subcategory of {\bf K}$\left( {\bf Diff}_A\right) $
whose objects are complexes of objects in ${\frak D}$).

A sequence $T_1\rightarrow T_2\rightarrow T_3$ of functors $T_i:{\frak D}%
{\bf \longrightarrow }{\frak D}$, $i=1,2,3$, (and functorial morphisms) $%
{\frak D}$ being an abelian subcategory of $A-{\bf Mod}$, is said {\em exact
in }$\left[ {\frak D},{\frak D}\right] $ if it is exact in ${\frak D}$ when
applied to any object of ${\frak D}$.

Let ${\bf N}_{+}^\infty =\underleftarrow{\lim }{\bf N}_{+}^k$ be the set of
infinite sequences of positive integers. If $\sigma \in {\bf N}_{+}^n$ (or $%
\sigma \in {\bf N}_{+}^\infty $) then $\sigma \left( r\right) \doteq \left(
\sigma _1,...,\sigma _r\right) $ for $r\leq n$ (or any $r\in {\bf N}_{+}$).
We denote by ${\bf 1}$ the element $(1,...,1,1,...1,...)\in {\bf N}%
_{+}^\infty $.

\section{Absolute and relative functors}

We recall here in a slightly different way the necessary definitions from 
\cite{Vi1} ,\cite{Vi2},\cite{Kr} (see also \cite{KLV}).

If $P$ and $Q$ are $A$-modules and $a\in A$ we define: 
\[
\begin{array}{c}
\delta _a:Hom_K(P,Q)\longrightarrow Hom_K(P,Q) \\ 
\Phi \longmapsto \left\{ \delta _a\Phi :p\longmapsto \Phi (ap)-a\Phi
(p)\right\} \quad ,\text{ }p\in P
\end{array}
\]
(where we use juxtaposition to indicate both $A$-module multiplications in $%
P $ and $Q$). For each $a\in A$ $\delta _a$ is a homomorphism of $K$%
-modules, and commutativity of $A$ implies that 
\[
\delta _{a_1}\circ \delta _{a_2}=\delta _{a_2}\circ \delta _{a_1}\quad
,\forall a_1,a_2\in A\text{.} 
\]
\ 

\begin{definition}
A $K-${\em differential operator }({\em DO}){\em \ of order} $\leq s$ from
an $A$-module $P$ to another $A$-module $Q$ is an element $\Delta \in
Hom_K(P,Q)$ such that: 
\[
\left[ \delta _{a_0}\circ \delta _{a_1}\circ ...\circ \delta _{a_s}\right]
(\Delta )=0,\text{ }\forall \left\{ a_0,a_1,...,a_s\right\} \subset A. 
\]
\end{definition}

We will write synthetically $\delta _{a_0,\ldots ,a_s}$ for $\delta
_{a_0}\circ \delta _{a_1}\circ ...\circ \delta _{a_s}$.

The {\em set }{\it Diff}$_k(P,Q)$ of differential operators of order $\leq k$
from $P$ to $Q$ is endowed naturally with two different $A$-module
structures:

(i) ({\it Diff}$_k(P,Q),\tau $) $\doteq {\rm Diff}_k(P,Q)$ (left)$,$

$\tau :A\times ${\it Diff}$_k(P,Q)\longrightarrow ${\it Diff}$%
_k(P,Q):(a,\Delta )\longmapsto \tau (a,\Delta ):p\longmapsto a\Delta (p)$

(ii) ({\it Diff}$_k(P,Q),\tau ^{+}$) $\doteq {\rm Diff}_k^{+}(P,Q)$ (right)$%
, $

$\tau ^{+}:A\times ${\it Diff}$_k(P,Q)\longrightarrow ${\it Diff}$%
_k(P,Q):(a,\Delta )\longmapsto \tau ^{+}(a,\Delta ):p\longmapsto \Delta
(ap). $

We will often write, to be concise, $\tau (a,\Delta )\equiv a\Delta $ and $%
\tau ^{+}(a,\Delta )\equiv a^{+}\Delta .$ It is easy to see that $(${\it Diff%
}$_k(P,Q),(\tau ,\tau ^{+}))\equiv \left( {\rm Diff}_k^{}(P,Q),{\rm Diff}%
_k^{+}(P,Q)\right) \doteq {\rm Diff}_k^{(+)}(P,Q)$ is an $A$-bimodule.\ 

\begin{remark}
\label{r1}Since 
\[
\delta _{a_0}(\Delta )\equiv 0\Leftrightarrow \Delta (a_0p)=a_0\Delta \left(
p\right) \text{ ,\quad }\forall \text{ }a_0\in A,\;\forall p\in P 
\]
${\rm Diff}_0(P,Q)$ and $Hom_A(P,Q)$ are identified as $A$-(bi)modules: 
\[
Hom_A(P,Q)\simeq {\rm Diff}_0(P,Q)\simeq {\rm Diff}_0^{+}(P,Q). 
\]
\end{remark}

The obvious inclusion (of sets): 
\[
{\rm Diff}_k(P,Q)\hookrightarrow {\rm Diff}_l(P,Q),\text{ }k\leq l 
\]
induces a monomorphism of $A$-bimodules: 
\[
{\rm Diff}_k^{(+)}(P,Q)\hookrightarrow {\rm Diff}_l^{(+)}(P,Q),\text{ }k\leq
l\text{;} 
\]
the direct limit of the system in $A-{\bf BiMo{\rm D}}$ :

\[
{\rm Diff}_0^{(+)}(P,Q)\hookrightarrow {\rm Diff}_1^{(+)}(P,Q)%
\hookrightarrow ...\hookrightarrow {\rm Diff}_n^{(+)}(P,Q)\hookrightarrow
... 
\]
is denoted by ${\rm Diff}^{(+)}(P,Q)=({\rm Diff}(P,Q),{\rm Diff}^{+}(P,Q))$.

With a given $A$-module $P$ the following three functors are associated

\[
\begin{array}{c}
{\rm Diff}_k:Q\longmapsto {\rm Diff}_k\left( P,Q\right) \\ 
{\rm Diff}_k^{+}:Q\longmapsto {\rm Diff}_k^{+}\left( P,Q\right) \\ 
{\rm Diff}_k^{\left( +\right) }:Q\longmapsto {\rm Diff}_k^{\left( +\right)
}\left( P,Q\right) .
\end{array}
\]
Let us put ${\rm Diff}_k^{\left( +\right) }\left( A,Q\right) \equiv {\rm Diff%
}_k^{\left( +\right) }Q$. Remark \ref{r1} implies ${\rm Diff}_0^{+}={\rm Diff%
}_0={\rm Id}_{A-{\bf MoD}}$.

To simplify notations we will write ${\rm Diff}_{\sigma _1,\ldots ,\sigma
_n}^{+}$ instead of ${\rm Diff}_{\sigma _1}^{+}\circ \cdots \circ {\rm Diff}%
_{\sigma _n}^{+}$.

\begin{definition}
For any $s,t\geq 0$,\label{alucce} define an $A$-module homomorphism 
\[
{\em C}_{s,t}(P):{\rm Diff}_s^{+}({\rm Diff}_t^{+}P)\longrightarrow {\rm Diff%
}_{s+t}^{+}P 
\]
by

\[
{\em C}_{s,t}(P)(\Delta ):a\longmapsto \Delta (a)(1)\text{,\quad }\Delta \in 
{\rm Diff}_s^{+}({\rm Diff}_t^{+}P). 
\]
Then $P\longmapsto {\em C}_{s,t}(P)$ defines a morphism ${\rm Diff}%
_{s,t}^{+}\longrightarrow {\rm Diff}_{s+t\text{ }}^{+}$ of functors called
the {\em composition} or {\em ''gluing'' morphism}.
\end{definition}

Note that ${\rm D}_{\left( k\right) }\left( Q\right) \doteq \left\{ \Delta
\in {\rm Diff}_kQ\mid \Delta \left( 1\right) =0\right\} $ is an $A$%
-submodule of ${\rm Diff}_kQ$ (but not of ${\rm Diff}_k^{+}Q$ !). The
functor ${\rm D}_{\left( k\right) }:Q\mapsto D_{\left( k\right) }\left(
Q\right) $ allows to form the following short exact sequence:

\begin{equation}
0\rightarrow {\rm D}_{\left( k\right) }\stackrel{i_k}{\longrightarrow }{\rm %
Diff}_k\stackrel{p_k}{\longrightarrow }{\rm Id}_{A-{\bf Mod}}\rightarrow 0
\label{du}
\end{equation}
in $\left[ A-{\bf Mod,}A-{\bf Mod}\right] $, where $i_k$ is the canonical
inclusion and $p_k$ is defined by:

\[
p_k\left( Q\right) :{\rm Diff}_kQ\rightarrow Q:\Delta \mapsto \Delta \left(
1\right) \text{\quad , }\Delta \in {\rm Diff}_kQ 
\]
for any $A$-module $Q$. The functor monomorphism ${\rm Id}_{A-{\bf Mod}%
}\equiv {\rm Diff}_0\hookrightarrow {\rm Diff}_k$ splits (1), so that ${\rm %
Diff}_k={\rm D}_{\left( k\right) }\bigoplus {\rm Id}_{A-{\bf Mod}}$. Note
that ${\rm D}_{\left( 1\right) }\left( Q\right) $ is nothing but the $A$%
-module of all $Q$-valued $K$-linear {\it derivations} on $A$, denoted in
the literature usually by ${\rm Der}_{A/K}\left( Q\right) $ (see for example 
\cite{Bou X}).

Let $Q$ be an $A$-module and $P$, $P^{+}$ be the left and right $A$-modules
corresponding to an $A$-bimodule $P^{\left( +\right) }\equiv \left(
P,P^{+}\right) $. Let's denote by ${\rm Diff}_k^{\bullet }\left(
Q,P^{+}\right) $ (resp. ${\rm D}_{\left( k\right) }^{\bullet }\left(
P^{+}\right) $) the $A$-module which coincides with ${\rm Diff}_k\left(
Q,P^{+}\right) $ (resp. ${\rm D}_{\left( k\right) }\left( P^{+}\right) $) as 
$K$-modules and whose $A$-module structures are inherited by that of $P$ 
\footnote{%
These $A$-module structures are well defined due to the fact that $\left(
P,P^{+}\right) \equiv P^{\left( +\right) }$ is a bimodule. Moreover one can
give similar definitions with $P^{+}$ replaced by $P$ i.e. using the
canonical involution $A-{\bf BiMod\longrightarrow }A-{\bf BiMod:}\left(
P,P^{+}\right) \rightarrow \left( P^{+},P\right) $.}:

\begin{eqnarray*}
\text{(mult. by }a\text{ in }{\rm Diff}_k^{\bullet }\left( Q,P^{+}\right) 
\text{)}\qquad \left( a^{\bullet }\Delta \right) \left( q\right) &\doteq
&a\Delta \left( q\right) \\
\text{(mult. by }a\text{ in }D_k^{\bullet }\left( P^{+}\right) \text{)}%
\qquad \left( a^{\bullet }\delta \right) \left( q\right) &\doteq &a\delta
\left( q\right)
\end{eqnarray*}
where both $a\Delta \left( q\right) $ and $a\delta \left( q\right) $ denote
the multiplication by $a$ in $P$. The correspondence

\[
{\rm D}_{(k)}^{(\bullet )}:P^{(+)}\longmapsto (D_{(k)}(P^{+}),\text{ }%
D_{(k)}^{\bullet }(P^{+})) 
\]
defines a functor $A-{\bf BiMod\longrightarrow }A-{\bf BiMod}$ in an obvious
way. If $Q=A$ we write ${\rm Diff}_k^{\bullet }\left( P^{+}\right) $ for $%
{\rm Diff}_k^{\bullet }\left( A,P^{+}\right) $. Obviously, ${\rm D}_{\left(
k\right) }^{\bullet }\left( P^{+}\right) $ is an $A$-submodule of ${\rm Diff}%
_k^{\bullet }(P^{+})$.

Let us introduce now the {\it category of biads}, $*A-{\bf BiMod}$, whose
objects are ordered triples of $A$-modules:

\[
(P,P^{+};Q) 
\]
with $P^{\left( +\right) }\doteq (P,P^{+})$ being an $A$-bimodule and $Q$ an 
$A$-submodule of $P$. The corresponding morphisms are those of underlying $A$%
-bimodules ''respecting'' the selected submodules, i.e.:

\[
f:(P,P^{+})\longrightarrow (\overline{P},\overline{P}^{+})\text{ such that }%
f(Q)\subset \overline{Q}. 
\]

\begin{examples}
\label{e1}(i) If $P$ is an $A$-module and $s\leq k$, then $({\rm Diff}_kP,%
{\rm Diff}_k^{+}P;{\rm D}_{(s)}(P))$ is a biad for each $s\leq k$.

(ii) If $^{}P^{\left( +\right) }\doteq (P,P^{+})$ is an $A$- bimodule, then
the following are biads: 
\[
\begin{array}{c}
({\rm Diff}_k^{\bullet }P^{+},{\rm Diff}_k^{+}P^{+};{\rm D}_{(k)}^{\bullet
}(P^{+})) \\ 
({\rm Diff}_kP,{\rm Diff}_k^{+}P;{\rm D}_{(k)}(P)) \\ 
({\rm Diff}_kP^{+},{\rm Diff}_k^{+}P^{+};{\rm D}_{(k)}(P^{+}))\text{.}
\end{array}
\]
\end{examples}

We associate to a biad $(P^{(+)};Q)$ the following $K$-modules:

\[
\begin{array}{c}
{\rm Diff}_k(Q\subset P^{+})\doteq \left\{ \Delta \in {\rm Diff}_kP^{+}\mid
im(\Delta )\subset Q\right\} \\ 
{\rm D}_{(k)}(Q\subset P^{+})\doteq \left\{ \Delta \in {\rm D}%
_{(k)}(P^{+})\mid im(\Delta )\subset Q\right\} .
\end{array}
\]
The $A$-module structure of ${\rm Diff}_k^{\bullet }(P^{+})$ (resp., of $%
{\rm Diff}_k^{+}(P^{+})$ or of ${\rm D}_{(k)}^{\bullet }(P^{+})$) induces an 
$A$-module structure on ${\rm Diff}_k(Q\subset P^{+})$ (resp., on ${\rm Diff}%
_k^{+}(Q\subset P^{+})$, resp. on ${\rm D}_{(k)}(Q\subset P^{+})$), which is
denoted by ${\rm Diff}_k^{\bullet }(Q\subset P^{+})$ (resp. ${\rm Diff}%
_k^{+}(Q\subset P^{+})$ or ${\rm D}_{(k)}(Q\subset P^{+})$). So, the
following inclusions take place in $A-{\bf Mod}$

\[
\begin{array}{c}
{\rm Diff}_k^{\bullet }(Q\subset P^{+})\subset {\rm Diff}_k^{\bullet }(P^{+})
\\ 
{\rm Diff}_k^{+}(Q\subset P^{+})\subset {\rm Diff}_k^{+}(P^{+}) \\ 
{\rm D}_{(k)}(Q\subset P^{+})\subset {\rm D}_{(k)}^{\bullet }(P^{+})\text{.}
\end{array}
\]
We have then a biad:

\[
({\rm Diff}_k^{\bullet }(Q\subset P^{+}),\text{ }{\rm Diff}_k^{+}(Q\subset
P^{+}),\text{ }{\rm D}_{(k)}(Q\subset P^{+})) 
\]
which is a sub-biad of $({\rm Diff}_k^{\bullet }(P^{+}),$ ${\rm Diff}%
_k^{+}(P^{+}),$ ${\rm D}_{(k)}^{\bullet }(P^{+}))$. The following result is
straightforward\ 

\begin{lemma}
\TeXButton{sissa}{\label{sissa}}\label{l2}\label{l2}\label{l2}\label{l2} If $%
\overline{P}^{(+)}$ is a sub-bimodule of $P^{(+)}$ and $Q\subset \overline{P}
$, then we have:

\[
\begin{array}{c}
{\rm D}_{(k)}(Q\subset \overline{P}^{+})={\rm D}_{(k)}(Q\subset P^{+}) \\ 
{\rm Diff}_k^{+}(Q\subset \overline{P}^{+})={\rm Diff}_k^{+}(Q\subset P^{+})
\\ 
{\rm Diff}_k^{\bullet }(Q\subset \overline{P}^{+})={\rm Diff}_k^{\bullet
}(Q\subset P^{+}).
\end{array}
\]
\end{lemma}

Canonical functors

\[
A-{\bf Mod}\stackrel{i}{\bf \hookrightarrow }A-{\bf BiMod} 
\begin{tabular}{l}
$\stackrel{j}{\bf \hookrightarrow }$ \\ 
$\stackunder{h}{\leftarrow }$%
\end{tabular}
{\bf *}A-{\bf BiMod} 
\]
are defined as follows 
\[
i:P\longmapsto \left( P,P\right) ,j:\left( P,P^{+}\right) \longmapsto \left(
P,P^{+};P\right) , 
\]
\[
h:\left( P,P^{+};Q\right) \longmapsto \left( P,P^{+}\right) \qquad \text{%
(forgetful).} 
\]
Obviously, $i$ is fully faithful while $h$ is left inverse and adjoint to $j$%
. So, $j$ is fully faithful. Moreover there are three obvious forgetful
functors $p_1,p_2,p_3:$ ${\bf *}A-{\bf BiMod\rightarrow }A-{\bf Mod}$ with $%
p_1\left( P,P^{+};Q\right) =P,$ $p_2\left( P,P^{+};Q\right) =P^{+},$ $%
p_3\left( P,P^{+};Q\right) =Q$.

We define now some {\em absolute functors} we need (see Definition \ref
{def
1.18} for their ''relative'' version).

\begin{definition}
\TeXButton{garp}{\label{garp}} For $k\geq 0$, $\underline{{\cal P}}_{(k)}:*A-%
{\bf BiMod}\longrightarrow *A-{\bf BiMod}$ is the functor

\[
\underline{{\cal P}}_{(k)}:(P^{(+)};Q)\mapsto \left( {\rm Diff}_k^{\bullet
}(Q\subset P^{+}),\text{ }{\rm Diff}_k^{+}(Q\subset P^{+});\text{ }{\rm D}%
_{(k)}(Q\subset P^{+})\right) 
\]
\ For $\sigma =\left( \sigma _1,...,\sigma _n\right) \in {\bf N}_{+}^n$,
define \underline{${\cal P}$}$_\sigma \equiv $\underline{${\cal P}$}$%
_{\left( \sigma _1,...,\sigma _n\right) }:*A-{\bf BiMod}\longrightarrow *A-%
{\bf BiMod}$ as the composition:

\[
\underline{{\cal P}}_{\left( \sigma _1,...,\sigma _n\right) }\doteq 
\underline{{\cal P}}_{\left( \sigma _1\right) }\circ \cdots \circ \underline{%
{\cal P}}_{\left( \sigma _n\right) } 
\]
We put 
\begin{eqnarray*}
{\cal P}_{\left( \sigma _1,...,\sigma _n\right) }^{\bullet }=p_1\circ 
\underline{{\cal P}}_{\left( \sigma _1,...,\sigma _n\right) } \\
{\cal P}_{\left( \sigma _1,...,\sigma _n\right) }^{+}=p_2\circ \underline{%
{\cal P}}_{\left( \sigma _1,...,\sigma _n\right) } \\
{\rm D}_{\left( \sigma _1,...,\sigma _n\right) }^{\bullet }=p_3\circ 
\underline{{\cal P}}_{\left( \sigma _1,...,\sigma _n\right) }.
\end{eqnarray*}

Thanks to the full faithfulness of $i$ and $j$ mentioned above, we can
simply write $\underline{{\cal P}}_{\left( \sigma _1,...,\sigma _n\right) }$
for both $\underline{{\cal P}}_{\left( \sigma _1,...,\sigma _n\right) }\circ
j$ and $\underline{{\cal P}}_{\left( \sigma _1,...,\sigma _n\right) }\circ
j\circ i$ by specifying the source category only in case of a possible
confusion.
\end{definition}

\begin{remark}
\TeXButton{rem 1.14}{\label{rem 1.14}} {\it (}{\em Words }{\rm \cite{Vi1}}%
{\it )} If $\tau \in {\bf N}_{+}^k$ and $\sigma \in {\bf N}_{+}^n$, the
functor

\[
\underline{{\cal P}}_\tau ^{}({\rm D}_\sigma ^{\bullet }\subset {\cal P}%
_\sigma ^{+})\equiv \left( {\cal P}_\tau ^{\bullet }({\rm D}_\sigma
^{\bullet }\subset {\cal P}_\sigma ^{+}),{\cal P}_\tau ^{+}({\rm D}_\sigma
^{\bullet }\subset {\cal P}_\sigma ^{+});{\rm D}_\tau ^{\bullet }({\rm D}%
_\sigma ^{\bullet }\subset {\cal P}_\sigma ^{+})\right) 
\]

is obviously identified with $\underline{{\cal P}}_{(\tau _1)}\circ \ldots
\circ \underline{{\cal P}}_{(\tau _k)}\circ \underline{{\cal P}}_\sigma $ ,
i.e.

\[
\underline{{\cal P}}_\tau ({\rm D}_\sigma \subset {\cal P}_\sigma ^{+})=%
\underline{{\cal P}}_{(\tau ,\sigma )} 
\]
where $(\tau ,\sigma )\doteq (\tau _1,\ldots ,\tau _k,\sigma _1,\ldots
,\sigma _n)\in {\bf N}_{+}^{n+k}$. This gives a ''closed'' functorial
meaning to ''words'' composed of functors defined above.
\end{remark}

If $P$ is an $A$-module we write ${\rm D}_{\left( \sigma _1,\ldots ,\sigma
_n\right) }(P)$ for ${\rm D}_{\left( \sigma _1,\ldots ,\sigma _n\right)
}^{\bullet }(j\circ i\left( P\right) )$). It is not difficult to see that 
\[
{\rm D}_{\left( \sigma _1,\ldots ,\sigma _n\right) }(P)={\rm D}_{\left(
\sigma _1\right) }\left( {\rm D}_{\left( \sigma _2,\ldots ,\sigma _n\right)
}(P)\subset {\rm Diff}_{\sigma _2,\ldots ,\sigma _n}^{+}(P)\right) . 
\]
Furthermore,

\[
{\rm D}_{\left( \sigma _1,\ldots ,\sigma _n\right) }(P)\hookrightarrow {\rm D%
}_{\left( \sigma _1\right) }^{\bullet }\left( {\rm Diff}_{\sigma _2,\ldots
,\sigma _n}^{+}(P)\right) \hookrightarrow {\rm Diff}_{\sigma _1}^{\bullet
}\left( {\rm Diff}_{\sigma _2,\ldots ,\sigma _n}^{+}(P)\right) 
\]
are inclusions in $A-{\bf Mod}$ , while the inclusion ${\rm D}_{\left(
\sigma _1,\ldots ,\sigma _n\right) }(P)\hookrightarrow {\rm Diff}_{\sigma
_1,\ldots ,\sigma _n}^{+}(P)$ is a DO of order $\leq \sigma _1$. So, the
functors ${\rm D}_\sigma $'s can be defined also inductively:

\begin{definition}
Let $\sigma =\left( \sigma _1,\sigma _2,\ldots \sigma _n,\ldots \right) \in 
{\bf N}_{+}^\infty $ and $\sigma \left( n\right) =\left( \sigma _1,\ldots
,\sigma _n\right) $. Then functors ${\rm D}_{\sigma \left( n\right) }:A-{\bf %
Mod\rightarrow }A-{\bf Mod}$ are defined by induction: 
\[
\begin{array}{c}
{\rm D}_{\sigma \left( 1\right) }\doteq {\rm D}_{\left( \sigma _1\right) }
\\ 
{\rm D}_{\sigma \left( n\right) }:P\longmapsto {\rm D}_{\left( \sigma
_1\right) }\left( {\rm D}_{\left( \sigma _2,\ldots ,\sigma _n\right) }\left(
P\right) \subset {\rm Diff}_{\sigma _2,\ldots ,\sigma _n}^{+}\left( P\right)
\right) .
\end{array}
\]
\end{definition}

If $\sigma =(1,..1,1,...)$ we also write ${\rm D}_n$ for ${\rm D}_{\sigma
\left( n\right) }$.\ 

For any $\sigma \in {\bf N}_{+}^\infty $ and $n\in {\bf N}_{+}$, the
sequence in $\left[ A-{\bf Mod},A-{\bf Mod}\right] $:

\begin{equation}
0\rightarrow {\rm D}_{\sigma \left( n\right) }\stackrel{{\rm I}_{\sigma
\left( n\right) }}{\longrightarrow }{\rm D}_{\sigma \left( n-1\right)
}^{\bullet }\circ {\rm Diff}_{\sigma _n}^{\left( +\right) }\stackrel{\pi
_{\sigma \left( n\right) }}{\longrightarrow }{\rm D}_{\left( \sigma
_1,\ldots ,\sigma _{n-2},\sigma _{n-1}+\sigma _n\right) }  \label{(2)}
\end{equation}
where ${\rm I}_{\sigma \left( n\right) }$ is the natural inclusion and $\pi
_{\sigma \left( n\right) \text{ }}$ is the composition (see Definition \ref
{alucce})

\[
{\rm D}_{\sigma \left( n-1\right) }^{\bullet }\circ {\rm Diff}_{\sigma
_n}^{\left( +\right) }\hookrightarrow {\rm D}_{\sigma \left( n-2\right)
}^{\bullet }\circ {\rm Diff}_{\sigma _{n-1}}^{\left( +\right) }\circ {\rm %
Diff}_{\sigma _n}^{+}\stackrel{{\rm D}_{\sigma \left( n-2\right) }^{\bullet
}\left( {\em C}_{\sigma _{n-1},\sigma _n}\right) }{\longrightarrow }{\rm D}%
_{\sigma \left( n-2\right) }^{\bullet }\circ {\rm Diff}_{\sigma
_{n-1}+\sigma _n}^{\left( +\right) }\text{ ,} 
\]
is exact.

\begin{remark}
\label{new}Let $\sigma =\left( \sigma _1,\sigma _2,\ldots \sigma _n\right)
\in {\bf N}_{+}^n$. We have a canonical split exact short exact sequence in $%
\left[ A-{\bf Mod},A-{\bf Mod}\right] $%
\[
0\rightarrow {\rm D}_\sigma \stackrel{l_\sigma }{\stackunder{\rho _\sigma }{%
\rightleftarrows }}{\cal P}_\sigma ^{\bullet }\stackrel{\psi _\sigma }{%
\longrightarrow }{\rm D}_{\left( \sigma _2,\ldots ,\sigma _n\right)
}\rightarrow 0 
\]
where $l_\sigma $ is the canonical inclusion, $\psi _\sigma $ is given by ($%
P $ is an $A$-module) 
\[
\psi _\sigma (P)(\Delta )=\Delta \left( 1\right) 
\]
and $\rho _\sigma $ by

\[
\rho _\sigma \left( P\right) \left( \Delta \right) =\Delta -\Delta \left(
1\right) . 
\]
Hence ${\cal P}_\sigma ^{\bullet }\simeq {\rm D}_{\left( \sigma _2,\ldots
,\sigma _n\right) }\oplus {\rm D}_\sigma $.
\end{remark}

\smallskip\ \ 

Now we define the {\em relative} (i.e. relative to an arbitrary $A$-module $%
P $) {\em functors}.

\begin{definition}
\TeXButton{def 1.18}{\label{def 1.18}}If $k\geq 0$ and $P$ is an $A$-module,
define 
\[
\underline{{\cal P}}_{\left( k\right) }\left[ P\right] :*A-{\bf %
BiMod\longrightarrow }A-{\bf BiMod} 
\]
to be the functor

\[
\left( Q,Q^{+};S\right) \longmapsto \left( {\rm Diff}_k^{\bullet }\left(
P,S\subset Q^{+}\right) ,{\rm Diff}_k^{+}\left( P,S\subset Q^{+}\right)
\right) . 
\]
If $\sigma \left( n\right) =\left( \sigma _1,\ldots ,\sigma _n\right) \in 
{\bf N}_{+}^n$, $n>1$, define the functor 
\[
\underline{{\cal P}}_{\sigma \left( n\right) }\left[ P\right] :*A-{\bf %
BiMod\longrightarrow }\left( A,A\right) -{\bf BiMod} 
\]
as the composition 
\[
\underline{{\cal P}}_{\sigma \left( n\right) }\left[ P\right] \doteq 
\underline{{\cal P}}_{\left( \sigma _1\right) }\left[ P\right] \circ 
\underline{{\cal P}}_{\left( \sigma _2,...,\sigma _n\right) }. 
\]
As in the absolute case (Def. \ref{garp}), we set 
\begin{eqnarray*}
\underline{{\cal P}}_{\sigma \left( n\right) }^{\bullet }\left[ P\right]
\doteq p_1\circ \underline{{\cal P}}_{\sigma \left( n\right) }\left[ P\right]
\\
\underline{{\cal P}}_{\sigma \left( n\right) }^{+}\left[ P\right] \doteq
p_2\circ \underline{{\cal P}}_{\sigma \left( n\right) }\left[ P\right]
\end{eqnarray*}
and still denote by $\underline{{\cal P}}_{\sigma \left( n\right) }\left[
P\right] $ both $\underline{{\cal P}}_{\sigma \left( n\right) }\left[
P\right] \circ j$ and $\underline{{\cal P}}_{\sigma \left( n\right) }\left[
P\right] \circ j\circ i$.
\end{definition}

By Lemma \ref{sissa}, we have $\left( {\cal P}_{\sigma \left( n\right)
}^{\bullet },{\cal P}_{\sigma \left( n\right) }^{+}\right) \equiv \left( 
{\cal P}_\sigma ^{\bullet }\left[ A\right] ,{\cal P}_\sigma ^{+}\left[
A\right] \right) $. Moreover, $\underline{{\cal P}}_{\sigma \left( n\right)
}\left[ P\right] $ is (controvariantly) functorial in $P$ and we have

\begin{lemma}
\TeXButton{split}{\label{split}}If $0\rightarrow P_1\stackrel{f}{%
\longrightarrow }P_2\stackrel{g}{\longrightarrow }P_3\rightarrow 0$ is exact
(resp. split exact) in $A-{\bf Mod}$, then 
\begin{eqnarray*}
0\rightarrow \underline{{\cal P}}_{\sigma \left( n\right) }^{\bullet }\left[
P_3\right] \stackrel{g^{\vee }}{\longrightarrow }\underline{{\cal P}}%
_{\sigma \left( n\right) }^{\bullet }\left[ P_2\right] \stackrel{f^{\vee }}{%
\longrightarrow }\underline{{\cal P}}_{\sigma \left( n\right) }^{\bullet
}\left[ P_1\right] \\
0\rightarrow \underline{{\cal P}}_{\sigma \left( n\right) }^{+}\left[
P_3\right] \stackrel{g^{\vee }}{\longrightarrow }\underline{{\cal P}}%
_{\sigma \left( n\right) }^{+}\left[ P_2\right] \stackrel{f^{\vee }}{%
\longrightarrow }\underline{{\cal P}}_{\sigma \left( n\right) }^{+}\left[
P_1\right]
\end{eqnarray*}
are exact (resp.

\begin{eqnarray*}
0\rightarrow \underline{{\cal P}}_{\sigma \left( n\right) }^{\bullet }\left[
P_3\right] \stackrel{g^{\vee }}{\longrightarrow }\underline{{\cal P}}%
_{\sigma \left( n\right) }^{\bullet }\left[ P_2\right] \stackrel{f^{\vee }}{%
\longrightarrow }\underline{{\cal P}}_{\sigma \left( n\right) }^{\bullet
}\left[ P_1\right] \rightarrow 0 \\
0\rightarrow \underline{{\cal P}}_{\sigma \left( n\right) }^{+}\left[
P_3\right] \stackrel{g^{\vee }}{\longrightarrow }\underline{{\cal P}}%
_{\sigma \left( n\right) }^{+}\left[ P_2\right] \stackrel{f^{\vee }}{%
\longrightarrow }\underline{{\cal P}}_{\sigma \left( n\right) }^{+}\left[
P_1\right] \rightarrow 0
\end{eqnarray*}
are exact).
\end{lemma}

\TeXButton{Proof}{\proof} Straightforward, by induction on $n$. 
\TeXButton{End Proof}{\endproof}

If $P$ is an $A$-module, there are exact sequences in $[A-{\bf Mod},A-{\bf %
Mod}]$

\begin{equation}
0\rightarrow {\cal P}_{\sigma \left( n\right) }^{\bullet }\left[ P\right]
\hookrightarrow \stackunder{}{\cal P}_{\sigma \left( n-1\right) }^{\bullet
}\left[ P\right] \circ {\rm Diff}_{\sigma _n}^{\left( +\right) }\stackrel{%
q_{\sigma \left( n\right) }}{\longrightarrow }{\cal P}_{\left( \sigma
_1,\ldots \sigma _{n-2},\sigma _{n-1}+\sigma _n\right) }^{\bullet }\left[
P\right]  \label{relHol}
\end{equation}
where the monomorphism is the natural inclusion while $q_{\sigma (n)}$ is
induced by the ''gluing'' morphism with respect to the pair of indexes $%
\left( \sigma _{n-1},\sigma _n\right) $, i.e.:

\begin{eqnarray*}
\stackunder{}{\cal P}_{\left( \sigma _1,\ldots ,\sigma _{n-1}\right)
}^{\bullet }\left[ P\right] \circ {\rm Diff}_{\sigma _n}^{\left( +\right)
}\left( Q\right) &\ni &\Delta \longmapsto q_{\sigma (n)}(\Delta )=\overline{%
\Delta }\in {\cal P}_{\left( \sigma _1,\ldots \sigma _{n-2},\sigma
_{n-1}+\sigma _n\right) }^{\bullet }\left[ P\right] \left( Q\right) \\
\left( \ldots \left( \left( \overline{\Delta }\left( p\right) \right) \left(
a_1\right) \right) \cdots \right) \left( a_{n-2}\right) &\doteq &\left(
\left( \ldots \left( \left( \Delta \left( p\right) \right) \left( a_1\right)
\right) \cdots \right) \left( a_{n-2}\right) \right) \left( 1\right) ,
\end{eqnarray*}
where $p\in P$ and $a_1,...,a_{n-2}\in A$. We have analogous exact sequences
in $[A-{\bf Mod},A-{\bf Mod}]$:

\begin{equation}
0\rightarrow {\cal P}_{\sigma \left( n\right) }^{\bullet }\left[ P\right]
\hookrightarrow {\rm Diff}_{\sigma _1}^{\bullet }\left( P,\cdot \right)
\circ {\cal P}_{\left( \sigma _2,\ldots ,\sigma n\right) }^{+}\left[
P\right] \stackrel{g_{\sigma (n)}}{\longrightarrow }{\cal P}_{\left( \sigma
_1+\sigma _2,\sigma _3,\ldots ,\sigma _n\right) }^{\bullet }\left[ P\right]
\label{relantiHol}
\end{equation}
where $g_{\sigma (n)}:\Delta \mapsto \widehat{\Delta }$ with $\widehat{%
\Delta }\left( p\right) \doteq \Delta \left( p\right) \left( 1\right) $, $%
p\in P$ (i.e. we ''glue'' with respect to the first two indexes); the upper
boldface dot in ${\rm Diff}_{\sigma _1}^{\bullet }\left( P,\cdot \right) $
denotes the $A$-module structure induced by ${\cal P}_{\left( \sigma
_2,\ldots ,\sigma _n\right) }^{\bullet }\left[ P\right] $.

The following definition will allow us to be concise in the next Section:

\begin{definition}
For any $n>0$ and any $\sigma \in {\bf N}_{+}^n$, the functors (in $[A-{\bf %
Mod},A-{\bf Mod}]$) ${\cal P}_\sigma ^{\bullet }$, $D_\sigma ^{\bullet }$
are called the {\em relevant absolute functors} while, if $P$ is an $A$%
-module, the functors ${\cal P}_\sigma ^{\bullet }\left[ P\right] $ , are
called the {\em relevant functors relative} to the{\it \ }$A$-module $P$.
\end{definition}

\section{Absolute and relative representative objects}

In this Section we consider (strict) representative objects of the functors
introduced in the previous Section. We obtain, as particular cases, the
standard modules of K\"ahler differential forms of Algebraic Geometry and
the de Rham forms of Differential Geometry. We emphasize that in our
approach all these (and not only those of degree one) are obtained as
representative objects of suitable functors. One of the major advantages of
this approach is to allow natural generalizations.

Let ${\frak D}$ be a full subcategory of $A-{\bf Mod}$. We denote by $A-{\bf %
BiMod}_{{\frak D}}$ the subcategory of $A-{\bf BiMod}$ whose objects are
couples of objects of ${\frak D}$ and by $*A-{\bf BiMod}_{{\frak D}}$ the
subcategory of $*A-{\bf BiMod}$ consisting of triples whose elements are
objects of ${\frak D}$ (Section 1).

\begin{definition}
\TeXButton{marker def 2.1}{\label{def 2.1}}A full {\em abelian} subcategory $%
{\frak D}$ of $A-{\bf Mod}$ is said to be {\em differentially closed} if the
following properties are satisfied:

{\it (a)} each functor defined in the previous Section, when restricted to $%
{\frak D}$ (resp. $A-{\bf BiMod}_{{\frak D}}$, resp. $*A-{\bf BiMod}_{{\frak %
D}}$) has values in ${\frak D}$ (resp. $A-{\bf BiMod}_{{\frak D}}$, resp. $%
*A-{\bf BiMod}_{{\frak D}}$);

{\it (b)} if $T:A-{\bf Mod\rightarrow }$ $A-{\bf Mod}$ is a {\em relevant
absolute} functor or a {\em relevant relative} functor, {\em relative to an
object of} ${\frak D}$, then $T_{\mid {\frak D}}:{\frak D}\longrightarrow 
{\frak D}$ is strictly representable in ${\frak D}$;

{\it (c)} $A\in Ob\left( {\frak D}\right) $;

{\it (d)} ${\frak D}$ is closed under tensor product (over $A$);

(e) ${\frak D}$ is closed under taking subobjects (i.e. if $P\subseteq Q$ in 
$A-{\bf Mod}$ and $Q$ is in ${\frak D}$ then $P$ is in ${\frak D}$).
\end{definition}

Condition {\it (a)} is needed to have an ambient category which is
''closed'' with respect to functorial differential calculus; as it will be
clear in the following, since among the functors of Section 1 there are also
compositions of relevant ''elementary'' ones, we would like that
representative objects of these nonelementary functors\footnote{%
Relevant or not.} (for example ${\rm D}_{(s)}^{\bullet }\circ {\rm Diff}%
_t^{\left( +\right) }$), if existing, {\it could be expressed in terms} of
representative objects of the relevant ''elementary'' ones (${\rm D}%
_{(s)}^{} $ and ${\rm Diff}_t$ in the example). Condition ({\it d}) makes it
possible.

${\frak D}$ being abelian and satisfying ({\it b}), exactness of sequences
of strictly representable functors yields exactness of the ''dual''
sequences of representative objects in ${\frak D}$. Condition (e) is related
to the existence of canonical generators for some representative objects and
will become clear in the sequel. Note also that (e) implies that if $f$ is a
morphism in ${\frak D}$, $im\left( f\right) $ (resp. $\ker \left( f\right) $%
) is the same when considered in ${\frak D}$ or in $A-{\bf Mod}$.

Let us recall some elementary facts about bimodules, mainly to fix our
notations.

If $P^{(+)}=(P,P^{+})$ is an $A$-bimodule and $a\in A$, we write $a$ for the
multiplication in $P$ and $a^{+}$ for the multiplication in $P^{+}$. If $Q$
is an $A$-module we denote by:

({\bf I}) $P^{+}\otimes _A^{\bullet }Q$ the $A$-module obtained by the
abelian group $P^{+}\otimes _AQ$ with multiplication by elements of $A$
defined as

\[
a^{\bullet }(p\otimes q)\doteq (ap)\otimes q\text{, \quad }a\in A,\text{ }%
p\in P^{+},\text{ }q\in Q 
\]
(note that $a^{\bullet }(p\otimes q)\neq p\otimes aq$). Then $P^{+}\otimes
_A^{(\bullet )}Q$ $\doteq \left( P^{+}\otimes _A^{\bullet }Q,P^{+}\otimes
_AQ\right) $ is an $A$-bimodule;

({\bf II}) $Hom_A^{\bullet }(Q,P^{+})$ the $A$-module obtained by the
abelian group $Hom_A(Q,P^{+})$ with multiplication by elements of $A$
defined as :

\[
\left[ a^{\bullet }f\right] (p)\doteq a_0\cdot (f(p))\text{, }\quad a\in A,%
\text{ }p\in P,\text{ }f\in Hom_A(Q,P^{+}). 
\]
Denote by $Hom_A^{(\bullet )}(Q,P^{+})\doteq \left( Hom_A^{\bullet
}(Q,P^{+}),Hom_A(Q,P^{+})\right) $ the corresponding $A$-bimodule;

({\bf III}) $Hom_A^{+}(P,Q)$ the $A$-module obtained by the abelian group $%
Hom_A(P,Q)$ with multiplication by elements of $A$: 
\[
\left[ a^{+}f\right] (p)\doteq f(a^{+}p). 
\]
Then 
\[
Hom_A^{(+)}(P,Q)\doteq \left( Hom_A(P,Q),Hom_A^{+}(P,Q)\right) 
\]
is an $A$-bimodule.

In the same way we can define the $A$-modules $P\otimes _A^{+}Q$ , $%
Hom_A^{+}(P,Q)$ and $Hom_A^{\bullet }(P^{+},Q)$.

\begin{example}
\TeXButton{marker es 2.2}{\label{es 2.2}}\label{c6}If $P$ and $Q$ are $A$%
-modules, we have an isomorphism in \ $\left[ A-{\bf Mod,}A-{\bf Mod}\right] 
$

\[
{\rm Diff}_s(\cdot ,Q)\simeq Hom_A^{\bullet }\left( \cdot ,{\rm Diff}%
_s^{+}Q\right) \text{.} 
\]
\end{example}

It is not diificult to prove the following\ 

\begin{lemma}
\label{l3}\label{l3}\TeXButton{marker lem 2.3}{\label{lem 2.3}}Let $R$ and $%
P $ be $A$-modules and $(Q,Q^{+})$ an $A$-bimodule. Then we have a canonical
isomorphism in $A-{\bf BiMod}$: 
\[
Hom_A^{\left( \bullet \right) }\left( R,Hom_A^{+}(Q,P)\right) \stackrel{%
\varphi }{\simeq }Hom_A^{(+)}(Q^{+}\otimes _A^{\bullet }R,P) 
\]
\end{lemma}

\begin{proposition}
\label{jet}Let $P$ be an $A$-module and $k\in {\bf N}_{+}$. Then ${\rm Diff}%
_k\left( P,\cdot \right) :A-{\bf Mod\rightarrow }$ $A-{\bf Mod}$ is strictly
representable by the so-called $k${\em -jet module} ${\bf J}^k\left(
P\right) $.
\end{proposition}

\TeXButton{Proof}{\proof} See \cite{KLV}, p. 12. \TeXButton{End Proof}
{\endproof}\ 

In other words, there exists a universal DO $j_k(P):P\longrightarrow {\bf J}%
^k(P)$, of order $\leq k$ (often denoted simply by $j_k$), such that for
each DO $\Delta :P\rightarrow Q$ of order $\leq k$, there is a unique $A$%
-homomorphism $f^\Delta $ and a commutative diagram

\[
\begin{array}{ccccc}
P & \stackrel{j_k(P)}{\longrightarrow } & {\bf J}^k(P) &  &  \\ 
& \stackunder{\Delta }{\searrow } & \downarrow & f^\Delta &  \\ 
&  & Q &  & 
\end{array}
\text{.} 
\]
The $A$-module ${\bf J}^k(P)$ is generated by $\left\{ j_k(p)\mid p\in
P\right\} $. Moreover, ${\bf J}^k(P)$ has a bimodule structure ${\bf J}%
_{\left( +\right) }^k\left( P\right) =\left( {\bf J}^k\left( P\right) ,{\bf J%
}_{+}^k\left( P\right) \right) $ which can be described in the following,
purely functorial, way. Suppose ${\frak D}\subseteq A-{\bf Mod}$ is a
subcategory such that $\forall P\in Ob\left( {\frak D}\right) $ the functor $%
{\rm Diff}_k\left( P,\cdot \right) $ when restricted to ${\frak D}$ has
values in ${\frak D}$ and is strictly representable in ${\frak D}$ (e.g. $%
{\frak D}=A-{\bf Mod}$ by proposition \ref{lem 2.3}). Let ${\bf J}_{{\frak D}%
}^k\left( P\right) $ be the corresponding representative object i.e. $%
Hom_A\left( {\bf J}_{{\frak D}}^k\left( P\right) ,\cdot \right) \simeq {\rm %
Diff}_k\left( P,\cdot \right) $ in $\left[ {\frak D},{\frak D}\right] $. If $%
j_k^{{\frak D}}:P\rightarrow {\bf J}_{{\frak D}}^k\left( P\right) $
corresponds to the identity morphism of ${\bf J}_{{\frak D}}^k\left(
P\right) $, we can define, for each $a\in A$, the DO 
\[
a^{+}:P\rightarrow {\bf J}_{{\frak D}}^k\left( P\right) :p\longmapsto j_k^{%
{\frak D}}\left( ap\right) \text{.\quad } 
\]
The corresponding $A$-endomorphism of ${\bf J}_{{\frak D}}^k$ is still
denoted by $a^{+}$ and gives the required second $A$-module structure ${\bf J%
}_{{\frak D}+}^k$ on the abelian group ${\bf J}_{{\frak D}}^k$.

Using the bimodule ${\bf J}_{(+)}^k$ and Proposition \ref{jet}, we get an
isomorphism in

$\left[ A-{\bf Mod,}A-{\bf Mod}\right] $

\begin{equation}
{\rm Diff}_k^{+}\simeq Hom_A^{+}({\bf J}^k,\cdot )\text{.}  \label{plus}
\end{equation}

If $P,Q$ are $A$-modules, then:

\[
Hom_A({\bf J}_{+}^k\otimes _A^{\bullet }P,Q)\simeq Hom_A^{\bullet
}(P,Hom_A^{+}({\bf J}^k,Q)) 
\]
by Lemma \ref{l3}; therefore

\[
Hom_A^{\bullet }(P,Hom_A^{+}({\bf J}^s,Q))\simeq Hom_A^{\bullet }(P,{\rm Diff%
}_s^{+}Q)\text{. } 
\]
and, by (\ref{plus}) and Example \ref{es 2.2}, we finally get:

\[
Hom_A^{\bullet }(P,{\rm Diff}_s^{+}Q)\simeq {\rm Diff}_s(P,Q)\text{. } 
\]
Since representative objects of the same functor are canonically isomorphic,
we have proved:

\begin{lemma}
\TeXButton{lem 2.5}{\label{lem 2.5}}\label{ozu}There are canonical
isomorphisms in $\left[ A-{\bf Mod,}A-{\bf Mod}\right] $: 
\[
\begin{array}{c}
{\bf J}^k(\cdot )\simeq {\bf J}_{+}^k\otimes _A^{\bullet }\left( \cdot
\right) \\ 
{\bf J}_{+}^k(\cdot )\simeq {\bf J}_{+}^k\otimes _A\left( \cdot \right) 
\text{ .}
\end{array}
\]
\end{lemma}

We are now able to prove a basic result

\begin{lemma}
\TeXButton{lem 2.6}{\label{lem 2.6}}({\bf a}) If $\tau \in {\bf N}%
_{+}^\infty $, $n>0$, $t\geq 0$ and ${\rm D}_{\tau \left( n\right) }$ is
strictly representable in $A-{\bf Mod}$ by ${\bf \Lambda }^{\tau \left(
n\right) }$, then ${\rm D}_{\tau \left( n\right) }^{\bullet }\circ {\rm Diff}%
_t^{\left( +\right) }$ is strictly representable in $A-{\bf Mod}$ by ${\bf J}%
^t\left( {\bf \Lambda }^{\tau \left( n\right) }\right) $;

\label{l5}({\bf b}) If $s,t\geq 0$, then 
\[
{\cal P}_{(s)}^{\bullet }\circ {\rm Diff}_t^{\left( +\right) }\equiv {\rm %
Diff}_s^{\bullet }\circ {\rm Diff}_t^{\left( +\right) }:A-{\bf %
Mod\rightarrow }A-{\bf Mod} 
\]
is strictly representable by ${\bf J}^t({\bf J}^s).$

({\bf c}) If $s,t\geq 0$ and $P$ is an $A$-module, then 
\[
{\cal P}_{(s)}^{\bullet }\left[ P\right] \circ {\rm Diff}_t^{\left( +\right)
}\equiv {\rm Diff}_s^{\bullet }\left( P,{\rm Diff}_t^{\left( +\right)
}\left( \cdot \right) \right) :A-{\bf Mod\rightarrow }A-{\bf Mod} 
\]
is strictly representable by ${\bf J}^t({\bf J}^s\left( P\right) )$.

({\bf d}) If ${\cal P}_{\sigma (n)}^{\bullet }$ is strictly representable by 
${\rm Hol}^{\sigma \left( n\right) }$, then 
\[
{\cal P}_{\sigma (n)}^{\bullet }\circ {\rm Diff}_k^{(+)}:A-{\bf %
Mod\rightarrow }A-{\bf Mod} 
\]
is strictly representable by ${\bf J}^k\left( {\rm Hol}^{\sigma \left(
n\right) }\right) $;

({\bf e}) If $P$ is an $A$-module and $\stackunder{}{\cal P}_{\sigma
(n)}^{\bullet }\left[ P\right] $ is strictly representable by ${\rm Hol}%
^{\sigma \left( n\right) }\left[ P\right] $, then 
\[
{\cal P}_{\sigma (n)}^{\bullet }\left[ P\right] \circ {\rm Diff}_k^{(+)}:A-%
{\bf Mod\rightarrow }A-{\bf Mod} 
\]
is strictly representable by ${\bf J}^k\left( {\rm Hol}^{\sigma \left(
n\right) }\left[ P\right] \right) $;
\end{lemma}

\TeXButton{Proof}{\proof} The proofs are very similar. We prove only \label
{l5}({\bf b}) and ({\bf d}).

({\bf b}) 
\[
{\rm Diff}_s^{\bullet }({\rm Diff}_t^{+}P)\simeq Hom_A^{\bullet }({\bf J}^s,%
{\rm Diff}_t^{+}P)\simeq Hom_A^{\bullet }({\bf J}^s,Hom_A^{+}({\bf J}^t,P))%
\text{ ;} 
\]
by (\ref{plus}) this is isomorphic to $Hom_A({\bf J}_{+}^t\otimes ^{\bullet }%
{\bf J}^s,P)$ and, finally by Proposition \ref{ozu}, to $Hom_A({\bf J}^t(%
{\bf J}^s),P)$.

({\bf d}) 
\begin{eqnarray*}
{\cal P}_{\sigma (n)}^{\bullet }\left( {\rm Diff}_k^{(+)}\left( P\right)
\right) &\simeq &Hom_A^{\bullet }({\rm Hol}^{\sigma \left( n\right) },{\rm %
Diff}_t^{+}P)\simeq Hom_A({\rm Hol}^{\sigma \left( n\right) },Hom_A^{+}({\bf %
J}^k,P))\simeq \\
\ &\simeq &Hom_A({\bf J}_{+}^k\otimes ^{\bullet }{\rm Hol}^{\sigma \left(
n\right) },P)\simeq Hom_A\left( {\bf J}^k\left( {\rm Hol}^{\sigma \left(
n\right) }\right) ,P\right) \text{.}
\end{eqnarray*}

\qquad \TeXButton{End Proof}{\endproof}

Note that if ${\rm D}_\sigma $ is representable for any $\sigma \in {\bf N}%
_{+}^n$, $n>0$ then ${\cal P}_\sigma ^{\bullet }$ is representable for any $%
\sigma \in {\bf N}_{+}^n$, $n>0$ by Remark \ref{new}.

\begin{remark}
Note that, for example, ${\rm Diff}_k^{\bullet }\circ {\rm Diff}_l^{\left(
+\right) }\circ {\rm Diff}_m^{+}$ is representable but not strictly
representable in $A-{\bf Mod}$:

\begin{eqnarray*}
{\rm Diff}_k^{\bullet }({\rm Diff}_l^{+}\left( {\rm Diff}_m^{+}P\right)
)\simeq Hom_A({\bf J}^l({\bf J}^k),{\rm Diff}_m^{+}P)\simeq \\
\ \simeq Hom_A\left( {\bf J}^l({\bf J}^k),Hom_A^{+}\left( {\bf J}^m,P\right)
\right) \simeq Hom_A^{+}\left( {\bf J}_{+}^m\otimes _A^{\bullet }{\bf J}^l(%
{\bf J}^k),P\right) \simeq \\
\ \simeq Hom_A^{+}\left( {\bf J}^m\left( {\bf J}^l({\bf J}^k)\right)
,P\right) \text{.}
\end{eqnarray*}
\end{remark}

\smallskip\ 

We conclude this preliminaries with the following elementary result

\begin{lemma}
{\rm (''Third-representable'' lemma)} \label{terzo}Let $0\rightarrow T_1%
\stackrel{i}{\rightarrow }T_2\stackrel{\varphi }{\rightarrow }T_3$ be an
exact sequence in $\left[ A-{\bf Mod,}A-{\bf Mod}\right] $, with $T_2$ and $%
T_3$ strictly representable by $\tau _2$ and $\tau _3$, respectively; then $%
T_1$ is also strictly representable by: 
\[
\frac{\tau _2}{\varphi ^{\vee }(\tau _3)}\text{,} 
\]
where $\varphi ^{\vee }:\tau _3\rightarrow \tau _2$ is the
dual-representative of $\varphi $.
\end{lemma}

\TeXButton{Proof}{\proof} If $P$ is an $A$-module, the morphism

\[
\begin{array}{c}
\chi _P:T_1(P)\rightarrow Hom_A(\frac{\tau _2}{\varphi ^{\vee }(\tau _3)},P)
\\ 
q\longmapsto \chi _p(q):\left[ t_2\right] _{{\bf mod}\varphi ^{\vee }(\tau
_3)}\longmapsto \stackrel{\wedge }{q}(t_2)
\end{array}
\]
where $\stackrel{\wedge }{q}\doteq i(P)(q)$, is well defined (since $%
\stackrel{\wedge }{q}\circ \varphi ^{\vee }=\varphi (P)(\stackrel{\wedge }{q}%
)=0$) and is an isomorphism, natural in $P$.\TeXButton{End Proof}{\endproof}

The next Theorem, collecting some of the results above, asserts that $A-{\bf %
Mod}${\it \ is itself differentially closed }(see Definition{\it \ }\ref
{def
2.1}).

\begin{theorem}
\TeXButton{thm 2.7}{\label{thm 2.7}}Let $P$ be an $A$-module, $\sigma \in $ $%
{\bf N}_{+}^\infty $ and $k\in {\bf N}$. Then:

({\bf i}) ${\rm Diff}_k\left( P,\cdot \right) $ is strictly representable in 
$A-{\bf Mod}$ by the $k${\em -jet module} ${\bf J}^k\left( P\right) $;

({\bf ii}) for each $n>0$, the functor ${\rm D}_{\sigma \left( n\right) }$
is strictly representable in $A-{\bf Mod}$ by the so-called {\em higher de
Rham forms}{\it ' module of type} $\sigma \left( n\right) $, ${\bf \Lambda }%
^{\sigma \left( n\right) }$;

({\bf iii}) ${\cal P}_{\sigma (n)}^{\bullet }$ and ${\cal P}_{\sigma
(n)}^{\bullet }\left[ P\right] $ $:A-{\bf Mod\rightarrow }$ $A-{\bf Mod}$
are strictly representable by the so-called {\em absolute holonomy module}
of type $\sigma \left( n\right) $, ${\rm Hol}^{\sigma \left( n\right) }$ and 
{\em relative holonomy module} of type $\sigma \left( n\right) $, ${\rm Hol}%
^{\sigma \left( n\right) }\left[ P\right] $.
\end{theorem}

\TeXButton{Proof}{\proof} ({\bf i}) is Proposition \ref{jet}.

({\bf ii})\ The strict representability of ${\rm D}_{\sigma \left( n\right)
} $ in $A-{\bf Mod}$ may be proved by induction on $n$. The case $n=1$
follows from the exact sequence (\ref{du}), the ''third-representable''
lemma and ({\bf i}). Now, suppose we have proved strict representability of $%
{\rm D}_{\tau (k)}$ for each $\tau \in {\bf N}_{+}^\infty $ and each $k\leq
n-1$. From the exact sequence in $\left[ A-{\bf Mod,}A-{\bf Mod}\right] $

\begin{equation}
0\rightarrow {\rm D}_{\sigma \left( n\right) }\hookrightarrow {\rm D}%
_{\sigma \left( n-1\right) }^{\bullet }\circ {\rm Diff}_{\sigma
_n}^{+}\longrightarrow {\rm D}_{(\sigma \left( n-2\right) ,\sigma
_{n-1}+\sigma _n)}\text{,}  \label{ghazzali}
\end{equation}
the last morphism being

\[
\Delta \longmapsto \widehat{\Delta } 
\]

\[
\left( \left( ...\left( \left( \widehat{\Delta }(a_1)\right) \left(
a_2\right) \right) ...\right) \left( a_{n-2}\right) \right) \left(
a_{n-1}\right) \doteq \left( \left( \left( ...\left( \left( \Delta
(a_1)\right) \left( a_2\right) \right) ...\right) \left( a_{n-2}\right)
\right) \left( a_{n-1}\right) \right) \left( 1\right) 
\]
(i.e. we use the ''gluing'' morphism of definition \ref{alucce} with respect
to the last two indexes), Lemma \ref{lem 2.6} ({\bf a}) and the
''third-representable'' lemma, we obtain strict representability for ${\rm D}%
_{\sigma \left( n\right) }$.

({\bf iii}) the case of ${\cal P}_{\sigma (n)}^{\bullet }$ follows, as for (%
{\bf ii}), by induction via Lemma \ref{lem 2.6} , Lemma \ref{terzo} and by
any of the following two exact sequences in $\left[ A-{\bf Mod,}A-{\bf Mod}%
\right] $:

\begin{equation}
0\rightarrow {\cal P}_{\sigma (n)}^{\bullet }\hookrightarrow {\cal P}%
_{\sigma (n-1)}^{\bullet }\circ {\rm Diff}_{\sigma _n}^{(+)}\stackrel{}{%
\longrightarrow }{\cal P}_{\left( \sigma (n-2),\sigma _{n-1}+\sigma
_n\right) }^{\bullet }  \label{Hol}
\end{equation}

\begin{equation}
0\rightarrow {\cal P}_{\sigma (n)}^{\bullet }\hookrightarrow {\rm Diff}%
_{\sigma _1}^{\bullet }\circ {\cal P}_{\left( \sigma _2,...,\sigma _n\right)
}^{+}\stackrel{{\em C}_{\sigma _1,\sigma _2}\left( {\rm Diff}_{\sigma
_3,...,\sigma _n}^{+}\right) }{\longrightarrow }{\cal P}_{\left( \sigma
_1+\sigma _2,\sigma _3,...,\sigma _n\right) }^{\bullet }  \label{antiHol}
\end{equation}

where:

(a) the upper boldface dot in ${\rm Diff}_{\sigma _1}^{\bullet }$ in (\ref
{antiHol}) refers to the $A$-bimodule structure $\left( {\cal P}_{\left(
\sigma _2,...,\sigma _n\right) }^{\bullet },{\cal P}_{\left( \sigma
_2,...,\sigma _n\right) }^{+}\right) $;

(b) the morphisms on the right are defined in the only natural way by using
the ''gluing'' morphism of Definition \ref{alucce} : for (\ref{Hol}) we
''glue'' with respect to the last two indexes while in (\ref{antiHol}) we
''glue'' with respect to the first two. The case of ${\cal P}_{\sigma
(n)}^{\bullet }\left[ P\right] $ is proved analogously, using (\ref{relHol})
instead of (\ref{Hol}) or (\ref{relantiHol}) in place of (\ref{antiHol}).
\qquad \TeXButton{End Proof}{\endproof}\ 

\begin{remark}
\TeXButton{rem 2.12}{\label{rem 2.12}} For any $k>0$, we have ${\bf \Lambda }%
^{\left( k\right) }\simeq \dfrac{{\rm I}}{{\rm I}^{k+1}}$ where ${\rm I}$ is
the kernel of the ring multiplication $A\otimes _KA\rightarrow A$; hence $%
{\bf \Lambda }^{\left( 1\right) }\simeq \Omega _{A/K}^1$ is just the $A$%
-module of {\em K\"ahler differentials} (relative to $K$). Moreover it is
not difficult to show (\cite{KLV} p. 17) that ${\bf \Lambda }^{\left(
1,...,1\right) }\simeq {\bf \Lambda }^n\doteq {\bf \Lambda }^1\wedge \ldots
\wedge {\bf \Lambda }^1$ ($n$ times) $\simeq \Omega _{A/K}^n$ and that for
each $k,l\in {\bf N}_{+}$ the map $\Delta \longmapsto \Delta $ induces a
monomorphism $\left[ A-{\bf Mod,}A-{\bf Mod}\right] $:

\[
{\rm D}_{{\bf 1}\left( k+l\right) }\hookrightarrow {\rm D}_{{\bf 1}\left(
k\right) }\circ {\rm D}_{{\bf 1}\left( l\right) }\text{, } 
\]
whose dual representative $A$-homomorphism is just the wedge product $\wedge
:{\bf \Lambda }^k(A)\otimes _A{\bf \Lambda }^l(A)\rightarrow {\bf \Lambda }%
^{k+l}(A)$.\ 
\end{remark}

If \ ${\frak D}$ is a differentially closed subcategory, we will denote the
strict representatives in ${\frak D}$ of the relevant functors by adding $%
{\frak D}$ as a subscript to the symbol used to denote the corresponding
representative object in $A-{\bf Mod}$; for example, we write ${\bf \Lambda }%
_{{\frak D}}^{\sigma \left( n\right) }$ for the representative object in $%
{\frak D}$ of the functor ${\rm D}_{\sigma (n)}:{\frak D}\rightarrow {\frak D%
}$.

\begin{remark}
\label{grate} As we did for ${\bf J}_{{\frak D}}^{\sigma _1}={\rm Hol}_{%
{\frak D}}^{\left( \sigma _1\right) }$, we can exhibit another compatible $A$%
-module structure on ${\rm Hol}_{{\frak D}}^{\left( \sigma _1,...,\sigma
_n\right) }$, $\forall n>0$. Let $a\in A$ and $\widehat{{\rm id}}\in {\cal P}%
_{\sigma (n)}^{\bullet }\left( {\rm Hol}_{{\frak D}}^{\left( \sigma
_1,...,\sigma _n\right) }\right) $ correspond to the identity of ${\rm Hol}_{%
{\frak D}}^{\left( \sigma _1,...,\sigma _n\right) }$ under the
representability isomorphism. Since ${\cal P}_{\sigma (n)}^{\bullet }\left( 
{\rm Hol}_{{\frak D}}^{\left( \sigma _1,...,\sigma _n\right) }\right) $ and $%
{\cal P}_{\sigma (n)}^{+}\left( {\rm Hol}_{{\frak D}}^{\left( \sigma
_1,...,\sigma _n\right) }\right) $ coincide as sets, we can consider $a^{+}%
\widehat{{\rm id}}$ (multiplication in ${\cal P}_{\sigma (n)}^{+}\left( {\rm %
Hol}_{{\frak D}}^{\left( \sigma _1,...,\sigma _n\right) }\right) $) as an $A$%
-endomorphism of ${\rm Hol}_{{\frak D}}^{\left( \sigma _1,...,\sigma
_n\right) }$. It is easy to verify that this choice defines another $A$%
-module structure on ${\rm Hol}_{{\frak D}}^{\left( \sigma _1,...,\sigma
_n\right) }$, denoted by ${\rm Hol}_{{\frak D},+}^{\left( \sigma
_1,...,\sigma _n\right) }$ and that $\left( {\rm Hol}_{{\frak D}}^{\left(
\sigma _1,...,\sigma _n\right) },{\rm Hol}_{{\frak D},+}^{\left( \sigma
_1,...,\sigma _n\right) }\right) $ is an $A$-bimodule.
\end{remark}

\begin{lemma}
\label{pimpa}If $P\in Ob\left( {\frak D}\right) $ we have a canonical
isomorphism ${\rm Hol}_{{\frak D}}^{\left( \sigma _1,...,\sigma _n\right)
}\left[ P\right] \simeq {\rm Hol}_{{\frak D},+}^{\left( \sigma _1,...,\sigma
_n\right) }\otimes ^{\bullet }P$ in ${\frak D}$.
\end{lemma}

\TeXButton{Proof}{\proof}Let $Q$ be an object in ${\frak D}$. By lemma \ref
{lem 2.3}, we have 
\[
Hom\left( {\rm Hol}_{{\frak D},+}^{\left( \sigma _1,...,\sigma _n\right)
}\otimes ^{\bullet }P,Q\right) \simeq Hom^{\bullet }\left( P,Hom^{+}\left( 
{\rm Hol}_{{\frak D}}^{\left( \sigma _1,...,\sigma _n\right) },Q\right)
\right) \simeq 
\]
\[
\ \simeq Hom^{\bullet }\left( P,{\cal P}_{\sigma (n)}^{+}\left( Q\right)
\right) 
\]
and by definition of ${\cal P}_{\sigma (n)}^{+}$ and example \ref{es 2.2} 
\[
Hom^{\bullet }\left( P,{\cal P}_{\sigma (n)}^{+}\left( Q\right) \right)
\simeq {\cal P}_{\left( \sigma _1\right) }^{\bullet }\left[ P\right] \left( 
{\cal P}_{\left( \sigma _2,...,\sigma _n\right) }^{}\left( Q\right) \right) =%
{\cal P}_{\left( \sigma _1,...,\sigma _n\right) }^{\bullet }\left[ P\right]
\left( Q\right) \text{.} 
\]
\TeXButton{End Proof}{\endproof}

\begin{remark}
\TeXButton{jam}{\label{jam}} For any differentiable closed subcategory $%
{\frak D}\subseteq A-{\bf Mod}$ it is still true, as in the case ${\frak D}%
=A-{\bf Mod}$, that ${\bf J}_{{\frak D}}^k\left( P\right) $ is generated as
an $A$-module by 
\[
\ \ \left\{ j_k\left( p\right) \mid p\in P\right\} . 
\]
In fact, let ${\bf J}_{{\frak D}}^k\left( P\right) ^{\sim }$ denote the $A$%
-submodule of ${\bf J}_{{\frak D}}^k\left( P\right) $ generated by $\left\{
j_k\left( p\right) \mid p\in P\right\} $; this is still an object of ${\frak %
D}$ by Definition \ref{def 2.1} (e). Now, the composition 
\[
P\stackrel{j_k^{{\frak D}}}{\longrightarrow }{\bf J}_{{\frak D}}^k\left(
P\right) \stackrel{\pi }{\twoheadrightarrow }\frac{{\bf J}_{{\frak D}%
}^k\left( P\right) }{{\bf J}_{{\frak D}}^k\left( P\right) ^{\sim }} 
\]
is the DO of order $\leq k$ corresponding to $\pi $ under the isomorphism 
\[
Hom_A\left( {\bf J}_{{\frak D}}^k\left( P\right) ,{\bf J}_{{\frak D}%
}^k\left( P\right) /{\bf J}_{{\frak D}}^k\left( P\right) ^{\sim }\right)
\simeq {\rm Diff}_k\left( P,{\bf J}_{{\frak D}}^k\left( P\right) /{\bf J}_{%
{\frak D}}^k\left( P\right) ^{\sim }\right) . 
\]
But $\pi \circ j_k^{{\frak D}}$ is zero hence $\pi =0$ and we conclude. Note
that this also shows that the canonical morphism 
\[
{\bf J}^k\left( P\right) \rightarrow {\bf J}_{{\frak D}}^k\left( P\right) 
\text{ } 
\]
is an $A$-epimorphism.
\end{remark}

If $P\in Ob({\frak D})$, the monomorphism in $\left[ {\frak D}{\bf ,}{\frak D%
}\right] $:

\[
{\rm Diff}_s(P,\cdot )\subset {\rm Diff}_t(P,\cdot )\text{,\qquad }t\geq
s,\quad 
\]
gives rise to a ${\frak D}$-epimorphism (also an $A$-epimorphism by Remark 
\ref{jam}) between representative objects:

\[
\pi _{t,s}(P):{\bf J}_{{\frak D}}^t(P)\rightarrow {\bf J}_{{\frak D}}^s(P) 
\]
which fits in the commutative diagram

\[
\begin{array}{ccccc}
P & \stackrel{j_s(P)}{\longrightarrow } & {\bf J}_{{\frak D}}^s(P) &  &  \\ 
& \stackunder{j_t(P)}{\searrow } & \uparrow & \pi _{t,s}(P) &  \\ 
&  & {\bf J}_{{\frak D}}^t(P) &  & 
\end{array}
\]
The rule $P\longmapsto {\bf J}_{{\frak D}}^t(P)$ defines in the obvious way
a (covariant) functor ${\frak D}\rightarrow {\frak D}$ (\cite{Kr} or \cite
{KLV}).\ 

The following example shows the importance\ of the appropriate choice of the
differentially closed subcategory of $A-{\bf Mod}$ in determining the ''
geometrical effectiveness'' and size of the representative objects of the
relevant functors.

\begin{example}
Let $M$ be a smooth real manifold (which we assume Hausdorff and with a
countable basis), $K={\Bbb R}$ and $A=C^\infty \left( M;{\Bbb R}\right) $.
Then (\cite{MVi}) ${\bf \Lambda }^{\sigma \left( n\right) }$ is in general
neither projective nor of finite type over $A$: in particular, when $\sigma
=\left( 1,...,1,...\right) $, {\it it does not} coincide with the $A$-module
of differential $n$-forms on the manifold $M$. To obtain these
''geometrical'' objects we must choose an appropriate subcategory of $A-{\bf %
Mod}$ : in our approach, choosing a ''geometry'' is equivalent to select a
differentially closed subcategory ${\frak D}$. For finite dimensional (real)
differential geometry we may choose ${\frak D}\doteq A-{\bf Mod}_{geom}$,
the full subcategory of geometric $A$-modules, i.e. of $A$-modules $P$ such
that $\bigcap\limits_{x\in M}I_xP=(0)$, $I_x$ being the maximal ideal of
smooth functions on $M$ vanishing at $x\in M$.

Note that $A-{\bf Mod}_{geom}\supseteq A-{\bf Mod}_{pr,\text{ }f.t.}$, the
full subcategory of projective $A$-modules of finite type, since $A$ itself
is a geometric $A$-module; however, $A-{\bf Mod}_{pr,\text{ }f.t.}$ is not
differentially closed because it is not abelian (and does not satisfy (e)).
Another reason that makes us prefer working with the bigger $A-{\bf Mod}%
_{geom}$ is its better functoriality with respect to change of algebras
induced by pull backs of smooth mappings of manifolds\footnote{%
If $f:M\rightarrow N$ is a smooth map and $P$ is a geometric $C^\infty (M)$%
-module then $P$ is still geometric when viewed as a $C^\infty (N)$-module
via the pull back $f^{*}:C^\infty (N)\longrightarrow C^\infty (M)$.
Projectivity is not preserved, instead.}. $A-{\bf Mod}_{geom}$ is
differentially closed due to the fact that the ''geometrization'' functor 
\begin{eqnarray*}
\left( \text{ }\cdot \text{ }\right) _{geom}:A-{\bf Mod}\longrightarrow A-%
{\bf Mod}_{geom} \\
P\longmapsto \text{ }P_{geom}\doteq \frac P{\bigcap\limits_{x\in M}I_xP}
\end{eqnarray*}
sends representative objects in $A-{\bf Mod}_{geom}$ to representative
objects in $A-{\bf Mod}_{geom}$ for all the relevant functors (\cite{KLV}).
''Geometrical'' objects are obtained as representative objects; for example $%
{\bf \Lambda }_{geom.}^{\left( 1,...,1\right) }$, with $(1,...,1)\in {\bf N}%
_{}^k$, is isomorphic to $\Gamma \left( \bigwedge\limits^kT^{*}M\right) $,
the module of sections of the $k$-th exterior power of the cotangent bundle
of $M$, i.e. the module of $k$-differential forms of $M$.
\end{example}

It is possible to encode the ''{\em smoothness}'' of the geometry we want to
describe, completely in the choice of the differentially closed subcategory:

\begin{definition}
A differentially closed subcategory ${\frak D}$ of $A-{\bf Mod}$ is called 
{\em smooth} if $\Lambda _{{\frak D}}^{\left( 1\right) }$ is a projective $A$%
-module of finite type.
\end{definition}

\begin{examples}
(i) If $K$ is an algebraically closed field of zero characteristic and $A$
is the coordinate ring of a {\it regular} affine $K$-variety, then ${\frak D}%
=A-{\bf Mod}$ is smooth (\cite{Ha}, II.8).

(ii) If $M$ is a smooth manifold and $A=C^\infty \left( M;{\bf R}\right) $
then $A-{\bf Mod}_{geom}$ is smooth while $A-{\bf Mod}$ is not.\ 
\end{examples}

It can be proved (as in the proof of Theorem \ref{thm 2.7}) that if ${\frak D%
}$ is smooth then {\em all the representative objects of relevant functors
are indeed projective and of finite type as }$A${\em -modules}. However, we
want to stress that since representative objects may be constructed also in
non-smooth cases, our approach works also in describing singular and even
infinite dimensional geometrical situations. However, to resort with useful
objects one has to make an adequate choice of ${\frak D}$.\ 

The following proposition will be useful in the next sections:

\begin{proposition}
\label{stop} Let ${\frak D}$ be smooth. Then

(i) ${\bf \Lambda }_{{\frak D}}^{{\bf 1}_n}=(0)$ for $n>>0$;

(ii) If $P\in Ob\left( {\frak D}\right) $, ${\rm Hol}_{{\frak D}}^{{\bf 1}%
_n}\left[ P\right] =(0)$ for $n>>0$.
\end{proposition}

\TeXButton{Proof}{\proof} By lemma \ref{pimpa} and remark \ref{new}, (ii)
follows from (i). Remark \ref{rem 2.12} together with the fact that ${\bf %
\Lambda }_{{\frak D}}^{(1)}$ is of finite type proves (i). 
\TeXButton{End Proof}{\endproof}

\section{Higher de Rham complexes}

In this Section we use the functors introduced in Section 1 and their
representative objects (Section 2) to build higher order analogs of the de
Rham complex. Their cohomology will be studied in Section 4.

Let ${\frak D}$ be a differentially closed subcategory of $A-{\bf Mod}$. The
dual representative of the monomorphism in $\left[ {\frak D},{\frak D}%
\right] $:

\[
{\rm D}_{(\sigma \left( n\right) ,k)}\hookrightarrow {\rm D}_{\sigma \left(
n\right) }^{\bullet }\circ {\rm Diff}_k^{(+)}\text{, \qquad }\sigma \left(
n\right) \in {\bf N}_{+}^n,\text{ }k\in {\bf N}_{+}\text{,} 
\]
is a ${\frak D}$-epimorphism:

\[
{\bf J}_{{\frak D}}^k({\bf \Lambda }_{{\frak D}}^{\sigma \left( n\right)
})\rightarrow {\bf \Lambda }_{{\frak D}}^{(\sigma \left( n\right) ,k)}\text{;%
} 
\]
define $d_{(\sigma \left( n\right) ,k)}^{{\frak D}}$ to be the composition

\begin{equation}
{\bf \Lambda }_{{\frak D}}^{\sigma \left( n\right) }\stackrel{j_k}{%
\longrightarrow }{\bf J}_{{\frak D}}^k({\bf \Lambda }_{{\frak D}}^{\sigma
\left( n\right) })\rightarrow {\bf \Lambda }_{{\frak D}}^{(\sigma \left(
n\right) ,k)}\text{.}  \label{scioglimento}
\end{equation}
Obviously, $d_{(\sigma \left( n\right) ,k)}^{{\frak D}}$ is a DO of order $%
\leq k$.

\smallskip\ 

\begin{definition}
If $\sigma \in {\bf N}_{+}^\infty $, the sequence in ${\bf DIFF}_A$%
\begin{equation}
0\rightarrow A\stackrel{d_{(\sigma _1)}^{{\frak D}}}{\rightarrow }\Lambda _{%
{\frak D}}^{(\sigma _1)}\stackrel{d_{\left( \sigma _1,\sigma _2\right) }^{%
{\frak D}}}{\longrightarrow }\Lambda _{{\frak D}}^{\sigma (2)}\rightarrow
\ldots \stackrel{d_{\sigma \left( k\right) }^{{\frak D}}}{\longrightarrow }%
\Lambda _{{\frak D}}^{\sigma \left( k\right) }\rightarrow \ldots
\label{dRsup}
\end{equation}
is called {\em higher de Rham sequence}{\it \ }of type $\sigma $ of the $K$%
-algebra $A$ and is denoted by ${\bf dR}_\sigma ^{{\frak D}}(A)$ or simply
by ${\bf dR}_\sigma ^{{\frak D}}$; each $d_{\sigma \left( k\right) }^{{\frak %
D}}$, $k>0$, is called {\em higher de Rham differential} and is a DO of
order $\leq \sigma _k$.
\end{definition}

\begin{remark}
When $\sigma ={\bf 1}\in {\bf N}_{+}^\infty $, the corresponding de Rham
sequence is called {\em ordinary}. In this case we write ${\bf \Lambda }_{%
{\frak D}}^k$ for ${\bf \Lambda }_{{\frak D}}^{\left( 1,...,1\right) }$ , $%
\left( 1,...,1\right) \in {\bf N}_{+}^k$, $\forall k>0$, so that:

\begin{equation}
{\bf dR}_{(\underline{1})}^{{\frak D}}\equiv {\bf dR}^{{\frak D}}{\bf :\quad 
}0\rightarrow A\stackrel{d}{\rightarrow }{\bf \Lambda }_{{\frak D}}^1%
\stackrel{d}{\rightarrow }{\bf \Lambda }_{{\frak D}}^2\rightarrow \ldots 
\stackrel{d}{\rightarrow }{\bf \Lambda }_{{\frak D}}^k\rightarrow \ldots
\label{dR}
\end{equation}
and each differential is a DO of order $\leq 1$.
\end{remark}

Each $d^{{\frak D}}$ in (\ref{dRsup}) is in fact a differential according to
the following:\ 

\begin{proposition}
$\forall \sigma \in {\bf N}_{+}^\infty $ the higher de Rham sequence ${\bf dR%
}_\sigma ^{{\frak D}}$ is a complex.\ 
\end{proposition}

\TeXButton{Proof}{\proof} Let $n\geq 0$ and consider the diagram defining
two consecutive higher de Rham differentials:

\[
\begin{array}{cccccc}
{\bf \Lambda }_{{\frak D}}^{\sigma \left( n\right) } & \stackrel{\doteq
d_{\sigma \left( n+1\right) }^{{\frak D}}}{\longrightarrow } & {\bf \Lambda }%
_{{\frak D}}^{\sigma \left( n+1\right) } & \stackrel{\doteq d_{\sigma \left(
n+2\right) }^{{\frak D}}}{\longrightarrow } & {\bf \Lambda }_{{\frak D}%
}^{\sigma \left( n+2\right) } &  \\ 
& \searrow ^{j_{\sigma _{n+1}}} & \nearrow _{\pi _1} & \searrow ^{j_{\sigma
_{n+2}}} & \nearrow _{\pi _2} &  \\ 
& {\bf J}_{{\frak D}}^{\sigma _{n+1}}({\bf \Lambda }_{{\frak D}}^{\sigma
\left( n\right) }) &  & {\bf J}_{{\frak D}}^{\sigma _{n+2}}({\bf \Lambda }_{%
{\frak D}}^{\sigma \left( n+1\right) }) &  & 
\end{array}
\text{.} 
\]
Since $d_{\sigma \left( n+2\right) }^{{\frak D}}\circ d_{\sigma \left(
n+1\right) }^{{\frak D}}\equiv $ $\pi _2\circ j_{\sigma _{n+2}}\circ \pi
_1\circ j_{\sigma _{n+1}}$ is a DO of order $\leq k+l$ , there exists a
unique $A$-homomorphism

\[
\varphi _{d_{\sigma \left( n+2\right) }^{{\frak D}}\circ d_{\sigma \left(
n+1\right) }^{{\frak D}}}:{\bf J}_{{\frak D}}^{\sigma _{n+2}+\sigma _{n+1}}(%
{\bf \Lambda }_{{\frak D}}^{\sigma \left( n\right) })\rightarrow {\bf %
\Lambda }_{{\frak D}}^{\sigma \left( n+2\right) } 
\]
which makes the following diagram commutative:

\[
\text{\TeXButton{Vtriangle1}
{\settriparms[1`1`-1;500]
\Vtriangle[\QTR{bf}{\Lambda }_{\QTR{frak}{D}}^{\sigma \left( n\right) }`\QTR{bf}{\Lambda }_{\QTR{frak}{D}}^{\sigma \left( n+2\right) }`\QTR{bf}{J}_{\QTR{frak}{D}}^{\sigma _{n+2}+\sigma _{n+1}}(\QTR{bf}{\Lambda }_{\QTR{frak}{D}}^{\sigma \left( n\right) });d_{\sigma \left( n+2\right) }^{\QTR{frak}{D}}\circ d_{\sigma \left( n+1\right) }^{\QTR{frak}{D}}`j_{\sigma _{n+2}+\sigma _{n+1}}`\varphi _{d_{\sigma \left( n+2\right) }^{\QTR{frak}{D}}\circ d_{\sigma \left( n+1\right) }^{\QTR{frak}{D}}}]}%
.} 
\]

It is not difficult to check that $\varphi _{d_{\sigma \left( n+2\right) }^{%
{\frak D}}\circ d_{\sigma \left( n+1\right) }^{{\frak D}}}$ is just the dual
representative of the composition:

\[
{\rm D}_{\sigma \left( n+2\right) }\hookrightarrow {\rm D}_{\sigma \left(
n+1\right) }^{\bullet }\circ {\rm Diff}_{\sigma _{n+2}}^{(+)}\longrightarrow 
{\rm D}_{\sigma \left( n\right) }^{\bullet }\circ {\rm Diff}_{\sigma
_{n+1}+\sigma _{n+2}}^{(+)} 
\]
which is immediately checked to be zero; therefore $\varphi _{d_{\sigma
\left( n+2\right) }^{{\frak D}}\circ d_{\sigma \left( n+1\right) }^{{\frak D}%
}}=0$ and $d_{\sigma \left( n+2\right) }^{{\frak D}}\circ d_{\sigma \left(
n+1\right) }^{{\frak D}}=0$ as well.\qquad \TeXButton{End Proof}{\endproof}\ 

\begin{remark}
\TeXButton{rem 3.6}{\label{rem 3.6}}({\bf i}) Let ${\frak D}=A-{\bf Mod}$.
By induction on $n$ we can prove (using for $n=1$ the explicit description
of ${\bf \Lambda }^{\left( \sigma _1\right) }$ given in Section 2) that
formula (\ref{scioglimento}) implies that ${\bf \Lambda }^{\sigma \left(
n\right) }$ is generated by the set 
\[
\left\{ d_{\sigma \left( n\right) }\left( a_1d_{\sigma \left( n-1\right)
}\left( a_2\ldots d_{\sigma \left( 1\right) }^{}\left( a_n\right) \ldots
\right) \right) \mid a_1,a_2,\ldots ,a_n\in A\right\} . 
\]
In fact, ${\bf J}^{\sigma _n}({\bf \Lambda }^{\sigma \left( n-1\right) })$
is known to be generated over $A$ by the elements $j_{\sigma _n}\left(
\omega \right) $, $\omega \in {\bf \Lambda }^{\sigma \left( n-1\right) }$
and ${\bf J}^{\sigma _n}({\bf \Lambda }^{\sigma \left( n-1\right)
})\rightarrow {\bf \Lambda }^{\sigma \left( n\right) }$ is an epimorphism.
This result still holds for ${\bf \Lambda }_{{\frak D}}^{\sigma \left(
n\right) }$ with $d$ replaced by $d^{{\frak D}}$, ${\frak D}\subseteq A-{\bf %
Mod}$ being any differentiable closed subcategory: the proof is analogous to
the argument used in Remark \ref{jam} (a). This also shows that the
canonical morphisms 
\[
\text{ \qquad }{\bf \Lambda }^{\sigma \left( n\right) }\rightarrow {\bf %
\Lambda }_{{\frak D}}^{\sigma \left( n\right) } 
\]
are $A$-epimorphisms.

({\bf ii}) In the ''ordinary'' case $\left( \sigma _1,...,\sigma _n\right)
=\left( 1,...,1\right) \equiv {\bf 1}\left( n\right) $, the $A$-module
structure of ${\rm Hol}_{{\frak D},+}^{{\bf 1}_n}$ (remark \ref{grate}) can
be expressed via the isomorphism (remark \ref{new}) 
\[
{\rm Hol}_{{\frak D}}^{\left( \sigma _1,...,\sigma _n\right) }\simeq {\bf %
\Lambda }_{{\frak D}}^{\left( \sigma _2,...,\sigma _n\right) }\oplus {\bf %
\Lambda }_{{\frak D}}^{\left( \sigma _1,...,\sigma _n\right) } 
\]
as $a^{+}\left( \rho ,\omega \right) =\left( a\rho ,a\omega +\left(
d_{\left( 1\right) }a\right) \wedge \rho \right) $.

({\bf iii}) If ${\frak D}=A-{\bf Mod}$, ${\bf dR}_{{\bf 1}}$ coincides with
the usual algebraic de Rham complex of the $K$-algebra $A$ (\cite{Bou III}
and \cite{Bou X}).

({\bf iv}) If $K={\Bbb R}$, $M$ is a smooth manifold, $A=C^\infty \left( M;%
{\Bbb R}\right) $ and ${\frak D}$ is the category of geometric $A$-modules
(see Section 2) then ${\bf dR}_{{\bf 1}}^{{\frak D}}$ is the geometric de
Rham complex on $M$. It turns out that any natural differential operator
occurring in differential geometry can be recovered functorially using our
approach: see \cite{VV} for the case of the Lie derivative and the
corresponding homotopy formula.
\end{remark}

If $\sigma ,\tau \in {\bf N}_{+}^\infty $ with $\sigma \geq \tau $ (i.e. $%
\sigma _i\geq \tau _i$, $\forall i\geq 1$), then for each $n>0$ we have a
monomorphism ${\rm D}_{\tau \left( n\right) }\hookrightarrow {\rm D}_{\sigma
\left( n\right) }$ in $\left[ {\frak D},{\frak D}\right] $; this induces a $%
{\frak D}$-epimorphism on representatives ${\bf \Lambda }_{{\frak D}%
}^{\sigma \left( n\right) }\rightarrow {\bf \Lambda }_{{\frak D}}^{\tau
\left( n\right) }$, $\forall n>0$. By \ref{rem 3.6}, this is also an $A$%
-epimorphism. All these epimorphisms commute with higher de Rham
differentials and therefore define a morphism of complexes

\begin{equation}
{\bf dR}_\sigma ^{{\frak D}}\rightarrow {\bf dR}_\tau ^{{\frak D}}
\label{luppolo}
\end{equation}
(if $\sigma \geq \tau $). So we can consider the ($A$-epimorphic) inverse
system $\left\{ {\bf dR}_\sigma ^{{\frak D}}\right\} _{\sigma \in {\bf N}%
_{+}^\infty }$ and give the following:

\begin{definition}
\TeXButton{infty}{\label{infty}}The {\it infinitely prolonged} (or, simply, 
{\it infinite}) {\it de Rham complex} of the $K$-algebra $A$, is the complex
in {\bf K}$\left( K-{\bf Mod}\right) $%
\begin{equation}
\begin{array}{c}
{\bf dR}_\infty ^{{\frak D}}(A)\doteq \underleftarrow{\lim }_{\sigma \in 
{\bf N}_{+}^\infty }{\bf dR}_\sigma ^{{\frak D}}(A), \\ 
{\bf dR}_\infty ^{{\frak D}}(A):\text{ }0\rightarrow A\stackrel{d_{\left(
\infty \right) }}{\longrightarrow }{\bf \Lambda }_{{\frak D}}^{\left( \infty
\right) }\stackrel{d_{\left( \infty ,\infty \right) }}{\longrightarrow }{\bf %
\Lambda }_{{\frak D}}^{\left( \infty ,\infty \right) _2}\rightarrow \cdots
\rightarrow {\bf \Lambda }_{{\frak D}}^{\left( \infty ,\ldots ,\infty
\right) _n}\rightarrow \cdots
\end{array}
\label{infinito}
\end{equation}
where ${\bf \Lambda }_{{\frak D}}^{\left( \infty ,\ldots ,\infty \right)
_n}\doteq \underleftarrow{\lim }_{\sigma \left( n\right) \in {\bf N}_{+}^n}%
{\bf \Lambda }_{{\frak D}}^{\sigma \left( n\right) }$, $\forall n>0$.
\end{definition}

\begin{remark}
\label{photos}\TeXButton{photos}{\label{photos}}{\em Three descriptions of
DO's between strict representative objects.}

We work in a fixed differentially closed subcategory ${\frak D}$ of $A-{\bf %
Mod}$ and all representative objects will be in ${\frak D}$.

Let $F_1$ and $F_2$ be representative objects of differential functors $%
{\cal F}_1$ and ${\cal F}_2$, respectively. Suppose that ${\cal F}_1$ has an
associated functor ${\cal F}_1^{\bullet }$ (having as domain $A-{\bf BiMod}_{%
{\frak D}}$) such that ${\cal F}_1^{\bullet }\left( {\rm Diff}_k^{\left(
+\right) }\right) $ is strictly representable by ${\bf J}^k\left( F_1\right) 
$: this is the case, for example, of ${\cal F}_1={\rm D}_{\sigma \left(
n\right) }$ or ${\rm Diff}_l$. Let

\begin{equation}
\Delta :F_1\longrightarrow F_2  \label{are}
\end{equation}
be a DO of order $\leq k$. Then, there exists a unique $A$-homomorphism (%
\cite{Kr}: jet-associated to $\Delta $)

\begin{equation}
f_\Delta :{\bf J}^k\left( F_1\right) \longrightarrow F_2  \label{ere}
\end{equation}
which represents $\Delta $ by duality: $\Delta =f_\Delta \circ j_k\left(
F_1\right) $. Since ${\bf J}^k\left( F_1\right) $ is the representative
object of ${\cal F}_1^{\bullet }\left( {\rm Diff}_k^{\left( +\right)
}\right) $, $f_\Delta $ defines a unique morphism in $\left[ {\frak D},%
{\frak D}\right] $: 
\begin{equation}
f^\Delta :{\cal F}_2\longrightarrow {\cal F}_1^{\bullet }\left( {\rm Diff}%
_k^{\left( +\right) }\right) \text{,}  \label{ire}
\end{equation}
called {\em generator} {\em morphism} of $\Delta $.

Formulas (\ref{ere}) and (\ref{ire}) give two different descriptions of a DO
between representative objects. Formula (\ref{ire}) allows one to identify
it with a functorial morphism which, as a rule, may be established in a
straightforward way and can then be used to define the corresponding natural
DO (\ref{are}). The following examples show this procedure at work in two
canonical cases; we assume for simplicity ${\frak D}=A-{\bf Mod}$.

(i) {\em Higher de Rham differential} $d_{\sigma \left( n\right) }$.

If ${\cal F}_2$ $={\rm D}_{\sigma \left( n\right) }$ , ${\cal F}_1={\rm D}%
_{\sigma \left( n-1\right) }$ , $k=\sigma _n$ and we take for (\ref{ire})
the natural inclusion 
\[
{\rm D}_{\sigma \left( n\right) }\hookrightarrow {\rm D}_{\sigma \left(
n-1\right) }^{\bullet }\left( {\rm Diff}_{\sigma _n}^{\left( +\right)
}\right) , 
\]
then $d_{\sigma \left( n\right) }:{\bf \Lambda }^{\sigma \left( n-1\right)
}\rightarrow {\bf \Lambda }^{\sigma \left( n\right) }$ is the corresponding
DO (\ref{are}).

(ii) {\em ''Absolute'' jet-operator} $j_k$.

In this almost tautological case, ${\cal F}_1=Hom_A\left( A,\cdot \right)
\equiv {\rm Diff}_0$ and ${\cal F}_2$ $={\rm Diff}_k\equiv $

$Hom_A^{\bullet }\left( A,\cdot \right) $ $\left( {\rm Diff}_k^{\left(
+\right) }\right) $ ; if we take (\ref{ire}) to be the identity 
\[
{\rm Id}:{\rm Diff}_k\rightarrow Hom_A^{\bullet }\left( A,\cdot \right)
\left( {\rm Diff}_k^{\left( +\right) }\right) \equiv {\rm Diff}_k, 
\]
then (\ref{are}) is just $j_k:A\rightarrow $ ${\bf J}^k$.
\end{remark}

\medskip 

\ Rigidity of higher de Rham cohomology

In this Section we prove the main result of this paper i.e. that in the
smooth case the higher-order de Rham cohomologies coincide with the ordinary
( i.e. lowest order) one. Essentially this amounts to a fairly intuitive
assert: raising the order of the natural DOs involved in the ${\bf dR}$%
-complexes does not change the cohomological information, provided the
situation in which we are working is smooth.

In this Section (and in the Appendix), $A$ is a $K$-algebra of zero
characteristic, containing $K$ as a subring and ${\frak D}$ a differentially
closed {\it smooth} subcategory of $A-{\bf Mod}$. As in the previous
Section, all representative objects, unless otherwise stated, will be
considered in ${\frak D}$.

Smoothness of ${\frak D}$ implies that for any $k,l\geq 0$, the gluing
morphism in $\left[ {\frak D},{\frak D}\right] $ 
\[
{\rm Diff}_k^{\bullet }\circ {\rm Diff}_l^{\left( +\right) }\stackrel{C_{k,l}%
}{\longrightarrow }{\rm Diff}_{k+l}^{} 
\]
is surjective i.e. that any DO can be expressed as composition of lower
order ones. This can be seen as follows. Let us fix $k$ and proceed by
induction on $l$. The case $l=0$ is trivial since ${\rm Diff}_0={\rm id}_{%
{\frak D}}$. To prove the inductive step let us consider the commutative
diagram 
\[
\begin{tabular}{lll}
${\rm Diff}_k^{\bullet }\circ {\rm Diff}_{l+1}^{\left( +\right) }$ & $%
\stackrel{C_{k,l+1}}{\longrightarrow }$ & ${\rm Diff}_{k+l+1}$ \\ 
$\qquad \uparrow $ &  & $\quad \uparrow $ \\ 
${\rm Diff}_k^{\bullet }\circ {\rm Diff}_l^{\left( +\right) }$ & $%
\stackunder{C_{k,l}}{\longrightarrow }$ & ${\rm Diff}_{k+l}$%
\end{tabular}
\]
(where the vertical arrows are natural inclusions) and suppose $C_{k,l}$ is
epic. Passing to the corresponding diagram of reperesentative objects
completing it with kernels, we get a commutative diagram (\cite{KLV} p. 52)
with exact columns 
\[
\begin{tabular}{lll}
$0$ &  & $0$ \\ 
$\downarrow $ &  & $\downarrow $ \\ 
${\rm S}^{k+l+1}\left( {\bf \Lambda }^1\right) $ & $\stackrel{\rho }{%
\longrightarrow }$ & ${\rm S}^{l+1}\left( {\bf \Lambda }^1\right) \otimes 
{\bf J}^k$ \\ 
$\downarrow $ &  & $\downarrow $ \\ 
${\bf J}^{k+l+1}$ & $\stackrel{C^{k,l+1}}{\longrightarrow }$ & ${\bf J}%
^{l+1}\left( {\bf J}^k\right) $ \\ 
$\downarrow $ &  & $\downarrow $ \\ 
${\bf J}^{k+l}$ & $\stackunder{C^{k,l}}{\longrightarrow }$ & ${\bf J}%
^l\left( {\bf J}^k\right) $ \\ 
$\downarrow $ &  & $\downarrow $ \\ 
$0$ &  & $0$%
\end{tabular}
\]
where ${\rm S}^r$ denotes the $r$-th symmetric power and $C^{s,t}$ is the
dual-representative of $C_{s,t}$. By duality it is enough to prove that $%
C^{k,l+1}$ is monic. By induction hypothesis $C^{k,l}$ is monic so we are
reduced to showing that $\rho $ is monic. It is not difficult to prove (e.g.
again by induction on $l$) that $\rho $ is just the composition 
\[
\begin{tabular}{lll}
${\rm S}^{k+l+1}\left( {\bf \Lambda }^1\right) $ & $\stackrel{\alpha }{%
\longrightarrow }\stackrel{}{{\rm S}^{l+1}\left( {\bf \Lambda }^1\right)
\otimes {\rm S}^k\left( {\bf \Lambda }^1\right) \stackrel{\beta }{%
\longrightarrow }}$ & ${\rm S}^{l+1}\left( {\bf \Lambda }^1\right) \otimes 
{\bf J}^k$%
\end{tabular}
\]
where $\alpha :\omega _1\cdots \omega _{k+l+1}\longmapsto \sum \left( \omega
_{i_1}\cdots \omega _{i_{l+1}}\right) \otimes \left( \omega _{j_1}\cdots
\omega _{j_k}\right) $ where the sum is extended to all partitions $\left(
\left( i_1,...,i_{l+1}\right) ,\left( j_1,...,j_k\right) \right) $ of $%
\left\{ 1,...,k+l+1\right\} $ of (ordered) length $(l+1,k)$ and $\beta ={\rm %
id}_{{\rm S}^{l+1}\left( {\bf \Lambda }^1\right) }\otimes i$ with 
\[
i:{\rm S}^k\left( {\bf \Lambda }^1\right) \hookrightarrow {\bf J}^k 
\]
the inclusion of the kernel of ${\bf J}^k\rightarrow {\bf J}^{k-1}$ (\cite
{KLV}, p. 52). $A$ is of zero characteristic hence $\alpha $ is well defined
and monic; ${\bf \Lambda }^1$ is projective hence $\beta $ is monic too.
Thus $\rho $ is monic and we conclude.

As a consequence $\forall n\geq 1$, we have the following short exact
sequence in $\left[ {\frak D},{\frak D}\right] $:

\begin{equation}
0\rightarrow {\rm D}_{\sigma \left( n\right) }\hookrightarrow {\rm D}%
_{\sigma \left( n-1\right) }^{\bullet }\circ {\rm Diff}_{\sigma _n}^{\left(
+\right) }\rightarrow {\rm D}_{\left( \sigma _1,...,\sigma _{n-2},\sigma
_{n-1}+\sigma _n\right) }\rightarrow 0  \label{(13)}
\end{equation}
(the new fact is that the last arrow of the sequence is epic since it is
induced by the gluing morphism).

The $n$-th cohomology $K$-module of the complex 
\[
{\bf dR}_\sigma :\quad 0\rightarrow A\stackrel{d_{\sigma \left( 1\right) }}{%
\longrightarrow }{\bf \Lambda }^{\sigma \left( 1\right) }\rightarrow \ldots
\rightarrow {\bf \Lambda }^{\sigma \left( n\right) }\stackrel{d_{\sigma
\left( n+1\right) }}{\longrightarrow }{\bf \Lambda }^{\sigma \left(
n+1\right) }\rightarrow \ldots 
\]
is denoted by:

\[
H_\sigma ^n\doteq \frac{\ker \left( d_{\sigma \left( n+1\right) }\right) }{%
im\left( d_{\sigma \left( n\right) }\right) }\equiv H^n\left( {\bf dR}%
_\sigma \right) \text{.} 
\]
Since $H_\sigma ^n$ only depends on $\sigma \left( n+1\right) $, we will
write also $H_{\sigma \left( n+1\right) }^n$ in place of $H_\sigma ^n$.

Note that in the situation of Remark \ref{rem 3.6} ({\bf iv}), $H_{\sigma
\left( n+1\right) }^n$ is the $n$-th de Rham cohomology ${\Bbb R}$-vector
space of the smooth manifold $M$.

\smallskip\ 

The rest of this Section will be devoted to proving the following result:

\begin{theorem}
\TeXButton{main}{\label{main}}\label{t1}{\rm (''Smooth'' rigidity of higher
de Rham cohomologies )}

If ${\frak D}$ is a smooth subcategory of $A-{\bf Mod}$, then, for each $%
\tau ,\sigma \in {\bf N}_{+}^\infty $ with $\tau \geq \sigma $, the
canonical ${\frak D}$-epimorphism (\ref{luppolo}):

\[
{\bf dR}_\tau \rightarrow {\bf dR}_\sigma 
\]
is a quasi-isomorphism; so:

\begin{equation}
H_\sigma ^n\simeq H_\tau ^n\text{, }\forall n\geq 0\text{.}  \label{(14)}
\end{equation}
\ 
\end{theorem}

\begin{corollary}
({\bf i}) If $M$ is a {\it smooth} manifold, $A=C^\infty \left( M;{\Bbb R}%
\right) $ and ${\frak D}=C^\infty \left( M;{\Bbb R}\right) -{\bf Mod}_{geom}$%
, then the higher de Rham cohomologies coincide with the standard de Rham
cohomology of $M$.

({\bf ii}) If $K$ is an algebraically closed field of zero characteristic
and $A$ is the coordinate ring of a {\it regular} affine variety over $K$,
then the higher de Rham cohomologies coincide with the standard algebraic
one.
\end{corollary}

Note that the last Corollary is false, in general, for a singular manifold
or a non-regular affine variety.

The strategy of the proof of Theorem \ref{t1} is the following.

Keeping $n\geq 0$ fixed, we prove the thesis by reducing, step by step, each
entry of $\sigma \left( n+1\right) $ to $1$, starting from $\sigma _{n+1}$,
i.e. we prove the chain of isomorphisms

\begin{equation}
H_{\sigma \left( n+1\right) }^n\simeq H_{\left( \sigma \left( n\right)
,1\right) }^n\simeq H_{\left( \sigma \left( n-1\right) ,1,1\right) }^n\simeq
\cdots \simeq H_{\left( \sigma _1,1,\ldots ,1\right) }^n\simeq H_{dR}^n
\label{sediz}
\end{equation}
where $H_{dR}^n$ stands for $H_{\left( 1,...,1\right) }^n$, $\left(
1,...,1\right) \in {\bf N}_{+}^{n+1}$ (the $n$-th ordinary de Rham
cohomology).

The first step in the chain (\ref{sediz}) is obtained via the following%
\footnote{%
This Lemma has been proved, independently, also by Yu. Torkhov.}:

\smallskip

\begin{lemma}
Let $n\in {\bf N}_{+}$. If $\sigma ,\tau \in {\bf N}_{+}^\infty $ are such
that $\sigma \left( n\right) =\tau \left( n\right) $, then:

(i) $\ker d_{\sigma \left( n+1\right) }=\ker d_{\tau \left( n+1\right) }$ ;

(ii) $im\left( d_{\sigma \left( n+1\right) }\right) \simeq im\left( d_{\tau
\left( n+1\right) }\right) $

(where $\simeq $ means $K-{\bf Mod}$-isomorphism).\ 
\end{lemma}

\TeXButton{Proof}{\proof} (ii) follows trivially from (i). Let $\sigma \in 
{\bf N}_{+}^n$ and $k>1$. Consider the short exact sequence:

\[
0\rightarrow {\rm D}_{\left( \sigma ,k-1,1\right) }\hookrightarrow {\rm D}%
_{\left( \sigma ,k-1\right) }^{\bullet }\circ {\rm Diff}_1^{\left( +\right) }%
\stackrel{i}{\longrightarrow }{\rm D}_{\left( \sigma ,k\right) }\rightarrow
0 
\]
whose dual-representative:

\begin{equation}
0\rightarrow {\bf \Lambda }^{\left( \sigma ,k\right) }\stackrel{i^{\vee }}{%
\longrightarrow }{\bf J}^1\left( {\bf \Lambda }^{\left( \sigma ,k-1\right)
}\right) \longrightarrow {\bf \Lambda }^{\left( \sigma ,k-1,1\right)
}\rightarrow 0  \label{garrafugio}
\end{equation}
is likewise exact (in ${\frak D}$). We embed the latter in the commutative
diagram:

\[
\begin{array}{ccccccccc}
0 & \rightarrow & 
\begin{array}{cc}
& 
\end{array}
{\bf \Lambda }^{\left( \sigma ,k\right) } & \stackrel{i^{\vee }}{%
\longrightarrow } & {\bf J}^1\left( {\bf \Lambda }^{\left( \sigma
,k-1\right) }\right) & \longrightarrow & {\bf \Lambda }^{\left( \sigma
,k-1,1\right) } & \rightarrow & 0 \\ 
&  & d_{\left( \sigma ,k\right) }\uparrow &  & \uparrow j_1 &  &  &  &  \\ 
&  & 
\begin{array}{cc}
& 
\end{array}
{\bf \Lambda }^\sigma & \stackunder{d_{\left( \sigma ,k-1\right) }}{%
\longrightarrow } & {\bf \Lambda }^{\left( \sigma ,k-1\right) } &  &  &  & 
\end{array}
\]

Now, if $\omega \in {\bf \Lambda }^\sigma $ is such that $d_{\left( \sigma
,k-1\right) }\left( \omega \right) =0$, then $j_1\left( d_{\left( \sigma
,k-1\right) }\left( \omega \right) \right) =0$ and, by commutativity, $%
\left( i^{\vee }\circ d_{\left( \sigma ,k\right) }\right) \left( \omega
\right) =0$. But $i^{\vee }$ is a monomorphism, so:

\[
\ker d_{\left( \sigma ,k-1\right) }\subseteq \ker d_{\left( \sigma ,k\right)
}\text{, }\forall k>1\text{.} 
\]
Since the inverse inclusion is obvious, (i) is proved.

\TeXButton{End Proof}{\endproof}\ 

To prove the ''$k$-th step'' of chain (\ref{sediz}), it is enough to show
that: 
\begin{equation}
H_{\left( \sigma \left( k-1\right) ,\sigma _k+1,1,\ldots ,1\right)
_{n+1}}^n\simeq H_{\left( \sigma \left( k-1\right) ,\sigma _k,1,\ldots
,1\right) _{n+1}}^n  \label{sizdiz}
\end{equation}
where we write $\left( \rho \right) _r$ if $\rho \in {\bf N}_{+}^r$.

To prove (\ref{sizdiz}) we construct an auxiliary complex.

\medskip\ 

Let $P$ be an object in${\frak D}$ and $\tau \in {\bf N}_{+}^\infty $. As
shown at the beginning of this section, smoothness of ${\frak D}$ implies
that $\forall n>0$ the ''relative'' sequence (\ref{relHol})

\begin{equation}
0\rightarrow {\cal P}_{\tau \left( n\right) }^{\bullet }\left[ P\right]
\longrightarrow \stackunder{}{\cal P}_{\tau \left( n-1\right) }^{\bullet
}\left[ P\right] \circ {\rm Diff}_{\tau _n}^{\left( +\right)
}\longrightarrow {\cal P}_{\left( \tau \left( n-2\right) ,\tau _{n-1}+\tau
_n\right) }^{\bullet }\left[ P\right] \rightarrow 0  \label{(18)}
\end{equation}
is exact also in the last term; hence, its dual representative:

\begin{equation}
0\rightarrow {\rm Hol}^{\left( \tau \left( n-2\right) ,\tau _{n-1}+\tau
_n\right) }\left[ P\right] \longrightarrow {\bf J}^{\tau _n}\left( {\rm Hol}%
^{\tau \left( n-1\right) }\left[ P\right] \right) \longrightarrow {\rm Hol}%
^{\tau \left( n\right) }\left[ P\right] \rightarrow 0  \label{dizet}
\end{equation}
is exact also in the first term. Furthermore, when $P$ varies in $Ob\left( 
{\frak D}\right) $, (\ref{dizet}) gives rise to a short exact sequence in $%
\left[ {\frak D},{\frak D}\right] $. We will refer to ${\rm Hol}^{\tau
\left( n\right) }\left[ P\right] $ as the ${\rm Hol}${\em -object of type }$%
\tau \left( n\right) $ of the $A$-module $P$; we have

\begin{equation}
{\rm Hol}^{\tau \left( n\right) }\left[ P\right] \simeq \frac{{\bf J}^{\tau
_n}\left( {\rm Hol}^{\tau \left( n-1\right) }\left[ P\right] \right) }{{\rm %
Hol}^{\left( \tau \left( n-2\right) ,\tau _{n-1}+\tau _n\right) }\left[
P\right] }  \label{dinoef}
\end{equation}
as $A$-modules. This allows us to give the following:

\begin{PropDef}
\TeXButton{plotigno}{\label{plotigno}}\label{p12} Let $P\in Ob\left( {\frak D%
}\right) $ and $\tau \in {\bf N}_{+}^\infty $. We define the sequence in $%
{\bf DIFF}_{A,{\frak D}}$: 
\[
\begin{array}{c}
{\bf Hol}^\tau \left( P\right) :\text{ }0\rightarrow {\rm Hol}^\emptyset
\left[ P\right] \doteq P\stackrel{\delta _{\tau \left( 1\right) }\left(
P\right) }{\longrightarrow }{\rm Hol}^{\tau \left( 1\right) }\left[ P\right]
\equiv {\bf J}^{\tau _1}\left( P\right) \rightarrow ... \\ 
\ldots \rightarrow {\rm Hol}^{\tau \left( n\right) }\left[ P\right] 
\stackrel{\delta _{\tau \left( n+1\right) }\left( P\right) }{\longrightarrow 
}{\rm Hol}^{\tau \left( n+1\right) }\left[ P\right] \rightarrow ...
\end{array}
\]
where, for each $n\geq 0$, 
\begin{equation}
\delta _{\tau \left( n+1\right) }\left[ P\right] :{\rm Hol}^{\tau \left(
n\right) }\left[ P\right] \longrightarrow {\rm Hol}^{\tau \left( n+1\right)
}\left[ P\right]  \label{(20)}
\end{equation}
is defined to be the DO whose description (\ref{ire}) of Remark \ref{photos}
is the canonical inclusion 
\[
{\cal P}_{\tau \left( n+1\right) }^{\bullet }\left[ P\right] \hookrightarrow 
\stackunder{}{\cal P}_{\tau \left( n\right) }^{\bullet }\left[ P\right]
\circ {\rm Diff}_{\tau _{n+1}}^{\left( +\right) }; 
\]
equivalently, $\delta _{\tau \left( n+1\right) }\left[ P\right] $ is the
composition:

\[
{\rm Hol}^{\tau \left( n\right) }\left[ P\right] \stackrel{j_{\tau
_{n+1}}\left( {\rm Hol}^{\tau \left( n\right) }\left[ P\right] \right) }{%
\longrightarrow }{\bf J}^{\tau _{n+1}}\left( {\rm Hol}^{\tau \left( n\right)
}\left[ P\right] \right) \stackrel{p_{\tau \left( n+1\right) }(P)}{%
\longrightarrow }\text{ }\frac{{\bf J}^{\tau _{n+1}}\left( {\rm Hol}^{\tau
\left( n\right) }\left[ P\right] \right) }{{\rm Hol}^{\left( \tau \left(
n-1\right) ,\tau _n+\tau _{n+1}\right) }\left[ P\right] }\simeq {\rm Hol}%
^{\tau \left( n+1\right) }\left[ P\right] 
\]
where $p_{\tau \left( n+1\right) }(P)$ is the canonical quotient projection. 
${\bf Hol}^\tau \left[ P\right] $ is a {\em complex} in ${\bf DIFF}_{A,%
{\frak D}}$, called ${\rm Hol}^\tau ${\em -complex of }$P$; moreover, ${\bf %
Hol}^\tau \left[ P\right] $ is natural in $P$ and defines a functor ${\bf Hol%
}^\tau :{\frak D}\rightarrow ${\bf K}$\left( {\bf DIFF}_{A,{\frak D}}\right) 
$\footnote{%
We recall that {\bf K}$\left( DIFF_{A,{\frak D}}\right) $ denotes the
category of complexes of differential operators formed by objects of ${\frak %
D}$.}.
\end{PropDef}

\TeXButton{Proof}{\proof} As always, it is better to work with functors
(i.e. differential operators) than with representative objects. In the
notations of Remark \ref{photos}, we have that $\varphi ^{\delta _{\tau
\left( n+1\right) }\circ \delta _{\tau \left( n\right) }}$\footnote{%
We write shortly $\delta _{\tau \left( k\right) }$ instead of $\delta _{\tau
\left( k\right) }\left[ P\right] $, for any $k\geq 0$.} coincides with the
composition:

\[
\begin{array}{cc}
&  \\ 
{\cal P}_{\tau \left( n+1\right) }^{\bullet }\left[ P\right] &  \\ 
\downarrow & \varphi ^{\delta _{\tau \left( n+1\right) }{}_{}} \\ 
{\cal P}_{\tau \left( n\right) }^{\bullet }\left[ P\right] \circ {\rm Diff}%
_{\tau _{n+1}}^{\left( +\right) } &  \\ 
\downarrow & \varphi ^{\delta _{\tau \left( n\right) }}\left( {\rm Diff}%
_{\tau _{n+1}}^{+}\right) \\ 
{\cal P}_{\tau \left( n-1\right) }^{\bullet }\left[ P\right] \circ \left( 
{\rm D}_{\left( \tau _n,\tau _{n+1}\right) }^{\bullet },{\rm D}_{\left( \tau
_n,\tau _{n+1}\right) }^{+}\right) &  \\ 
\downarrow & {\cal P}_{\tau \left( n-1\right) }^{\bullet }\left[ P\right]
\left( {\em C}_{\tau _n,\tau _{n+1}}\right) \\ 
{\cal P}_{\tau \left( n-1\right) }^{\bullet }\left[ P\right] \circ {\rm Diff}%
_{\tau _{n+1}+\tau _n}^{\left( +\right) } & 
\end{array}
\]
where the first two arrows are monomorphisms and the last is the ''gluing''
morphism with respect to the indexes $\left( \tau _n,\tau _{n+1}\right) $.
This composition is zero. In fact, if $Q\in Ob\left( {\frak D}\right) $ and $%
\Delta \in \left[ {\cal P}_{\tau \left( n+1\right) }^{\bullet }\left[
P\right] \right] \left( Q\right) $, then the image $\overline{\Delta }$ of $%
\Delta $ via this composition, is defined by:

\[
\left( \overline{\Delta }\left( p\right) \right) \left( a_1\right) \cdots
\left( a_{n-1}\right) \doteq \left( \left( \Delta \left( p\right) \right)
\left( a_1\right) \cdots \left( a_{n-1}\right) \right) \left( 1\right) \text{%
,} 
\]
and is zero because $\left( \Delta \left( p\right) \right) \left( a_1\right)
\cdots \left( a_{n-1}\right) \in {\rm D}_{\tau _{n+1}}\left( Q\right) $, for
each $p\in P$, $a_1,\ldots ,a_{n-1}\in A$.\qquad \TeXButton{End Proof}
{\endproof}

Now we show that if $\tau \in {\bf N}_{+}^\infty $ is {\it regular}, then,
for any object $P$ in ${\frak D}$, ${\bf Hol}^\tau \left[ P\right] $ is {\it %
acyclic}. In order to do this, we will exhibit (functorially in $P$) a {\it %
trivializing homotopy}.

Define: 
\begin{eqnarray*}
\varphi _\emptyset \left( P\right) &:&{\cal P}_{\tau \left( -1\right)
}^{\bullet }\left[ P\right] \doteq 0\longrightarrow {\cal P}_{\tau \left(
0\right) }^{\bullet }\left[ P\right] \doteq Hom_A\left( P,\cdot \right) \\
\varphi _{\tau \left( 1\right) }\left( P\right) &:&{\cal P}_{\tau \left(
0\right) }^{\bullet }\left[ P\right] \doteq Hom_A\left( P,\cdot \right)
\hookrightarrow {\cal P}_{\tau \left( 1\right) }^{\bullet }\left[ P\right]
\equiv {\rm Diff}_{\tau _1}\left( P,\cdot \right)
\end{eqnarray*}
which are morphisms in $\left[ {\frak D},{\frak D}\right] $; then define, by
induction on $n$,

\[
\varphi ^{\tau \left( n+1\right) }\left( P\right) \doteq i_{\tau \left(
n+1\right) }\left( P\right) -\widehat{\varphi }_{\tau \left( n\right)
}\left( P\right) \text{,} 
\]
where

\[
i_{\tau \left( n+1\right) }\left( P\right) :{\cal P}_{\tau \left( n\right)
}^{\bullet }\left[ P\right] \simeq {\cal P}_{\tau \left( n\right) }^{\bullet
}\left[ P\right] \circ {\rm Diff}_0^{\left( +\right) }\hookrightarrow {\cal P%
}_{\tau \left( n\right) }^{\bullet }\left[ P\right] \circ {\rm Diff}_{\tau
_{n+1}}^{\left( +\right) } 
\]

\begin{eqnarray*}
\widehat{\varphi }_{\tau \left( n\right) }\left( P\right) &:&{\cal P}_{\tau
\left( n\right) }^{\bullet }\left[ P\right] \hookrightarrow {\cal P}_{\tau
\left( n-1\right) }^{\bullet }\left[ P\right] \circ {\rm Diff}_{\tau
_n}^{\left( +\right) }\stackrel{\varphi _{\tau \left( n+1\right) }\left(
P\right) \left( {\rm Diff}_{\tau _n}^{+}\right) }{\longrightarrow } \\
\ &\rightarrow &{\cal P}_{\tau \left( n\right) }^{\bullet }\left[ P\right]
\circ {\rm Diff}_{\tau _n}^{\left( +\right) }\hookrightarrow {\cal P}_{\tau
\left( n\right) }^{\bullet }\left[ P\right] \circ {\rm Diff}_{\tau
_{n+1}}^{\left( +\right) }.
\end{eqnarray*}
With this definition, $\varphi _{\tau \left( n+1\right) }\left( P\right) :%
{\cal P}_{\tau \left( n\right) }^{\bullet }\left[ P\right] \rightarrow {\cal %
P}_{\tau \left( n\right) }^{\bullet }\left[ P\right] \circ {\rm Diff}_{\tau
_{n+1}}^{\left( +\right) }$, but it is easy to resolve the inductive
definition in the following one:

\begin{eqnarray}
\left\{ \left[ \varphi _{\tau \left( n+1\right) }\left( P\right) \left(
Q\right) \right] \left( \Delta \right) \right\} \left( p\right) \left(
a_1\right) \cdots \left( a_n\right) &=&a_n\Delta \left( p\right) \left(
a_1\right) \cdots \left( a_{n-1}\right) +  \label{lilla} \\
&&\ \sum\limits_{k=1}^{n-1}\left( -1\right) ^{n-k}\Delta \left( p\right)
\left( a_1\right) \cdots \left( a_ka_{k+1}\right) \cdots \left( a_n\right) +
\nonumber \\
&&\ +\left( -1\right) ^n\Delta \left( a_1p\right) \left( a_2\right) \cdots
\left( a_n\right)  \nonumber
\end{eqnarray}
($Q\in Ob\left( {\frak D}\right) $, $p\in P$, $a_1,\ldots ,a_n\in A$ and $%
\Delta \in {\cal P}_{\tau \left( n\right) }^{\bullet }\left[ P\right] \left(
Q\right) $). This shows that actually 
\[
\varphi _{\tau \left( n+1\right) }\left( P\right) :{\cal P}_{\tau \left(
n\right) }^{\bullet }\left[ P\right] \rightarrow {\cal P}_{\tau \left(
n+1\right) }^{\bullet }\left[ P\right] . 
\]
Therefore, we get a family

\[
\left\{ \varphi _{\tau \left( n\right) }\left( P\right) :{\cal P}_{\tau
\left( n-1\right) }^{\bullet }\left[ P\right] \rightarrow {\cal P}_{\tau
\left( n\right) }^{\bullet }\left[ P\right] \right\} _{n>0} 
\]
of morphisms in $\left[ {\frak D},{\frak D}\right] $. Of course, formula (%
\ref{lilla}) can equally be taken as the definition of the family $\left\{
\varphi _{\tau \left( n\right) }\left( P\right) \right\} _{n>0}$ but the
inductive definition can be ''dualized'', to representative objects, to give
the following (keeping the notations of Proposition \ref{plotigno}):

\begin{eqnarray*}
\varphi ^\emptyset \left( P\right) &:&{\rm Hol}^{\tau \left( 0\right)
}\left[ P\right] \doteq P\longrightarrow {\rm Hol}^{\tau \left( -1\right)
}\left[ P\right] \doteq 0 \\
\varphi ^{\tau \left( 1\right) }\left( P\right) &:&{\rm Hol}^{\tau \left(
1\right) }\left( P\right) \doteq {\bf J}^{\tau _1}\left( P\right)
\longrightarrow {\rm Hol}^{\tau \left( 0\right) }\left[ P\right] \doteq P%
\text{ (natural projection)} \\
\varphi ^{\tau \left( n+1\right) }\left( P\right) &:&{\rm Hol}^{\tau \left(
n+1\right) }\left[ P\right] \longrightarrow {\rm Hol}^{\tau \left( n\right)
}\left[ P\right]
\end{eqnarray*}
where $\varphi ^{\tau \left( n+1\right) }\left( P\right) $ is the only $%
{\frak D}$-morphism corresponding to the DO of order $\leq \tau _{n+1}$

\[
\Delta _{\tau \left( n\right) }\doteq {\rm id}_{{\rm Hol}^{\tau \left(
n\right) }\left[ P\right] }-\delta _{\tau \left( n\right) }\left( P\right)
\circ \varphi ^{\tau \left( n\right) }\left( P\right) :{\rm Hol}^{\tau
\left( n\right) }\left[ P\right] \longrightarrow {\rm Hol}^{\tau \left(
n\right) }\left[ P\right] \text{.} 
\]
As before, but dually\footnote{%
Subfunctors of strictly representable functors correspond to quotient
objects of the representatives.}, this definition gives apparently a ${\frak %
D}$-morphism:

\[
\widehat{\varphi }^{\tau \left( n+1\right) }\left( P\right) :{\bf J}^{\tau
_{n+1}}\left( {\rm Hol}^{\tau \left( n\right) }\left[ P\right] \right)
\longrightarrow {\rm Hol}^{\tau \left( n\right) }\left[ P\right] 
\]
($\tau $ being regular) but formula (\ref{lilla}) shows that actually

\[
\ker \left( \widehat{\varphi }^{\tau \left( n+1\right) }\left( P\right)
\right) \supseteq {\rm Hol}^{\left( \tau \left( n-1\right) ,\tau _n+\tau
_{n+1}\right) }\left[ P\right] 
\]
so that, by formula (\ref{dinoef}), $\widehat{\varphi }^{\tau \left(
n+1\right) }\left( P\right) $ induces, by passing to the quotient, the
morphism $\varphi ^{\tau \left( n+1\right) }\left( P\right) $ we wanted.

Now we have a family of ${\frak D}$-morphisms $\left\{ \varphi ^{\tau \left(
n\right) }\left( P\right) :{\rm Hol}^{\tau \left( n\right) }\left[ P\right]
\rightarrow {\rm Hol}^{\tau \left( n-1\right) }\left[ P\right] \right\}
_{n\geq 0}$ dual to $\left\{ \varphi _{\tau \left( n\right) }\left( P\right)
:{\cal P}_{\tau \left( n\right) }^{\bullet }\left[ P\right] \rightarrow 
{\cal P}_{\tau \left( n-1\right) }^{\bullet }\left[ P\right] \right\}
_{n\geq 0}$.

\begin{proposition}
\TeXButton{dan}{\label{dan}}For each object $P$ in ${\frak D}$ and for each
regular $\tau \in {\bf N}_{+}^\infty $, $\left\{ \varphi ^{\tau \left(
n\right) }\left( P\right) \right\} _{n\geq 0}$ is a trivializing homotopy
for ${\bf Hol}^\tau \left( P\right) $. Furthermore, $\left\{ \varphi ^{\tau
\left( n\right) }\left( P\right) \right\} _{n\geq 0}$ is natural in $P$.
\end{proposition}

\TeXButton{Proof}{\proof} We must show that the {\it sum} $L+R$ of the two
compositions:

\begin{eqnarray*}
L &:&{\cal P}_{\tau \left( n\right) }^{\bullet }\left[ P\right]
\hookrightarrow \stackunder{}{\cal P}_{\tau \left( n-1\right) }^{\bullet
}\left[ P\right] \circ {\rm Diff}_{\tau _n}^{\left( +\right) }\stackrel{%
\varphi _{\tau \left( n\right) }\left( P\right) \left( {\rm Diff}_{\tau
_n}^{+}\right) }{\longrightarrow } \\
&\rightarrow &\stackunder{}{\cal P}_{\tau \left( n\right) }^{\bullet }\left[
P\right] \circ {\rm Diff}_{\tau _n}^{\left( +\right) }\hookrightarrow 
\stackunder{}{\cal P}_{\tau \left( n\right) }^{\bullet }\left[ P\right]
\circ {\rm Diff}_{\tau _{n+1}}^{\left( +\right) }
\end{eqnarray*}

\[
R:{\cal P}_{\tau \left( n\right) }^{\bullet }\left[ P\right] \stackrel{%
\varphi _{\tau \left( n+1\right) }\left( P\right) }{\longrightarrow }{\cal P}%
_{\tau \left( n+1\right) }^{\bullet }\left[ P\right] \hookrightarrow 
\stackunder{}{\cal P}_{\tau \left( n\right) }^{\bullet }\left[ P\right]
\circ {\rm Diff}_{\tau _{n+1}}^{\left( +\right) } 
\]
equals ${\rm id}_{{\cal P}_{\tau \left( n\right) }^{\bullet }\left[ P\right]
}$ (which is then homotopic to the zero map) or, equivalently, that the
diagram:

\begin{equation}
\begin{tabular}{lll}
$\quad {\cal P}_{\tau \left( n\right) }^{\bullet }\left[ P\right] $ & $%
\stackrel{L+R}{\longrightarrow }$ & $\stackunder{}{\cal P}_{\tau \left(
n\right) }^{\bullet }\left[ P\right] \circ {\rm Diff}_{\tau _{n+1}}^{\left(
+\right) }$ \\ 
$\qquad \quad \quad {\rm id}$ & $\searrow $ & $\quad \quad \cup $ \\ 
&  & $\stackunder{}{\cal P}_{\tau \left( n\right) }^{\bullet }\left[
P\right] \circ {\rm Diff}_0^{\left( +\right) }\quad $%
\end{tabular}
\label{tomaso}
\end{equation}
is commutative. For $Q\in Ob\left( {\frak D}\right) $ and $\Delta \in {\cal P%
}_{\tau \left( n\right) }^{\bullet }\left[ P\right] \left( Q\right) $, we
have by (\ref{lilla}):

\[
L\left( \Delta \right) \left( p\right) \left( a_1\right) \cdots \left(
a_n\right) = 
\]
\[
=\left[ a_{n-1}^{+}\Delta \left( p\right) \left( a_1\right) \cdots \left(
a_{n-2}\right) +\sum\limits_{s=1}^{n-2}\left( -1\right) ^{n-1-s}\Delta
\left( p\right) \left( a_1\right) \cdots \left( a_sa_{s+1}\right) \cdots
\left( a_{n-1}\right) +\right. 
\]
\[
\left. +\left( -1\right) ^{n-1}\Delta \left( a_1p\right) \left( a_2\right)
\cdots \left( a_{n-1}\right) \right] \left( a_n\right)
=-\sum\limits_{s=1}^{n-1}\left( -1\right) ^{n-s}\Delta \left( p\right)
\left( a_1\right) \cdots \left( a_sa_{s+1}\right) \cdots \left( a_n\right) + 
\]

\[
+\left( -1\right) _{}^{n-1}\Delta \left( a_1p\right) \left( a_2\right)
\cdots \left( a_n\right) 
\]
while

\begin{eqnarray*}
R\left( \Delta \right) \left( p\right) \left( a_1\right) \cdots \left(
a_n\right) &=&a_n\Delta \left( p\right) \left( a_1\right) \cdots \left(
a_{n-1}\right) +  \label{lilla} \\
&&\ \sum\limits_{k=1}^{n-1}\left( -1\right) ^{n-k}\Delta \left( p\right)
\left( a_1\right) \cdots \left( a_ka_{k+1}\right) \cdots \left( a_n\right) +
\\
&&\ +\left( -1\right) ^n\Delta \left( a_1p\right) \left( a_2\right) \cdots
\left( a_n\right)
\end{eqnarray*}
so that

\[
\left( L+R\right) \left( \Delta \right) \left( p\right) \left( a_1\right)
\cdots \left( a_n\right) =a_n\Delta \left( p\right) \left( a_1\right) \cdots
\left( a_{n-1}\right) 
\]
i.e. $\left( L+R\right) \left( \Delta \right) $ coincides with the image of $%
\Delta $ via the inclusion

\[
{\cal P}_{\tau \left( n\right) }^{\bullet }\left[ P\right] \left( Q\right)
\simeq \stackunder{}{\cal P}_{\tau \left( n\right) }^{\bullet }\left[
P\right] \circ {\rm Diff}_0^{(+)}\left( Q\right) \hookrightarrow \stackunder{%
}{\cal P}_{\tau \left( n\right) }^{\bullet }\left[ P\right] \circ {\rm Diff}%
_{\tau _{n+1}}^{(+)}\left( Q\right) \text{.} 
\]
\TeXButton{End Proof}{\endproof}

\medskip\ 

We now use acyclicity of the ${\bf Hol}^{{\bf 1}}$-complex,${\bf 1}%
=(1,...,1,1,...,1,...)\in {\bf N}_{+}^\infty $, to prove the ''$k$-th step''
i.e. formula (\ref{sizdiz}). Let $\sigma =(\sigma _1,...,\sigma _k+1)\in 
{\bf N}_{+}^k$, $(\sigma ,{\bf 1})=(\sigma _1,...,\sigma
_k+1,1,1,...,1,...)\in {\bf N}_{+}^\infty $ and 
\[
{\bf K}_{(\sigma ,{\bf 1})}^{\left( k\right) }\doteq \ker \left( {\bf dR}%
_{(\sigma ,{\bf 1})}\rightarrow {\bf dR}_{(\sigma _1,...,\sigma
_k,1,1,...,1,...)}\right) . 
\]
For each $\left( \mu \right) _s=$ $\left( \mu _1,...,\mu _s\right) \in {\bf N%
}_{+}^s$, $1\leq r\leq s$, $r,s\in {\bf N}_{+}$, we put:

\[
K_{(\mu )_s}^{\left( r\right) }\doteq \ker \left( {\bf \Lambda }^{\left( \mu
_1,...,\mu _r,...\mu _s\right) }\rightarrow {\bf \Lambda }^{\left( \mu
_1,...,\mu _r-1,...\mu _s\right) }\right) \text{.} 
\]
To prove the ''$k$-th step'' it is enough to show acyclicity of ${\bf K}%
_{(\sigma ,{\bf 1})}^{\left( k\right) }$. We claim that there exists a
resolution of ${\bf K}_{(\sigma ,{\bf 1})}^{\left( k\right) }$ of the form: 
\[
\cdots \rightarrow {\bf Hol}^{{\bf 1}}\left[ K_{\left( \sigma _1,...,\sigma
_k+l+1\right) }^{(k)}\right] \left[ -k-l\right] \stackrel{\psi _l\left[
-k+l\right] }{\longrightarrow }{\bf Hol}^{{\bf 1}}\left[ K_{\left( \sigma
_1,...,\sigma _k+l\right) }^{(k)}\right] \left[ -k-l+1\right] \rightarrow
\cdots 
\]
\begin{equation}
\cdots \stackunder{}{\stackunder{\psi _2\left[ -\left( k+2\right) \right] }{%
\rightarrow }}{\bf Hol}^{{\bf 1}}\left[ K_{\left( \sigma _1,...,\sigma
_k+2\right) }^{(k)}\right] \left[ -k-1\right] \stackunder{\psi _1\left[
-k+1\right] }{\longrightarrow }{\bf Hol}^{{\bf 1}}\left[ K_{\left( \sigma
_1,...,\sigma _k+1\right) }^{(k)}\right] \left[ -k\right] \stackunder{\rho }{%
\longrightarrow }{\bf K}_{(\sigma ,{\bf 1})}^{\left( k\right) }\rightarrow 0
\label{res}
\end{equation}
where if $r\in {\bf N}_{+}$, $\left( \cdot \right) \left[ r\right] $
denotes, as usual, the $r$-shift both for complexes and morphisms of
complexes. We postpone in the Appendix the definition of the maps of
complexes 
\begin{eqnarray*}
\psi _l &:&{\bf Hol}^{{\bf 1}}\left[ K_{\left( \sigma _1,...,\sigma
_k+l+1\right) }^{(k)}\right] \stackrel{}{\longrightarrow }{\bf Hol}^{{\bf 1}%
}\left[ K_{\left( \sigma _1,...,\sigma _k+l\right) }^{(k)}\right] \left[
1\right] ,\text{ \quad }l\in {\bf N}_{+}, \\
\rho &:&{\bf Hol}^{{\bf 1}}\left[ K_{\left( \sigma _1,...,\sigma _k+1\right)
}^{(k)}\right] \left[ -k\right] \stackunder{}{\longrightarrow }{\bf K}%
_{(\sigma ,{\bf 1})}^{\left( k\right) }
\end{eqnarray*}
and the proof that (\ref{res}) is actually a resolution.

Assuming the existence of resolution (\ref{res}), the acyclicity of ${\bf K}%
_{(\sigma ,{\bf 1})}^{\left( k\right) }$ is then an immediate consequence of
acyclicity of ${\bf Hol}^{{\bf 1}}$-complexes together with the following
elementary fact

\begin{lemma}
Let $C^{\cdot },$ $P_i^{\cdot },$ $i>0,$ be cochain complexes in $A-{\bf Mod}
$ and 
\[
\cdots \longrightarrow P_n^{\cdot }\longrightarrow P_{n-1}^{\cdot
}\longrightarrow \cdots \longrightarrow P_1^{\cdot }\longrightarrow C^{\cdot
}\rightarrow 0 
\]
be a resolution of $C^{\cdot }.$ Suppose that $\forall i\geq 1$, $P_i^k=(0)$ 
$\forall k<0$ (so that $C^k=(0)$ $\forall k<0$ too). If each $P_i^{\cdot }$
is acyclic then so is $C^{\cdot }$.
\end{lemma}

\TeXButton{Proof}{\proof} It follows from the hypotheses that $C^{\cdot }$
is isomorphic in the derived category $D^{+}\left( A-{\bf Mod}\right) $ to
the total complex associated to the double complex induced by 
\[
\cdots \longrightarrow P_n^{\cdot }\longrightarrow P_{n-1}^{\cdot
}\longrightarrow \cdots \longrightarrow P_1^{\cdot }\text{ } 
\]
which is acyclic.

For those readers who feel uncomfortable with derived categories, here is a
more ''step-by-step'' proof. By the usual sign-trick we can associate to the
given resolution a $1^{st}$- quadrant double complex $R_{\cdot }^{\cdot
}=\left( R_q^p\right) $ which is mixed: homological in the vertical (i.e.
with $p$ fixed) direction and cohomological in the horizontal (i.e. with $p$
fixed) direction. We turn it into a 2$^{nd}$-quadrant homological double
complex $\widehat{R}_{\cdot \text{ }\cdot }=\left( \widehat{R}_{pq}\right) $
with $\widehat{R}_{pq}\doteq R_q^{-p}$.

Consider the spectral sequence induced by the ''filtration by rows'' on $%
\widehat{R}_{\cdot \text{ }\cdot }$ (e.g. \cite{We} p. 142):

\[
^{II}E_{pq}^0=\widehat{R}_{qp} 
\]
with $^{II}d_{pq}^0$ given by the horizontal differential in $\widehat{R}%
_{\cdot \text{ }\cdot }$ . Then $\left\{ ^{II}E_{pq}^1,^{II}d_{pq}^1\right\} 
$ is just:

\begin{eqnarray*}
&& 
\begin{tabular}{l}
$q=n+1$ \\ 
$q=n$ \\ 
$\quad \vdots $ \\ 
$q=2$ \\ 
$q=1$ \\ 
$q=0$%
\end{tabular}
\begin{tabular}{|l|l|l|l|l|l|}
\hline
$\leftarrow 0$ & $\leftarrow H^{n+1}\left( C^{\cdot }\right) $ & $\leftarrow
0$ & $\leftarrow 0$ &  & $\leftarrow 0$ \\ \hline
$\leftarrow 0$ & $\leftarrow H^n\left( C^{\cdot }\right) $ & $\leftarrow 0$
& $\leftarrow 0$ &  & $\leftarrow 0$ \\ \hline
$\quad \vdots $ & $\qquad \vdots $ & \quad $\vdots $ & \quad $\vdots $ &  & 
\quad $\vdots $ \\ \hline
$\leftarrow 0$ & $\leftarrow H^2\left( C^{\cdot }\right) $ & $\leftarrow 0$
& $\leftarrow 0$ &  & $\leftarrow 0$ \\ \hline
$\leftarrow 0$ & $\leftarrow H^1\left( C^{\cdot }\right) $ & $\leftarrow 0$
& $\leftarrow 0$ &  & $\leftarrow 0$ \\ \hline
$\leftarrow 0$ & $\leftarrow H^0(C^{\cdot })$ & $\leftarrow 0$ & $\leftarrow
0$ &  & $\leftarrow 0$ \\ \hline
\end{tabular}
\\
&& 
\begin{tabular}{lllllll}
$\quad \quad \quad p=-1\quad $ & $\quad 0\quad $ & $\quad \qquad 1\quad $ & $%
\quad 2$ & $\cdots $ & $\quad n$ & 
\end{tabular}
\end{eqnarray*}
(with differential induced by vertical differential in $\widehat{R}_{\cdot 
\text{ }\cdot }$). So $^{II}E$ degenerates at $^{II}E^1$. But $\widehat{R}%
_{\cdot \text{ }\cdot }$ is a 2$^{nd}$-quadrant double complex hence $^{II}E$
converges to $H^{*}\left( Tot\left( \widehat{R}_{\cdot \text{ }\cdot
}\right) \right) $, $Tot\left( \widehat{R}_{\cdot \text{ }\cdot }\right) $
being the total complex associated to $\widehat{R}_{\cdot \text{ }\cdot }$,
which is zero since by hypothesis $\widehat{R}_{\cdot \text{ }\cdot }$ has
exact columns. So $^{II}E_{pq}^1=(0)$ $\forall p,q$ and $C^{\cdot }$ is
acyclic. \TeXButton{End Proof}{\endproof}

\begin{corollary}
\TeXButton{ty}{\label{ty}}Let ${\frak D}\subseteq A-{\bf Mod}$ be a
differentially closed smooth subcategory. If we define the {\em stable}
infinite de Rham complex to be 
\[
{\bf dR}_\infty ^{\text{st}}\doteq \underleftarrow{\lim }_{k>0}{\bf dR}_{%
{\bf k}} 
\]
(where ${\bf k}\doteq \left( k,...,k,k,...k,...\right) $ ) then the
canonical morphism ${\bf dR}_\infty ^{\text{st}}\longrightarrow {\bf dR}_{%
{\bf 1}}$ is a quasi-isomorphism.
\end{corollary}

\TeXButton{Proof}{\proof} We use the following facts:

(i) the index category of the inverse system which defines ${\bf dR}_\infty
^{\text{st}}$ is countable;

(ii) the canonical ${\frak D}$-morphisms ${\bf dR}_{\underline{k^{\prime }}%
}\rightarrow {\bf dR}_{\underline{k}}$, $k^{\prime }\geq k$, in the inverse
system are epimorphisms.

If we denote by $\underleftarrow{\lim }_k^1$ the first right derived functor
of $\underleftarrow{\lim }_k$, (i) and (ii) imply, via standard spectral
sequence's arguments (e.g. \cite{Lu} Cor. 1.1, p. 535), that there is a
short exact sequence

\[
0\rightarrow \underleftarrow{\lim }_{k>0}^1\left[ H^{n-1}\left( {\bf dR}_{%
{\bf k}}\right) \right] \rightarrow H^n\left( {\bf dR}_\infty ^{\text{st}%
}\right) \rightarrow \underleftarrow{\lim }_{k>0}\left[ H^n\left( {\bf dR}_{%
{\bf k}}\right) \right] \rightarrow 0. 
\]
By theorem \ref{main}, the term on the right is isomorphic to $H_{dR}^n$, so
we are left to prove that $\underleftarrow{\lim }_{k>0}\left[ H^{n-1}\left( 
{\bf dR}_{\underline{k}}\right) \right] =\left( 0\right) $. But $%
\underleftarrow{\lim }_{k>0}^1$is right exact and

\[
H^{n-1}\left( {\bf dR}_{{\bf k}}\right) \equiv H_{{\bf k}}^{n-1}\cong
H_{dR}^{n-1}\doteq H_{{\bf 1}}^{n-1}\text{, }\forall n\geq 1\text{, } 
\]
by theorem \ref{main}, therefore it will be enough to prove the vanishing of 
$\underleftarrow{\lim }_{k>0}^1$ for the constant inverse system 
\[
\cdots \rightarrow H_{dR}^{n-1}\stackrel{id}{\longrightarrow }H_{dR}^{n-1}%
\stackrel{id}{\longrightarrow }H_{dR}^{n-1}\rightarrow \cdots 
\]

But this is an easy consequence of the following description (due to
Eilenberg) of $\underleftarrow{\lim }_{k>0}^1$for constant systems.

If we define

\[
D_0:\prod\limits_{k\in {\bf N}_{+}}H_{dR}^{n-1}\longrightarrow
\prod\limits_{k\in {\bf N}_{+}}H_{dR}^{n-1}:\left( \alpha _k\right) _{k\in 
{\bf N}_{+}}\longmapsto \left( \alpha _{k+1}-\alpha _k\right) _{k\in {\bf N}%
_{+}}\text{ ;} 
\]
then

\[
{\rm co}\ker \left( D_0\right) =\underleftarrow{\lim }_{k>0}^1\left[
H_{dR}^{n-1}\right] \text{.}\label{4.30} 
\]
Let $\left( \omega _k\right) _{k>0}\in \prod\limits_{k\in {\bf N}%
_{+}}H_{dR}^{n-1}$ and define $\left( \overline{\omega }_k\right) _{k>0}$ as 
$\overline{\omega }_k\doteq \sum\limits_{l=1}^{k-1}\omega _l$; then 
\[
D_0\left( \left( \overline{\omega }_k\right) _{k>0}\right) =\left( \overline{%
\omega }_{k+1}-\overline{\omega }_k\right) _{k>0}=\left( \omega _k\right)
_{k>0}. 
\]
Therefore $D_0$ is surjective and we conclude. \TeXButton{End Proof}
{\endproof}

\begin{corollary}
\label{markandeja} Let ${\frak D}\subseteq A-{\bf Mod}$ be a differentially
closed smooth subcategory such that $\forall \sigma \in {\bf N}_{+}^\infty $%
, $\exists n_\sigma \in {\bf N}_{+}$ :

\[
{\bf \Lambda }^{\sigma \left( r\right) }=(0)\text{, }\forall r>n_\sigma 
\text{.} 
\]
Then the canonical morphism ( Def. \ref{infty})

\[
{\bf dR}_\infty \rightarrow {\bf dR}_\sigma 
\]
is a quasi-isomorphism $\forall \sigma \in {\bf N}_{+}^\infty $.
\end{corollary}

\TeXButton{Proof}{\proof} Under our hypotheses

\[
\widehat{{\bf N}}_{+}\doteq \bigsqcup\limits_{k\in {\bf N}_{+}}\left\{ 
\underline{k}\in {\bf N}_{+}^\infty \mid \underline{k}\left( n\right) \equiv
\left( k,...,k\right) \in {\bf N}_{+}^n\text{, }\forall n\in {\bf N}%
_{+}\right\} 
\]
is cofinal in the index category of the system $\left\{ {\bf dR}_\sigma
\right\} $; hence ${\bf dR}_\infty \simeq \underleftarrow{\lim }_{k>0}{\bf dR%
}_{{\bf k}}$ and the thesis follows from corollary \ref{ty}. 
\TeXButton{End Proof}{\endproof}

\section{Appendix}

This appendix is devoted to defining the maps in the sequence (\ref{res})
and to showing that (\ref{res}) is exact.

There are two kinds of maps of complexes to be defined: 
\[
\rho \equiv \left( \rho ^n:{\bf Hol}^{{\bf 1}_{n-k}}\left[ K_{\left( \sigma
_1,...,\sigma _k+1\right) }^{\left( k\right) }\right] \rightarrow K_{\left(
\sigma _1,...,\sigma _k+1,1,...,1\right) _n}^{\left( k\right) }\right)
_{n\geq 0} 
\]
\[
\psi _l\equiv \left( \psi _l^n\right) _{n\geq 0}:{\bf Hol}^{{\bf 1}}\left[
K_{\left( \sigma _1,...,\sigma _k+l+1\right) }^{\left( k\right) }\right]
\rightarrow {\bf Hol}^{{\bf 1}}\left[ K_{\left( \sigma _1,...,\sigma
_k+l\right) }^{\left( k\right) }\right] \left[ 1\right] . 
\]

Let us first define $\rho $. We will define a functorial morphism 
\[
\Theta :{\rm D}_{\left( \sigma _1,...,\sigma _k+1,1,...,1\right)
_n}\longrightarrow {\cal P}_{{\bf 1}_{n-k}}^{\bullet }\left[ K_{\left(
\sigma _1,...,\sigma _k+1\right) }^{\left( k\right) }\right] 
\]
and show that the sequence 
\[
0\rightarrow {\rm D}_{\left( \sigma _1,...,\sigma _k,1,...,1\right) _n}%
\stackrel{}{\stackrel{\Phi }{\hookrightarrow }}{\rm D}_{\left( \sigma
_1,...,\sigma _k+1,1,...,1\right) _n}\stackrel{\Theta }{\longrightarrow }%
{\cal P}_{{\bf 1}_{n-k}}^{\bullet }\left[ K_{\left( \sigma _1,...,\sigma
_k+1\right) }^{\left( k\right) }\right] 
\]
is exact so that the dual representative of $\Theta $ will pass to the
quotient defining our surjective $\rho $.

Define $\Theta $ to be the following composition: 
\begin{eqnarray*}
{\rm D}_{\left( \sigma _1,...,\sigma _k+1,1,...,1\right) _n}\simeq {\rm D}%
_{\left( \sigma _1,...,\sigma _k+1,1\right) _{k+1}}^{\bullet }\left( {\rm D}%
_{{\bf 1}_{n-k-1}}\subset {\rm Diff}_{{\bf 1}_{n-k-1}}^{+}\right) \simeq
\end{eqnarray*}
\begin{eqnarray*}
\ &\simeq &Hom_A^{\bullet }\left( {\bf \Lambda }^{\left( \sigma
_1,...,\sigma _k+1,1\right) },{\rm D}_{{\bf 1}_{n-k-1}}\subset {\rm Diff}_{%
{\bf 1}_{n-k-1}}^{+}\right) \stackrel{\circ \text{ }d_{\left( \sigma
_1,...,\sigma _k+1,1\right) }}{\longrightarrow } \\
\ &\longrightarrow &{\rm Diff}_1^{\bullet }\left( {\bf \Lambda }^{\left(
\sigma _1,...,\sigma _k+1\right) },{\rm D}_{{\bf 1}_{n-k-1}}\subset {\rm Diff%
}_{{\bf 1}_{n-k-1}}^{+}\right) \simeq
\end{eqnarray*}
\[
\simeq {\cal P}_{{\bf 1}_{n-k}}^{\bullet }\left[ {\bf \Lambda }^{\left(
\sigma _1,...,\sigma _k+1\right) }\right] \stackrel{{\cal P}_{{\bf 1}%
_{n-k}}^{\bullet }\left[ j\right] }{\longrightarrow }{\cal P}_{{\bf 1}%
_{n-k}}^{\bullet }\left[ K_{\left( \sigma _1,...,\sigma _k+1\right)
}^{\left( k\right) }\right] 
\]
where 
\[
0\rightarrow K_{\left( \sigma _1,...,\sigma _k+1\right) }^{\left( k\right) }%
\stackrel{j}{\longrightarrow }{\bf \Lambda }^{\left( \sigma _1,...,\sigma
_k+1\right) }\stackrel{\pi }{\longrightarrow }{\bf \Lambda }^{\left( \sigma
_1,...,\sigma _k\right) }\rightarrow 0. 
\]
Then, $\Phi \circ \Theta $ coincides with the following 
\begin{eqnarray*}
{\rm D}_{\left( \sigma _1,...,\sigma _k,1,...,1\right) _n}\simeq {\rm D}%
_{\left( \sigma _1,...,\sigma _k,1\right) _{k+1}}^{\bullet }\left( {\rm D}_{%
{\bf 1}_{n-k-1}}\subset {\rm Diff}_{{\bf 1}_{n-k-1}}^{+}\right) \simeq
\end{eqnarray*}
\begin{eqnarray*}
\ &\simeq &Hom_A^{\bullet }\left( {\bf \Lambda }^{\left( \sigma
_1,...,\sigma _k,1\right) },{\rm D}_{{\bf 1}_{n-k-1}}\subset {\rm Diff}_{%
{\bf 1}_{n-k-1}}^{+}\right) \stackrel{\circ \text{ }d_{\left( \sigma
_1,...,\sigma _k,1\right) }}{\longrightarrow } \\
\ &\longrightarrow &{\rm Diff}_1^{\bullet }\left( {\bf \Lambda }^{\left(
\sigma _1,...,\sigma _k\right) },{\rm D}_{{\bf 1}_{n-k-1}}\subset {\rm Diff}%
_{{\bf 1}_{n-k-1}}^{+}\right) \simeq
\end{eqnarray*}
\[
\simeq {\cal P}_{{\bf 1}_{n-k}}^{\bullet }\left[ {\bf \Lambda }^{\left(
\sigma _1,...,\sigma _k\right) }\right] \stackrel{{\cal P}_{{\bf 1}%
_{n-k}}^{\bullet }\left[ \pi \right] }{\longrightarrow }{\cal P}_{{\bf 1}%
_{n-k}}^{\bullet }\left[ {\bf \Lambda }^{\left( \sigma _1,...,\sigma
_k+1\right) }\right] \stackrel{{\cal P}_{{\bf 1}_{n-k}}^{\bullet }\left[
j\right] }{\longrightarrow }{\cal P}_{{\bf 1}_{n-k}}^{\bullet }\left[
K_{\left( \sigma _1,...,\sigma _k+1\right) }^{\left( k\right) }\right] ; 
\]
but $\pi \circ j=0$ hence $im\left( \Phi \right) \subseteq \ker \left(
\Theta \right) $. We prove the reverse inclusion.

Let $P$ be an object in ${\frak D}$ and $h\in Hom_A^{\bullet }\left( {\bf %
\Lambda }^{\left( \sigma _1,...,\sigma _k+1,1\right) },{\rm D}_{{\bf 1}%
_{n-k-1}}\left( P\right) \subset {\rm Diff}_{{\bf 1}_{n-k-1}}^{+}\left(
P\right) \right) \simeq {\rm D}_{\left( \sigma _1,...,\sigma
_k+1,1,...,1\right) _n}\left( P\right) $ be such that 
\[
\Theta \left( h\right) =h\circ d_{\left( \sigma _1,...,\sigma _k+1,1\right)
}\circ j=0; 
\]
we claim that $h\in im\left( \Phi \right) $. Now 
\[
h\circ d_{\left( \sigma _1,...,\sigma _k+1,1\right) }\circ j=h\circ
j^{\prime }\circ d_{\left( \sigma _1,...,\sigma _k+1,1\right) }\mid
_{K_{\left( \sigma _1,...,\sigma _k+1\right) }^{\left( k\right) }}\text{ } 
\]
where 
\[
0\rightarrow K_{\left( \sigma _1,...,\sigma _k+1,1\right) }^{\left( k\right)
}\stackrel{j^{\prime }}{\longrightarrow }{\bf \Lambda }^{\left( \sigma
_1,...,\sigma _k+1,1\right) }\stackrel{\pi ^{\prime }}{\longrightarrow }{\bf %
\Lambda }^{\left( \sigma _1,...,\sigma _k,1\right) }\rightarrow 0. 
\]
But $h\in im\left( \Phi \right) $ iff $h\circ j^{\prime }=0$ so it is enough
to show that $im\left( d_{\left( \sigma _1,...,\sigma _k+1,1\right) }\mid
_{K_{\left( \sigma _1,...,\sigma _k+1\right) }^{\left( k\right) }}\right) $
generates $K_{\left( \sigma _1,...,\sigma _k+1,1\right) }^{\left( k\right) }$
over $A$ (since both $h$ and $j^{\prime }$ are $A$-homomorphisms). We know
that $im\left( j_1:Q\rightarrow {\bf J}^1\left( Q\right) \right) $ generates 
${\bf J}^1\left( Q\right) $ over $A$ for any object $Q$ in ${\frak D}$
(Section 2). Moreover, the 3$\times $3 lemma\footnote{%
Smoothness of ${\frak D}$ enters here.} gives us an exact commutative
diagram 
\[
\begin{tabular}{llllllllll}
&  & $0$ &  & $0$ &  & $0$ &  &  &  \\ 
&  & $\downarrow $ &  & $\downarrow $ &  & $\downarrow $ &  &  &  \\ 
$0$ & $\rightarrow $ & $K_{\left( \sigma _1,...,\sigma _k+2\right) }^{\left(
k\right) }$ & $\rightarrow $ & ${\bf J}^1\left( K_{\left( \sigma
_1,...,\sigma _k+1\right) }^{\left( k\right) }\right) $ & $\stackrel{t}{%
\rightarrow }$ & $K_{\left( \sigma _1,...,\sigma _k+1,1\right) }^{\left(
k\right) }$ & $\rightarrow $ & $0$ &  \\ 
&  & $\downarrow $ &  & $\downarrow $ &  & $\downarrow $ &  &  &  \\ 
$0$ & $\rightarrow $ & ${\bf \Lambda }^{\left( \sigma _1,...,\sigma
_k+2\right) }$ & $\rightarrow $ & ${\bf J}^1\left( {\bf \Lambda }^{\left(
\sigma _1,...,\sigma _k+1\right) }\right) $ & $\rightarrow $ & ${\bf \Lambda 
}^{\left( \sigma _1,...,\sigma _k+1,1\right) }$ & $\rightarrow $ & $0$ &  \\ 
&  & $\downarrow $ &  & $\downarrow $ &  & $\downarrow $ &  &  &  \\ 
$0$ & $\rightarrow $ & ${\bf \Lambda }^{\left( \sigma _1,...,\sigma
_k+1\right) }$ & $\rightarrow $ & ${\bf J}^1\left( {\bf \Lambda }^{\left(
\sigma _1,...,\sigma _k\right) }\right) $ & $\rightarrow $ & ${\bf \Lambda }%
^{\left( \sigma _1,...,\sigma _k,1\right) }$ & $\rightarrow $ & $0$ &  \\ 
&  & $\downarrow $ &  & $\downarrow $ &  & $\downarrow $ &  &  &  \\ 
&  & $0$ &  & $0$ &  & $0$ &  &  &  \\ 
&  &  &  &  &  &  &  &  & 
\end{tabular}
\]
(where we used the fact that the functor ${\bf J}^k\left( \cdot \right) $ is
exact if ${\frak D}$ is smooth: this follows from lemma \ref{lem 2.5} since $%
{\bf J}^k$ is projective); but $d_{\left( \sigma _1,...,\sigma _k+1,1\right)
}\mid _{K_{\left( \sigma _1,...,\sigma _k+1\right) }^{\left( k\right)
}}=t\circ j_1$ (by definition of $d$), hence $im\left( d_{\left( \sigma
_1,...,\sigma _k+1,1\right) }\mid _{K_{\left( \sigma _1,...,\sigma
_k+1\right) }^{\left( k\right) }}\right) $ generates $K_{\left( \sigma
_1,...,\sigma _k+1,1\right) }^{\left( k\right) }$ over $A$ and we have
finished.

\smallskip\ 

Now let's turn ourselves to the definition of 
\[
\psi _l\equiv \left( \psi _l^n:{\bf Hol}^{{\bf 1}_n}\left[ K_{\left( \sigma
_1,...,\sigma _k+l+1\right) }^{\left( k\right) }\right] \rightarrow {\bf Hol}%
^{{\bf 1}_{n+1}}\left[ K_{\left( \sigma _1,...,\sigma _k+l\right) }^{\left(
k\right) }\right] \right) _{n\geq 0}\text{ .} 
\]
First of all 
\[
\psi _l^0=0:\left( {\bf Hol}^{{\bf 1}}\left[ K_{\left( \sigma _1,...,\sigma
_k+l+1\right) }^{\left( k\right) }\right] \right) ^0=\left( 0\right)
\rightarrow \left( {\bf Hol}^{{\bf 1}}\left[ K_{\left( \sigma _1,...,\sigma
_k+l\right) }^{\left( k\right) }\right] \right) ^1={\bf J}^1\left( K_{\left(
\sigma _1,...,\sigma _k+l\right) }^{\left( k\right) }\right) . 
\]
For $n>0$%
\[
\psi _l^n:{\rm Hol}^{{\bf 1}_n}\left[ K_{\left( \sigma _1,...,\sigma
_k+l+1\right) }^{\left( k\right) }\right] \rightarrow {\rm Hol}^{{\bf 1}%
_{n+1}}\left[ K_{\left( \sigma _1,...,\sigma _k+l\right) }^{\left( k\right)
}\right] 
\]
will be defined as the dual representative of a functorial morphism 
\[
\psi _n^l:{\cal P}_{{\bf 1}_{n+1}}^{\bullet }\left[ K_{\left( \sigma
_1,...,\sigma _k+l\right) }^{\left( k\right) }\right] \rightarrow {\cal P}_{%
{\bf 1}_n}^{\bullet }\left[ K_{\left( \sigma _1,...,\sigma _k+l+1\right)
}^{\left( k\right) }\right] . 
\]
From the exact sequence 
\begin{eqnarray*}
0 &\rightarrow &K_{\left( \sigma _1,...,\sigma _k+l\right) }^{\left(
k\right) }\stackrel{i}{\longrightarrow }{\bf \Lambda }^{\left( \sigma
_1,...,\sigma _k+l\right) }\stackrel{p}{\longrightarrow }{\bf \Lambda }%
^{\left( \sigma _1,...,\sigma _k+l-1\right) }\rightarrow 0 \\
\text{(resp. }0 &\rightarrow &K_{\left( \sigma _1,...,\sigma _k+l+1\right)
}^{\left( k\right) }\stackrel{i^{\prime }}{\longrightarrow }{\bf \Lambda }%
^{\left( \sigma _1,...,\sigma _k+l+1\right) }\stackrel{p^{\prime }}{%
\longrightarrow }{\bf \Lambda }^{\left( \sigma _1,...,\sigma _k+l\right)
}\rightarrow 0\text{ )}
\end{eqnarray*}
and the fact that ${\bf \Lambda }^{\left( \sigma _1,...,\sigma _k+l-1\right)
}$ (resp. ${\bf \Lambda }^{\left( \sigma _1,...,\sigma _k+l\right) }$) is
projective we get (Prop. \ref{split}) an exact sequence of functors ${\frak D%
}\rightarrow {\frak D}$%
\begin{eqnarray*}
0 &\rightarrow &{\cal P}_{{\bf 1}_{n+1}}^{\bullet }\left[ {\bf \Lambda }%
^{\left( \sigma _1,...,\sigma _k+l-1\right) }\right] \stackrel{\epsilon }{%
\rightarrow }{\cal P}_{{\bf 1}_{n+1}}^{\bullet }\left[ {\bf \Lambda }%
^{\left( \sigma _1,...,\sigma _k+l\right) }\right] \stackrel{\eta }{%
\rightarrow }{\cal P}_{{\bf 1}_{n+1}}^{\bullet }\left[ K_{\left( \sigma
_1,...,\sigma _k+l\right) }^{\left( k\right) }\right] \rightarrow 0 \\
\text{(resp. }0 &\rightarrow &{\cal P}_{{\bf 1}_n}^{\bullet }\left[ {\bf %
\Lambda }^{\left( \sigma _1,...,\sigma _k+l\right) }\right] \stackrel{%
\epsilon ^{\prime }}{\rightarrow }{\cal P}_{{\bf 1}_n}^{\bullet }\left[ {\bf %
\Lambda }^{\left( \sigma _1,...,\sigma _k+l+1\right) }\right] \stackrel{\eta
^{\prime }}{\rightarrow }{\cal P}_{{\bf 1}_n}^{\bullet }\left[ K_{\left(
\sigma _1,...,\sigma _k+l+1\right) }^{\left( k\right) }\right] \rightarrow 0%
\text{ )}
\end{eqnarray*}
with $\epsilon ={\cal P}_{{\bf 1}_{n+1}}^{\bullet }\left[ p\right] $ and $%
\eta ={\cal P}_{{\bf 1}_{n+1}}^{\bullet }\left[ i\right] $ (resp. $\epsilon
^{\prime }={\cal P}_{{\bf 1}_{n+1}}^{\bullet }\left[ p^{\prime }\right] $
and $\eta ^{\prime }={\cal P}_{{\bf 1}_{n+1}}^{\bullet }\left[ i^{\prime
}\right] $). To define $\psi _n^l$ it will be then enough to define 
\[
\overline{\psi }_n^l:{\cal P}_{{\bf 1}_{n+1}}^{\bullet }\left[ {\bf \Lambda }%
^{\left( \sigma _1,...,\sigma _k+l\right) }\right] \rightarrow {\cal P}_{%
{\bf 1}_n}^{\bullet }\left[ K_{\left( \sigma _1,...,\sigma _k+l+1\right)
}^{\left( k\right) }\right] 
\]
and show that $\overline{\psi }_n^l\circ \epsilon =0$. We know (Section 2 )
that there is an exact sequence 
\[
0\rightarrow {\bf \Lambda }^{\left( \sigma _1,...,\sigma _k+l+1\right) }%
\stackrel{s}{\longrightarrow }{\bf J}^1\left( {\bf \Lambda }^{\left( \sigma
_1,...,\sigma _k+l\right) }\right) \stackrel{q}{\longrightarrow }{\bf %
\Lambda }^{\left( \sigma _1,...,\sigma _k+l,1\right) }\rightarrow 0 
\]
and this (Prop. \ref{split}) gives us the exact sequence 
\[
0\rightarrow {\cal P}_{{\bf 1}_n}^{\bullet }\left[ {\bf \Lambda }^{\left(
\sigma _1,...,\sigma _k+l,1\right) }\right] \stackrel{\alpha }{\rightarrow }%
{\cal P}_{{\bf 1}_n}^{\bullet }\left[ {\bf J}^1\left( {\bf \Lambda }^{\left(
\sigma _1,...,\sigma _k+l\right) }\right) \right] \stackrel{\beta }{%
\rightarrow }{\cal P}_{{\bf 1}_n}^{\bullet }\left[ {\bf \Lambda }^{\left(
\sigma _1,...,\sigma _k+l+1\right) }\right] \rightarrow 0 
\]
(with $\alpha ={\cal P}_{{\bf 1}_n}^{\bullet }\left[ q\right] $ and $\beta =%
{\cal P}_{{\bf 1}_n}^{\bullet }\left[ s\right] $). Then we take $\overline{%
\psi }_n^l$ to be the composition 
\begin{eqnarray*}
{\cal P}_{{\bf 1}_{n+1}}^{\bullet }\left[ {\bf \Lambda }^{\left( \sigma
_1,...,\sigma _k+l\right) }\right] \simeq {\cal P}_{\left( 1,1\right)
}^{\bullet }\left[ {\bf \Lambda }^{\left( \sigma _1,...,\sigma _k+l\right)
}\right] \left( {\rm D}_{{\bf 1}_{n-1}}\subset {\rm Diff}_{{\bf 1}%
_{n-1}}^{+}\right) \simeq
\end{eqnarray*}
\[
\simeq Hom_A^{\bullet }\left( {\rm Hol}^{(1,1)}\left[ {\bf \Lambda }^{\left(
\sigma _1,...,\sigma _k+l\right) }\right] ,{\rm D}_{{\bf 1}_{n-1}}\subset 
{\rm Diff}_{{\bf 1}_{n-1}}^{+}\right) \stackrel{\circ \text{ }\delta
_{\left( 1,1\right) }\left( {\bf \Lambda }^{\left( \sigma _1,...,\sigma
_k+l\right) }\right) }{\longrightarrow } 
\]
\begin{eqnarray*}
\ &\rightarrow &{\rm Diff}_1^{\bullet }\left( {\rm Hol}^{(1)}\left[ {\bf %
\Lambda }^{\left( \sigma _1,...,\sigma _k+l\right) }\right] ,{\rm D}_{{\bf 1}%
_{n-1}}\subset {\rm Diff}_{{\bf 1}_{n-1}}^{+}\right) \simeq \\
\ &\simeq &{\cal P}_{\left( 1\right) }^{\bullet }\left[ {\bf J}^1\left( {\bf %
\Lambda }^{\left( \sigma _1,...,\sigma _k+l\right) }\right) \right] \left( 
{\rm D}_{{\bf 1}_{n-1}}\subset {\rm Diff}_{{\bf 1}_{n-1}}^{+}\right) \simeq
\end{eqnarray*}
\[
\simeq {\cal P}_{{\bf 1}_n}^{\bullet }\left[ {\bf J}^1\left( {\bf \Lambda }%
^{\left( \sigma _1,...,\sigma _k+l\right) }\right) \right] \stackrel{\beta }{%
\longrightarrow }{\cal P}_{{\bf 1}_n}^{\bullet }\left[ {\bf \Lambda }%
^{\left( \sigma _1,...,\sigma _k+l+1\right) }\right] \stackrel{\eta ^{\prime
}}{\longrightarrow }{\cal P}_{{\bf 1}_n}^{\bullet }\left[ K_{\left( \sigma
_1,...,\sigma _k+l+1\right) }^{\left( k\right) }\right] . 
\]
Now we show that $\overline{\psi }_n^l\circ \epsilon =0$.

Note that using the identifications 
\begin{eqnarray*}
{\cal P}_{{\bf 1}_n}^{\bullet }\left[ {\bf J}^1\left( {\bf \Lambda }^{\left(
\sigma _1,...,\sigma _k+l\right) }\right) \right] &\simeq &{\rm Diff}%
_1^{\bullet }\left( {\bf J}^1\left( {\bf \Lambda }^{\left( \sigma
_1,...,\sigma _k+l\right) }\right) ,{\rm D}_{{\bf 1}_{n-1}}\subset {\rm Diff}%
_{{\bf 1}_{n-1}}^{+}\right) \simeq \\
\ &\simeq &Hom_A^{\bullet }\left( {\bf J}^1\left( {\bf J}^1\left( {\bf %
\Lambda }^{\left( \sigma _1,...,\sigma _k+l\right) }\right) \right) ,{\rm D}%
_{{\bf 1}_{n-1}}\subset {\rm Diff}_{{\bf 1}_{n-1}}^{+}\right)
\end{eqnarray*}
and 
\[
Hom_A^{\bullet }\left( {\bf J}^1\left( {\bf \Lambda }^{\left( \sigma
_1,...,\sigma _k+l+1\right) }\right) ,{\rm D}_{{\bf 1}_{n-1}}\subset {\rm %
Diff}_{{\bf 1}_{n-1}}^{+}\right) \simeq {\cal P}_{{\bf 1}_n}^{\bullet
}\left[ {\bf J}^1\left( {\bf \Lambda }^{\left( \sigma _1,...,\sigma
_k+l+1\right) }\right) \right] 
\]
(resp. the identifications 
\[
{\cal P}_{{\bf 1}_n}^{\bullet }\left[ {\bf \Lambda }^{\left( \sigma
_1,...,\sigma _k+l+1\right) }\right] \simeq Hom_A^{\bullet }\left( {\bf J}%
^1\left( {\bf \Lambda }^{\left( \sigma _1,...,\sigma _k+l+1\right) }\right) ,%
{\rm D}_{{\bf 1}_{n-1}}\subset {\rm Diff}_{{\bf 1}_{n-1}}^{+}\right) 
\]
and 
\[
Hom_A^{\bullet }\left( {\bf J}^1\left( K_{\left( \sigma _1,...,\sigma
_k+l+1\right) }^{\left( k\right) }\right) ,{\rm D}_{{\bf 1}_{n-1}}\subset 
{\rm Diff}_{{\bf 1}_{n-1}}^{+}\right) \simeq {\cal P}_{{\bf 1}_n}^{\bullet
}\left[ K_{\left( \sigma _1,...,\sigma _k+l+1\right) }^{\left( k\right)
}\right] \text{ )} 
\]
$\beta $ (resp. $\eta ^{\prime }$ ) is given by 
\begin{eqnarray*}
&&\ \ \ \ Hom_A^{\bullet }\left( {\bf J}^1\left( {\bf J}^1\left( {\bf %
\Lambda }^{\left( \sigma _1,...,\sigma _k+l\right) }\right) \right) ,{\rm D}%
_{{\bf 1}_{n-1}}\subset {\rm Diff}_{{\bf 1}_{n-1}}^{+}\right) \stackrel{%
\circ \text{ }{\bf J}^1\left( s\right) }{\longrightarrow } \\
\ &\longrightarrow &Hom_A^{\bullet }\left( {\bf J}^1\left( {\bf \Lambda }%
^{\left( \sigma _1,...,\sigma _k+l+1\right) }\right) ,{\rm D}_{{\bf 1}%
_{n-1}}\subset {\rm Diff}_{{\bf 1}_{n-1}}^{+}\right)
\end{eqnarray*}
(resp. by 
\begin{eqnarray*}
&&\ \ \ \ \ \ Hom_A^{\bullet }\left( {\bf J}^1\left( {\bf \Lambda }^{\left(
\sigma _1,...,\sigma _k+l+1\right) }\right) ,{\rm D}_{{\bf 1}_{n-1}}\subset 
{\rm Diff}_{{\bf 1}_{n-1}}^{+}\right) \stackrel{\circ \text{ }{\bf J}%
^1\left( i^{\prime }\right) }{\longrightarrow } \\
\ &\rightarrow &Hom_A^{\bullet }\left( {\bf J}^1\left( K_{\left( \sigma
_1,...,\sigma _k+l+1\right) }^{\left( k\right) }\right) ,{\rm D}_{{\bf 1}%
_{n-1}}\subset {\rm Diff}_{{\bf 1}_{n-1}}^{+}\right) \text{ ).}
\end{eqnarray*}
With a similar analysis we see that $\epsilon $, viewed as a morphism 
\begin{eqnarray*}
&&\ \ \ \ \ \ Hom_A^{\bullet }\left( {\rm Hol}^{\left( 1,1\right) }\left( 
{\bf \Lambda }^{\left( \sigma _1,...,\sigma _k+l-1\right) }\right) ,{\rm D}_{%
{\bf 1}_{n-1}}\subset {\rm Diff}_{{\bf 1}_{n-1}}^{+}\right) \stackrel{%
\epsilon }{\rightarrow } \\
\ &\rightarrow &Hom_A^{\bullet }\left( {\rm Hol}^{\left( 1,1\right) }\left( 
{\bf \Lambda }^{\left( \sigma _1,...,\sigma _k+l\right) }\right) ,{\rm D}_{%
{\bf 1}_n}\subset {\rm Diff}_{{\bf 1}_n}^{+}\right)
\end{eqnarray*}
is given by taking the composition with $\left[ {\bf J}^1\left( {\bf J}%
^1\left( p\right) \right) \right] $ where 
\begin{eqnarray*}
\left[ {\bf J}^1\left( {\bf J}^1\left( p\right) \right) \right] &:&\frac{%
{\bf J}^1\left( {\bf J}^1\left( {\bf \Lambda }^{\left( \sigma _1,...,\sigma
_k+l-1\right) }\right) \right) }{{\bf J}^2\left( {\bf \Lambda }^{\left(
\sigma _1,...,\sigma _k+l-1\right) }\right) }\simeq {\rm Hol}^{\left(
1,1\right) }\left( {\bf \Lambda }^{\left( \sigma _1,...,\sigma _k+l-1\right)
}\right) \rightarrow \\
\ &\rightarrow &{\rm Hol}^{\left( 1,1\right) }\left( {\bf \Lambda }^{\left(
\sigma _1,...,\sigma _k+l\right) }\right) \simeq \frac{{\bf J}^1\left( {\bf J%
}^1\left( {\bf \Lambda }^{\left( \sigma _1,...,\sigma _k+l\right) }\right)
\right) }{{\bf J}^2\left( {\bf \Lambda }^{\left( \sigma _1,...,\sigma
_k+l\right) }\right) }
\end{eqnarray*}
is the quotient map of ${\bf J}^1\left( {\bf J}^1\left( p\right) \right) :%
{\bf J}^1\left( {\bf J}^1\left( {\bf \Lambda }^{\left( \sigma _1,...,\sigma
_k+l-1\right) }\right) \right) \rightarrow {\bf J}^1\left( {\bf J}^1\left( 
{\bf \Lambda }^{\left( \sigma _1,...,\sigma _k+l\right) }\right) \right) $.

Therefore $\overline{\psi }_n^l\circ \epsilon $, viewed as a morphism 
\[
Hom_A^{\bullet }\left( {\rm Hol}^{\left( 1,1\right) }\left( {\bf \Lambda }%
^{\left( \sigma _1,...,\sigma _k+l-1\right) }\right) ,{\rm D}_{{\bf 1}%
_{n-1}}\subset {\rm Diff}_{{\bf 1}_{n-1}}^{+}\right) \longrightarrow 
\]
\[
\longrightarrow \ \ \ Hom_A^{\bullet }\left( {\bf J}^1\left( K_{\left(
\sigma _1,...,\sigma _k+l+1\right) }^{\left( k\right) }\right) ,{\rm D}_{%
{\bf 1}_{n-1}}\subset {\rm Diff}_{{\bf 1}_{n-1}}^{+}\right) 
\]
is given by 
\[
f\longmapsto f\circ \left[ {\bf J}^1\left( {\bf J}^1\left( p\right) \right)
\right] \circ \xi _{\left( \sigma _1,...,\sigma _k+l\right) }\circ {\bf J}%
^1\left( s\right) \circ {\bf J}^1\left( i^{\prime }\right) 
\]
where 
\[
\xi _{\left( \sigma _1,...,\sigma _k+l\right) }:{\bf J}^1\left( {\bf J}%
^1\left( {\bf \Lambda }^{\left( \sigma _1,...,\sigma _k+l\right) }\right)
\right) \rightarrow \frac{{\bf J}^1\left( {\bf J}^1\left( {\bf \Lambda }%
^{\left( \sigma _1,...,\sigma _k+l\right) }\right) \right) }{{\bf J}^2\left( 
{\bf \Lambda }^{\left( \sigma _1,...,\sigma _k+l\right) }\right) } 
\]
is the natural projection. (Recall from Section 4 that $\delta _{\left(
1,1\right) }\left( {\bf \Lambda }^{\left( \sigma _1,...,\sigma _k+l\right)
}\right) $ is given by the composition 
\[
{\bf J}^1\left( {\bf \Lambda }^{\left( \sigma _1,...,\sigma _k+l\right)
}\right) \stackrel{j_1\left( {\bf \Lambda }^{\left( \sigma _1,...,\sigma
_k+l\right) }\right) }{\longrightarrow }{\bf J}^1\left( {\bf J}^1\left( {\bf %
\Lambda }^{\left( \sigma _1,...,\sigma _k+l\right) }\right) \right) 
\stackrel{\xi _{\left( \sigma _1,...,\sigma _k+l\right) }}{\longrightarrow }%
\frac{{\bf J}^1\left( {\bf J}^1\left( {\bf \Lambda }^{\left( \sigma
_1,...,\sigma _k+l\right) }\right) \right) }{{\bf J}^2\left( {\bf \Lambda }%
^{\left( \sigma _1,...,\sigma _k+l\right) }\right) }\text{ ).} 
\]
But $\left[ {\bf J}^1\left( {\bf J}^1\left( p\right) \right) \right] \circ
\xi $ coincides with 
\[
{\bf J}^1\left( {\bf J}^1\left( {\bf \Lambda }^{\left( \sigma _1,...,\sigma
_k+l\right) }\right) \right) \stackrel{{\bf J}^1\left( {\bf J}^1\left(
p\right) \right) }{\longrightarrow }{\bf J}^1\left( {\bf J}^1\left( {\bf %
\Lambda }^{\left( \sigma _1,...,\sigma _k+l-1\right) }\right) \right) 
\stackrel{\xi _{\left( \sigma _1,...,\sigma _k+l-1\right) }}{\longrightarrow 
}\frac{{\bf J}^1\left( {\bf J}^1\left( {\bf \Lambda }^{\left( \sigma
_1,...,\sigma _k+l\right) }\right) \right) }{{\bf J}^2\left( {\bf \Lambda }%
^{\left( \sigma _1,...,\sigma _k+l\right) }\right) } 
\]
($\xi _{\left( \sigma _1,...,\sigma _k+l-1\right) }$ being again the natural
projection), so that 
\[
\begin{tabular}{lll}
$\overline{\psi }_n^l\circ \epsilon :f$ & $\longmapsto $ & $f\circ \xi
_{\left( \sigma _1,...,\sigma _k+l\right) }\circ {\bf J}^1\left( {\bf J}%
^1\left( p\right) \right) \circ {\bf J}^1\left( s\right) \circ {\bf J}%
^1\left( i^{\prime }\right) =$ \\ 
&  & $=f\circ \xi _{\left( \sigma _1,...,\sigma _k+l\right) }\circ {\bf J}%
^1\left( {\bf J}^1\left( p\right) \circ s\circ i^{\prime }\right) .$ \\ 
&  & 
\end{tabular}
\]
Again as above, the 3$\times $3 lemma gives us an exact commutative diagram: 
\[
\begin{tabular}{ccccccccc}
&  & $0$ &  & $0$ &  & $0$ &  &  \\ 
&  & $\downarrow $ &  & $\downarrow $ &  & $\downarrow $ &  &  \\ 
$0$ & $\rightarrow $ & $K_{\left( \sigma _1,...,\sigma _k+l+1\right)
}^{\left( k\right) }$ & $\rightarrow $ & ${\bf J}^1\left( K_{\left( \sigma
_1,...,\sigma _k+l\right) }^{\left( k\right) }\right) $ & $\stackrel{}{%
\rightarrow }$ & $K_{\left( \sigma _1,...,\sigma _k+l,1\right) }^{\left(
k\right) }$ & $\rightarrow $ & $0$ \\ 
&  & $i^{\prime }\downarrow \quad $ &  & $\downarrow $ &  & $\downarrow $ & 
&  \\ 
$0$ & $\rightarrow $ & ${\bf \Lambda }^{\left( \sigma _1,...,\sigma
_k+l+1\right) }$ & $\stackrel{s}{\rightarrow }$ & ${\bf J}^1\left( {\bf %
\Lambda }^{\left( \sigma _1,...,\sigma _k+l\right) }\right) $ & $\rightarrow 
$ & ${\bf \Lambda }^{\left( \sigma _1,...,\sigma _k+l,1\right) }$ & $%
\rightarrow $ & $0$ \\ 
&  & $\downarrow $ &  & $_{{\bf J}^1\left( p\right) }\downarrow $ $\qquad $
&  & $\downarrow $ &  &  \\ 
$0$ & $\rightarrow $ & ${\bf \Lambda }^{\left( \sigma _1,...,\sigma
_k+l\right) }$ & $\rightarrow $ & ${\bf J}^1\left( {\bf \Lambda }^{\left(
\sigma _1,...,\sigma _k+l-1\right) }\right) $ & $\rightarrow $ & ${\bf %
\Lambda }^{\left( \sigma _1,...,\sigma _k+l-1,1\right) }$ & $\rightarrow $ & 
$0$ \\ 
&  & $\downarrow $ &  & $\downarrow $ &  & $\downarrow $ &  &  \\ 
&  & $0$ &  & $0$ &  & $0$ &  & 
\end{tabular}
\]
which finally shows that ${\bf J}^1\left( p\right) \circ s\circ i^{\prime }$ 
$=0$ and hence $\overline{\psi }_n^l\circ \epsilon =0$.

Therefore 
\[
\psi _n^l:{\cal P}_{{\bf 1}_{n+1}}^{\bullet }\left[ K_{\left( \sigma
_1,...,\sigma _k+l\right) }^{\left( k\right) }\right] \rightarrow {\cal P}_{%
{\bf 1}_n}^{\bullet }\left[ K_{\left( \sigma _1,...,\sigma _k+l+1\right)
}^{\left( k\right) }\right] 
\]
is well defined as well as its dual representative 
\[
\psi _l^n:{\rm Hol}^{{\bf 1}_n}\left[ K_{\left( \sigma _1,...,\sigma
_k+l+1\right) }^{\left( k\right) }\right] \rightarrow {\rm Hol}^{{\bf 1}%
_{n+1}}\left[ K_{\left( \sigma _1,...,\sigma _k+l\right) }^{\left( k\right)
}\right] 
\]
as we wanted.

Just as in the case of $\rho $, an easy application of the 3$\times $3 lemma
proves that $im\left( \psi _{n+1}^{l-1}\right) =\ker \left( \psi _n^l\right) 
$.

It is easy to verify that $\rho $ and $\psi _l$ so defined are maps of
complexes; therefore (\ref{res}) is a resolution of ${\bf K}_{\left( \sigma
_1,...,\sigma _k+1,{\bf 1}\right) }^{(k)}$ as desired.

\begin{center}
\smallskip 
\ 
\end{center}

\end{document}

%% file: diagram.tex
\makeatletter

\def\diagram{\leftwidth=\z@ \rightwidth=\z@ \topheight=\z@
\botheight=\z@ \setbox\@picbox\hbox\bgroup}

\def\enddiagram{\egroup\wd\@picbox\rightwidth\unitlength
\ht\@picbox\topheight\unitlength \dp\@picbox\botheight\unitlength
\hskip\leftwidth\unitlength\box\@picbox}

\def\bfig{\begin{diagram}}
\def\efig{\end{diagram}}
\newcount\wideness \newcount\leftwidth \newcount\rightwidth
\newcount\highness \newcount\topheight \newcount\botheight

\def\ratchet#1#2{\ifnum#1<#2 \global #1=#2 \fi}

\def\putbox(#1,#2)#3{%
\horsize{\wideness}{#3} \divide\wideness by 2
{\advance\wideness by #1 \ratchet{\rightwidth}{\wideness}}
{\advance\wideness by -#1 \ratchet{\leftwidth}{\wideness}}
\vertsize{\highness}{#3} \divide\highness by 2
{\advance\highness by #2 \ratchet{\topheight}{\highness}}
{\advance\highness by -#2 \ratchet{\botheight}{\highness}}
\put(#1,#2){\makebox(0,0){$#3$}}}

\def\putlbox(#1,#2)#3{%
\horsize{\wideness}{#3}
{\advance\wideness by #1 \ratchet{\rightwidth}{\wideness}}
{\ratchet{\leftwidth}{-#1}}
\vertsize{\highness}{#3} \divide\highness by 2
{\advance\highness by #2 \ratchet{\topheight}{\highness}}
{\advance\highness by -#2 \ratchet{\botheight}{\highness}}
\put(#1,#2){\makebox(0,0)[l]{$#3$}}}

\def\putrbox(#1,#2)#3{%
\horsize{\wideness}{#3}
{\ratchet{\rightwidth}{#1}}
{\advance\wideness by -#1 \ratchet{\leftwidth}{\wideness}}
\vertsize{\highness}{#3} \divide\highness by 2
{\advance\highness by #2 \ratchet{\topheight}{\highness}}
{\advance\highness by -#2 \ratchet{\botheight}{\highness}}
\put(#1,#2){\makebox(0,0)[r]{$#3$}}}

\def\adjust[#1]{} % For compatibility

\newcount \coefa
\newcount \coefb
\newcount \coefc
\newcount\tempcounta
\newcount\tempcountb
\newcount\tempcountc
\newcount\tempcountd
\newcount\xext
\newcount\yext
\newcount\xoff
\newcount\yoff
\newcount\gap%
\newcount\arrowtypea
\newcount\arrowtypeb
\newcount\arrowtypec
\newcount\arrowtyped
\newcount\arrowtypee
\newcount\height
\newcount\width
\newcount\xpos
\newcount\ypos
\newcount\run
\newcount\rise
\newcount\arrowlength
\newcount\halflength
\newcount\arrowtype
\newdimen\tempdimen
\newdimen\xlen
\newdimen\ylen
\newsavebox{\tempboxa}%
\newsavebox{\tempboxb}%
\newsavebox{\tempboxc}%

\newdimen\w@dth

\def\setw@dth#1#2{\setbox\z@\hbox{$#1$}\w@dth=\wd\z@
\setbox\@ne\hbox{$#2$}\ifnum\w@dth<\wd\@ne \w@dth=\wd\@ne \fi
\advance\w@dth by 1.2em}

\def\t@^#1_#2{\def\n@one{#1}\def\n@two{#2}\mathrel{\setw@dth{#1}{#2}
\mathop{\hbox to \w@dth{\rightarrowfill}}\limits
\ifx\n@one\empty\else ^{\box\z@}\fi
\ifx\n@two\empty\else _{\box\@ne}\fi}}
\def\t@@^#1{\@ifnextchar_ {\t@^{#1}}{\t@^{#1}_{}}}
\def\to{\@ifnextchar^ {\t@@}{\t@@^{}}}

\def\t@left^#1_#2{\def\n@one{#1}\def\n@two{#2}\mathrel{\setw@dth{#1}{#2}
\mathop{\hbox to \w@dth{\leftarrowfill}}\limits
\ifx\n@one\empty\else ^{\box\z@}\fi
\ifx\n@two\empty\else _{\box\@ne}\fi}}
\def\t@@left^#1{\@ifnextchar_ {\t@left^{#1}}{\t@left^{#1}_{}}}
\def\toleft{\@ifnextchar^ {\t@@left}{\t@@left^{}}}

\def\two@^#1_#2{\def\n@one{#1}\def\n@two{#2}\mathrel{\setw@dth{#1}{#2}
\mathop{\vcenter{\hbox to \w@dth{\rightarrowfill}\kern-1.7ex
                 \hbox to \w@dth{\rightarrowfill}}%
       }\limits
\ifx\n@one\empty\else ^{\box\z@}\fi
\ifx\n@two\empty\else _{\box\@ne}\fi}}
\def\tw@@^#1{\@ifnextchar_ {\two@^{#1}}{\two@^{#1}_{}}}
\def\two{\@ifnextchar^ {\tw@@}{\tw@@^{}}}

\def\tofr@^#1_#2{\def\n@one{#1}\def\n@two{#2}\mathrel{\setw@dth{#1}{#2}
\mathop{\vcenter{\hbox to \w@dth{\rightarrowfill}\kern-1.7ex
                 \hbox to \w@dth{\leftarrowfill}}%
       }\limits
\ifx\n@one\empty\else ^{\box\z@}\fi
\ifx\n@two\empty\else _{\box\@ne}\fi}}
\def\t@fr@^#1{\@ifnextchar_ {\tofr@^{#1}}{\tofr@^{#1}_{}}}
\def\tofro{\@ifnextchar^ {\t@fr@}{\t@fr@^{}}}

\def\mon{\mathop{\m@th\hbox to
      14.6\P@{\lasyb\char'51\hskip-2.1\P@$\arrext$\hss
$\mathord\rightarrow$}}\limits} % width of \epi
\def\leftmono{\mathrel{\m@th\hbox to
14.6\P@{$\mathord\leftarrow$\hss$\arrext$\hskip-2.1\P@\lasyb\char'50%
}}\limits} % width of \epi
\mathchardef\arrext="0200       % amr minus for arrow extension (see \into)

\def\settypes(#1,#2,#3){\arrowtypea#1 \arrowtypeb#2 \arrowtypec#3}
\def\settoheight#1#2{\setbox\@tempboxa\hbox{#2}#1\ht\@tempboxa\relax}%
\def\settodepth#1#2{\setbox\@tempboxa\hbox{#2}#1\dp\@tempboxa\relax}%
\def\settokens[#1`#2`#3`#4]{%
     \def\tokena{#1}\def\tokenb{#2}\def\tokenc{#3}\def\tokend{#4}}
\def\setsqparms[#1`#2`#3`#4;#5`#6]{%
\arrowtypea #1
\arrowtypeb #2
\arrowtypec #3
\arrowtyped #4
\width #5
\height #6
}
\def\setpos(#1,#2){\xpos=#1 \ypos#2}

\def\settriparms[#1`#2`#3;#4]{\settripairparms[#1`#2`#3`1`1;#4]}%

\def\settripairparms[#1`#2`#3`#4`#5;#6]{%
\arrowtypea #1
\arrowtypeb #2
\arrowtypec #3
\arrowtyped #4
\arrowtypee #5
\width #6
\height #6
}

\def\resetparms{\settripairparms[1`1`1`1`1;500]\width 500}%default values%

\resetparms

\def\mvector(#1,#2)#3{%%
\put(0,0){\vector(#1,#2){#3}}%
\put(0,0){\vector(#1,#2){26}}%
}
\def\evector(#1,#2)#3{{%%
\arrowlength #3
\put(0,0){\vector(#1,#2){\arrowlength}}%
\advance \arrowlength by-30
\put(0,0){\vector(#1,#2){\arrowlength}}%
}}

\def\horsize#1#2{%
\settowidth{\tempdimen}{$#2$}%
#1=\tempdimen
\divide #1 by\unitlength
}

\def\vertsize#1#2{%
\settoheight{\tempdimen}{$#2$}%
#1=\tempdimen
\settodepth{\tempdimen}{$#2$}%
\advance #1 by\tempdimen
\divide #1 by\unitlength
}

\def\putvector(#1,#2)(#3,#4)#5#6{{%
\ifnum3<\arrowtype
\putdashvector(#1,#2)(#3,#4)#5\arrowtype
\else
\ifnum\arrowtype<-3
\putdashvector(#1,#2)(#3,#4)#5\arrowtype
\else
\xpos=#1
\ypos=#2
\run=#3
\rise=#4
\arrowlength=#5
\ifnum \arrowtype<0
    \ifnum \run=0
        \advance \ypos by-\arrowlength
    \else
        \tempcounta \arrowlength
        \multiply \tempcounta by\rise
        \divide \tempcounta by\run
        \ifnum\run>0
            \advance \xpos by\arrowlength
            \advance \ypos by\tempcounta
        \else
            \advance \xpos by-\arrowlength
            \advance \ypos by-\tempcounta
        \fi
    \fi
    \multiply \arrowtype by-1
    \multiply \rise by-1
    \multiply \run by-1
\fi
\ifcase \arrowtype
\or \put(\xpos,\ypos){\vector(\run,\rise){\arrowlength}}%
\or \put(\xpos,\ypos){\mvector(\run,\rise)\arrowlength}%
\or \put(\xpos,\ypos){\evector(\run,\rise){\arrowlength}}%
\fi\fi\fi
}}

\def\putsplitvector(#1,#2)#3#4{%%
\xpos #1
\ypos #2
\arrowtype #4
\halflength #3
\arrowlength #3
\gap 140
\advance \halflength by-\gap
\divide \halflength by2
\ifnum\arrowtype>0
   \ifcase \arrowtype
   \or \put(\xpos,\ypos){\line(0,-1){\halflength}}%
       \advance\ypos by-\halflength
       \advance\ypos by-\gap
       \put(\xpos,\ypos){\vector(0,-1){\halflength}}%
   \or \put(\xpos,\ypos){\line(0,-1)\halflength}%
       \put(\xpos,\ypos){\vector(0,-1)3}%
       \advance\ypos by-\halflength
       \advance\ypos by-\gap
       \put(\xpos,\ypos){\vector(0,-1){\halflength}}%
   \or \put(\xpos,\ypos){\line(0,-1)\halflength}%
       \advance\ypos by-\halflength
       \advance\ypos by-\gap
       \put(\xpos,\ypos){\evector(0,-1){\halflength}}%
   \fi
\else \arrowtype=-\arrowtype
   \ifcase\arrowtype
   \or \advance \ypos by-\arrowlength
       \put(\xpos,\ypos){\line(0,1){\halflength}}%
       \advance\ypos by\halflength
       \advance\ypos by\gap
       \put(\xpos,\ypos){\vector(0,1){\halflength}}%
   \or \advance \ypos by-\arrowlength
       \put(\xpos,\ypos){\line(0,1)\halflength}%
       \put(\xpos,\ypos){\vector(0,1)3}%
       \advance\ypos by\halflength
       \advance\ypos by\gap
       \put(\xpos,\ypos){\vector(0,1){\halflength}}%
   \or \advance \ypos by-\arrowlength
       \put(\xpos,\ypos){\line(0,1)\halflength}%
       \advance\ypos by\halflength
       \advance\ypos by\gap
       \put(\xpos,\ypos){\evector(0,1){\halflength}}%
   \fi
\fi
}

\def\putmorphism(#1)(#2,#3)[#4`#5`#6]#7#8#9{{%
\run #2
\rise #3
\ifnum\rise=0
  \puthmorphism(#1)[#4`#5`#6]{#7}{#8}#9%
\else\ifnum\run=0
  \putvmorphism(#1)[#4`#5`#6]{#7}{#8}#9%
\else
\setpos(#1)%
\arrowlength #7
\arrowtype #8
\ifnum\run=0
\else\ifnum\rise=0
\else
\ifnum\run>0
    \coefa=1
\else
   \coefa=-1
\fi
\ifnum\arrowtype>0
   \coefb=0
   \coefc=-1
\else
   \coefb=\coefa
   \coefc=1
   \arrowtype=-\arrowtype
\fi
\width=2
\multiply \width by\run
\divide \width by\rise
\ifnum \width<0  \width=-\width\fi
\advance\width by60
\if l#9 \width=-\width\fi
\putbox(\xpos,\ypos){#4}%            %node 1
{\multiply \coefa by\arrowlength%      %node 2
\advance\xpos by\coefa
\multiply \coefa by\rise
\divide \coefa by\run
\advance \ypos by\coefa
\putbox(\xpos,\ypos){#5} }%
{\multiply \coefa by\arrowlength%      %label
\divide \coefa by2
\advance \xpos by\coefa
\advance \xpos by\width
\multiply \coefa by\rise
\divide \coefa by\run
\advance \ypos by\coefa
\if l#9%
   \putrbox(\xpos,\ypos){#6}%
\else\if r#9%
   \putlbox(\xpos,\ypos){#6}%
\fi\fi }%
{\multiply \rise by-\coefc%             %arrow
\multiply \run by-\coefc
\multiply \coefb by\arrowlength
\advance \xpos by\coefb
\multiply \coefb by\rise
\divide \coefb by\run
\advance \ypos by\coefb
\multiply \coefc by70
\advance \ypos by\coefc
\multiply \coefc by\run
\divide \coefc by\rise
\advance \xpos by\coefc
\multiply \coefa by140
\multiply \coefa by\run
\divide \coefa by\rise
\advance \arrowlength by\coefa
\ifcase\arrowtype
\or \put(\xpos,\ypos){\vector(\run,\rise){\arrowlength}}%
\or \put(\xpos,\ypos){\mvector(\run,\rise){\arrowlength}}%
\or \put(\xpos,\ypos){\evector(\run,\rise){\arrowlength}}%
\fi}\fi\fi\fi\fi}}

\newcount\numbdashes \newcount\lengthdash \newcount\increment

\def\howmanydashes{% Actually returns both number and length
\numbdashes=\arrowlength \lengthdash=40
\divide\numbdashes by \lengthdash
\lengthdash=\arrowlength
\divide\lengthdash by \numbdashes
\increment=\lengthdash
\multiply\lengthdash by 3
\divide\lengthdash by 5
}

\def\putdashvector(#1)(#2,#3)#4#5{%
\ifnum#3=0 \putdashhvector(#1){#4}#5
\else
\ifnum#2=0
\putdashvvector(#1){#4}#5\fi\fi}

\def\putdashhvector(#1,#2)#3#4{{%
\arrowlength=#3 \howmanydashes
\multiput(#1,#2)(\increment,0){\numbdashes}%
{\vrule height .4pt width \lengthdash\unitlength}
\arrowtype=#4 \xpos=#1
\ifnum\arrowtype<0 \advance\arrowtype by 7 \fi
\ifcase\arrowtype
\or \advance\xpos by 10
    \put(\xpos,#2){\vector(-1,0){\lengthdash}}
    \advance\xpos by 40
    \put(\xpos,#2){\vector(-1,0){\lengthdash}}
\or \advance \xpos by 10
    \put(\xpos,#2){\vector(-1,0){\lengthdash}}
    \advance\xpos by  \arrowlength
    \advance\xpos by  -50
    \put(\xpos,#2){\vector(-1,0){\lengthdash}}
\or \advance\xpos by 10
    \put(\xpos,#2){\vector(-1,0){\lengthdash}}
\or \advance\xpos by \arrowlength
    \advance\xpos by -\lengthdash
    \put(\xpos,#2){\vector(1,0){\lengthdash}}
\or {\advance\xpos by 10
    \put(\xpos,#2){\vector(1,0){\lengthdash}}}
    \advance\xpos by \arrowlength
    \advance\xpos by -\lengthdash
    \put(\xpos,#2){\vector(1,0){\lengthdash}}
\or \advance\xpos by \arrowlength
    \advance\xpos by -\lengthdash
    \put(\xpos,#2){\vector(1,0){\lengthdash}}
    \advance\xpos by -40
    \put(\xpos,#2){\vector(1,0){\lengthdash}}
   \fi
}}

\def\putdashvvector(#1,#2)#3#4{{%
\arrowlength=#3 \howmanydashes
\ypos=#2 \advance\ypos by -\arrowlength
\multiput(#1,#2)(0,\increment){\numbdashes}%
    {\vrule width .4pt height \lengthdash\unitlength}
\arrowtype=#4 \ypos=#2
\ifnum\arrowtype<0 \advance\arrowtype by 7 \fi
\ifcase\arrowtype
\or \advance\ypos by \arrowlength \advance\ypos by -40
    \put(#1,\ypos){\vector(0,1){\lengthdash}}
    \advance\ypos by -40
    \put(#1,\ypos){\vector(0,1){\lengthdash}}
\or \advance\ypos by 10
    \put(#1,\ypos){\vector(0,1){\lengthdash}}
    \advance\ypos by \arrowlength \advance\ypos by -40
    \put(#1,\ypos){\vector(0,1){\lengthdash}}
\or \advance\ypos by \arrowlength \advance\ypos by -40
    \put(#1,\ypos){\vector(0,1){\lengthdash}}
\or \advance\ypos by 10
    \put(#1,\ypos){\vector(0,-1){\lengthdash}}
\or \advance\ypos by 10
    \put(#1,\ypos){\vector(0,-1){\lengthdash}}
    \advance\ypos by \arrowlength \advance\ypos by -40
    \put(#1,\ypos){\vector(0,-1){\lengthdash}}
\or \advance\ypos by 10
    \put(#1,\ypos){\vector(0,-1){\lengthdash}}
    \advance\ypos by 40
    \put(#1,\ypos){\vector(0,-1){\lengthdash}}
\fi
}}

\def\puthmorphism(#1,#2)[#3`#4`#5]#6#7#8{{%
\xpos #1
\ypos #2
\width #6
\arrowlength #6
\arrowtype=#7
\putbox(\xpos,\ypos){#3\vphantom{#4}}%
{\advance \xpos by\arrowlength
\putbox(\xpos,\ypos){\vphantom{#3}#4}}%
\horsize{\tempcounta}{#3}%
\horsize{\tempcountb}{#4}%
\divide \tempcounta by2
\divide \tempcountb by2
\advance \tempcounta by30
\advance \tempcountb by30
\advance \xpos by\tempcounta
\advance \arrowlength by-\tempcounta
\advance \arrowlength by-\tempcountb
\putvector(\xpos,\ypos)(1,0)\arrowlength\arrowtype
\divide \arrowlength by2
\advance \xpos by\arrowlength
\vertsize{\tempcounta}{#5}%
\divide\tempcounta by2
\advance \tempcounta by20
\if a#8 %
   \advance \ypos by\tempcounta
   \putbox(\xpos,\ypos){#5}%
\else
   \advance \ypos by-\tempcounta
   \putbox(\xpos,\ypos){#5}%
\fi}}

\def\putvmorphism(#1,#2)[#3`#4`#5]#6#7#8{{%
\xpos #1
\ypos #2
\arrowlength #6
\arrowtype #7
\settowidth{\xlen}{$#5$}%
\putbox(\xpos,\ypos){#3}%
{\advance \ypos by-\arrowlength
\putbox(\xpos,\ypos){#4}}%
{\advance\arrowlength by-140
\advance \ypos by-70
\ifdim\xlen>0pt
   \if m#8%
      \putsplitvector(\xpos,\ypos)\arrowlength\arrowtype
   \else
   \putvector(\xpos,\ypos)(0,-1)\arrowlength\arrowtype
   \fi
\else
   \putvector(\xpos,\ypos)(0,-1)\arrowlength\arrowtype
\fi}%
\ifdim\xlen>0pt
   \divide \arrowlength by2
   \advance\ypos by-\arrowlength
   \if l#8%
      \advance \xpos by-40
      \putrbox(\xpos,\ypos){#5}%
   \else\if r#8%
      \advance \xpos by40
      \putlbox(\xpos,\ypos){#5}%
   \else
      \putbox(\xpos,\ypos){#5}%
   \fi\fi
\fi
}}

\def\putsquarep<#1>(#2)[#3;#4`#5`#6`#7]{{%
\setsqparms[#1]%
\setpos(#2)%
\settokens[#3]%
\puthmorphism(\xpos,\ypos)[\tokenc`\tokend`{#7}]{\width}{\arrowtyped}b%
\advance\ypos by \height
\puthmorphism(\xpos,\ypos)[\tokena`\tokenb`{#4}]{\width}{\arrowtypea}a%
\putvmorphism(\xpos,\ypos)[``{#5}]{\height}{\arrowtypeb}l%
\advance\xpos by \width
\putvmorphism(\xpos,\ypos)[``{#6}]{\height}{\arrowtypec}r%
}}

\def\putsquare{\@ifnextchar <{\putsquarep}{\putsquarep%
   <\arrowtypea`\arrowtypeb`\arrowtypec`\arrowtyped;\width`\height>}}
\def\square{\@ifnextchar< {\squarep}{\squarep
   <\arrowtypea`\arrowtypeb`\arrowtypec`\arrowtyped;\width`\height>}}
\def\squarep<#1>[#2`#3`#4`#5;#6`#7`#8`#9]{{%       %     #2------>#3
\setsqparms[#1]%                                   %      |       |
\diagram%                                          %      |       |
\putsquarep<\arrowtypea`\arrowtypeb`\arrowtypec`%  %    #7|       |#8
\arrowtyped;\width`\height>%                       %      |       |
(0,0)[#2`#3`#4`{#5};#6`#7`#8`{#9}]%                %      |       |
\enddiagram%                                       %      v       v
}}                                                 %     #4------>#5
\def\putptrianglep<#1>(#2,#3)[#4`#5`#6;#7`#8`#9]{{%
\settriparms[#1]%
\xpos=#2 \ypos=#3
\advance\ypos by \height
\puthmorphism(\xpos,\ypos)[#4`#5`{#7}]{\height}{\arrowtypea}a%
\putvmorphism(\xpos,\ypos)[`#6`{#8}]{\height}{\arrowtypeb}l%
\advance\xpos by\height
\putmorphism(\xpos,\ypos)(-1,-1)[``{#9}]{\height}{\arrowtypec}r%
}}

\def\putptriangle{\@ifnextchar <{\putptrianglep}{\putptrianglep
   <\arrowtypea`\arrowtypeb`\arrowtypec;\height>}}
\def\ptriangle{\@ifnextchar <{\ptrianglep}{\ptrianglep
   <\arrowtypea`\arrowtypeb`\arrowtypec;\height>}}
\def\ptrianglep<#1>[#2`#3`#4;#5`#6`#7]{{%%    %      #2----->#3
\settriparms[#1]%                             %      |      /
\diagram%                                     %      |     /
\putptrianglep<\arrowtypea`\arrowtypeb`%      %    #6|    /#7
\arrowtypec;\height>%                         %      |   /
(0,0)[#2`#3`#4;#5`#6`{#7}]%                   %      |  /
\enddiagram%%                                 %      v v
}}                                            %      #4

\def\putqtrianglep<#1>(#2,#3)[#4`#5`#6;#7`#8`#9]{{%
\settriparms[#1]%
\xpos=#2 \ypos=#3
\advance\ypos by\height
\puthmorphism(\xpos,\ypos)[#4`#5`{#7}]{\height}{\arrowtypea}a%
\putmorphism(\xpos,\ypos)(1,-1)[``{#8}]{\height}{\arrowtypeb}l%
\advance\xpos by\height
\putvmorphism(\xpos,\ypos)[`#6`{#9}]{\height}{\arrowtypec}r%
}}

\def\putqtriangle{\@ifnextchar <{\putqtrianglep}{\putqtrianglep
   <\arrowtypea`\arrowtypeb`\arrowtypec;\height>}}
\def\qtriangle{\@ifnextchar <{\qtrianglep}{\qtrianglep
   <\arrowtypea`\arrowtypeb`\arrowtypec;\height>}}
\def\qtrianglep<#1>[#2`#3`#4;#5`#6`#7]{{%%    %        #2----->#3
\settriparms[#1]%                             %         \      |
\width=\height                                %          \     |
\diagram%                                     %         #6\    |#7
\putqtrianglep<\arrowtypea`\arrowtypeb`%      %            \   |
\arrowtypec;\height>%                         %             \  |
(0,0)[#2`#3`#4;#5`#6`{#7}]%                   %              v v
\enddiagram%%                                 %               #4
}}

\def\putdtrianglep<#1>(#2,#3)[#4`#5`#6;#7`#8`#9]{{%
\settriparms[#1]%
\xpos=#2 \ypos=#3
\puthmorphism(\xpos,\ypos)[#5`#6`{#9}]{\height}{\arrowtypec}b%
\advance\xpos by \height \advance\ypos by\height
\putmorphism(\xpos,\ypos)(-1,-1)[``{#7}]{\height}{\arrowtypea}l%
\putvmorphism(\xpos,\ypos)[#4``{#8}]{\height}{\arrowtypeb}r%
}}

\def\putdtriangle{\@ifnextchar <{\putdtrianglep}{\putdtrianglep
   <\arrowtypea`\arrowtypeb`\arrowtypec;\height>}}
\def\dtriangle{\@ifnextchar <{\dtrianglep}{\dtrianglep
   <\arrowtypea`\arrowtypeb`\arrowtypec;\height>}}
\def\dtrianglep<#1>[#2`#3`#4;#5`#6`#7]{{%%    %                  / |
\settriparms[#1]%                             %                 /  |
\width=\height                                %              #5/   |#6
\diagram%                                     %               /    |
\putdtrianglep<\arrowtypea`\arrowtypeb`%      %              /     |
\arrowtypec;\height>%                         %             v      v
(0,0)[#2`#3`#4;#5`#6`{#7}]%                   %            #3----->#4
\enddiagram%%                                 %                #7
}}

\def\putbtrianglep<#1>(#2,#3)[#4`#5`#6;#7`#8`#9]{{%
\settriparms[#1]%
\xpos=#2 \ypos=#3
\puthmorphism(\xpos,\ypos)[#5`#6`{#9}]{\height}{\arrowtypec}b%
\advance\ypos by\height
\putmorphism(\xpos,\ypos)(1,-1)[``{#8}]{\height}{\arrowtypeb}r%
\putvmorphism(\xpos,\ypos)[#4``{#7}]{\height}{\arrowtypea}l%
}}

\def\putbtriangle{\@ifnextchar <{\putbtrianglep}{\putbtrianglep
   <\arrowtypea`\arrowtypeb`\arrowtypec;\height>}}
\def\btriangle{\@ifnextchar <{\btrianglep}{\btrianglep
   <\arrowtypea`\arrowtypeb`\arrowtypec;\height>}}
\def\btrianglep<#1>[#2`#3`#4;#5`#6`#7]{{%%   %              | \
\settriparms[#1]%                            %              |  \
\width=\height                               %            #5|   \#6
\diagram%                                    %              |    \
\putbtrianglep<\arrowtypea`\arrowtypeb`%     %              |     \
\arrowtypec;\height>%                        %              v      v
(0,0)[#2`#3`#4;#5`#6`{#7}]%                  %              #3----->#4
\enddiagram%%                                %                 #7
}}

\def\putAtrianglep<#1>(#2,#3)[#4`#5`#6;#7`#8`#9]{{%
\settriparms[#1]%
\xpos=#2 \ypos=#3
{\multiply \height by2
\puthmorphism(\xpos,\ypos)[#5`#6`{#9}]{\height}{\arrowtypec}b}%
\advance\xpos by\height \advance\ypos by\height
\putmorphism(\xpos,\ypos)(-1,-1)[#4``{#7}]{\height}{\arrowtypea}l%
\putmorphism(\xpos,\ypos)(1,-1)[``{#8}]{\height}{\arrowtypeb}r%
}}

\def\putAtriangle{\@ifnextchar <{\putAtrianglep}{\putAtrianglep
   <\arrowtypea`\arrowtypeb`\arrowtypec;\height>}}
\def\Atriangle{\@ifnextchar <{\Atrianglep}{\Atrianglep
   <\arrowtypea`\arrowtypeb`\arrowtypec;\height>}}
\def\Atrianglep<#1>[#2`#3`#4;#5`#6`#7]{{%%         %         /   \
\settriparms[#1]%                                  %        /     \
\width=\height                                     %     #5/       \#6
\diagram%                                          %      /         \
\putAtrianglep<\arrowtypea`\arrowtypeb`%           %     /           \
\arrowtypec;\height>%                              %    v             v
(0,0)[#2`#3`#4;#5`#6`{#7}]%                        %   #3------------>#4
\enddiagram%%                                      %          #7
}}

\def\putAtrianglepairp<#1>(#2)[#3;#4`#5`#6`#7`#8]{{%
\settripairparms[#1]%
\setpos(#2)%
\settokens[#3]%
\puthmorphism(\xpos,\ypos)[\tokenb`\tokenc`{#7}]{\height}{\arrowtyped}b%
\advance\xpos by\height
\puthmorphism(\xpos,\ypos)[\phantom{\tokenc}`\tokend`{#8}]%
{\height}{\arrowtypee}b%
\advance\ypos by\height
\putmorphism(\xpos,\ypos)(-1,-1)[\tokena``{#4}]{\height}{\arrowtypea}l%
\putvmorphism(\xpos,\ypos)[``{#5}]{\height}{\arrowtypeb}m%
\putmorphism(\xpos,\ypos)(1,-1)[``{#6}]{\height}{\arrowtypec}r%
}}

\def\putAtrianglepair{\@ifnextchar <{\putAtrianglepairp}{\putAtrianglepairp%
   <\arrowtypea`\arrowtypeb`\arrowtypec`\arrowtyped`\arrowtypee;\height>}}
\def\Atrianglepair{\@ifnextchar <{\Atrianglepairp}{\Atrianglepairp%
   <\arrowtypea`\arrowtypeb`\arrowtypec`\arrowtyped`\arrowtypee;\height>}}

\def\Atrianglepairp<#1>[#2;#3`#4`#5`#6`#7]{{%           %  #2a
\settripairparms[#1]%                         %           / | \
\settokens[#2]%                               %          /  |  \
\width=\height                                %       #3/  #4   \#5
\diagram%                                     %        /    |    \
\putAtrianglepairp                            %       /     |     \
<\arrowtypea`\arrowtypeb`\arrowtypec`%        %      v      v      v
\arrowtyped`\arrowtypee;\height>%             %     #2b---->#2c---->#2d
(0,0)[{#2};#3`#4`#5`#6`{#7}]%                 %         #6     #7
\enddiagram%%
}}

\def\putVtrianglep<#1>(#2,#3)[#4`#5`#6;#7`#8`#9]{{%
\settriparms[#1]%
\xpos=#2 \ypos=#3
\advance\ypos by\height
{\multiply\height by2
\puthmorphism(\xpos,\ypos)[#4`#5`{#7}]{\height}{\arrowtypea}a}%
\putmorphism(\xpos,\ypos)(1,-1)[`#6`{#8}]{\height}{\arrowtypeb}l%
\advance\xpos by\height
\advance\xpos by\height
\putmorphism(\xpos,\ypos)(-1,-1)[``{#9}]{\height}{\arrowtypec}r%
}}

\def\putVtriangle{\@ifnextchar <{\putVtrianglep}{\putVtrianglep
   <\arrowtypea`\arrowtypeb`\arrowtypec;\height>}}
\def\Vtriangle{\@ifnextchar <{\Vtrianglep}{\Vtrianglep
   <\arrowtypea`\arrowtypeb`\arrowtypec;\height>}}
\def\Vtrianglep<#1>[#2`#3`#4;#5`#6`#7]{{%%     %        #2------------->#3
\settriparms[#1]%                              %         \             /
\width=\height                                 %          \           /
\diagram%                                      %         #6\         /#7
\putVtrianglep<\arrowtypea`\arrowtypeb`%       %            \       /
\arrowtypec;\height>%                          %             \     /
(0,0)[#2`#3`#4;#5`#6`{#7}]%                    %              v   v
\enddiagram%%                                  %               #4
}}

\def\putVtrianglepairp<#1>(#2)[#3;#4`#5`#6`#7`#8]{{
\settripairparms[#1]%
\setpos(#2)%
\settokens[#3]%
\advance\ypos by\height
\putmorphism(\xpos,\ypos)(1,-1)[`\tokend`{#6}]{\height}{\arrowtypec}l%
\puthmorphism(\xpos,\ypos)[\tokena`\tokenb`{#4}]{\height}{\arrowtypea}a%
\advance\xpos by\height
\puthmorphism(\xpos,\ypos)[\phantom{\tokenb}`\tokenc`{#5}]%
{\height}{\arrowtypeb}a%
\putvmorphism(\xpos,\ypos)[``{#7}]{\height}{\arrowtyped}m%
\advance\xpos by\height
\putmorphism(\xpos,\ypos)(-1,-1)[``{#8}]{\height}{\arrowtypee}r%
}}

\def\putVtrianglepair{\@ifnextchar <{\putVtrianglepairp}{\putVtrianglepairp%
    <\arrowtypea`\arrowtypeb`\arrowtypec`\arrowtyped`\arrowtypee;\height>}}
\def\Vtrianglepair{\@ifnextchar <{\Vtrianglepairp}{\Vtrianglepairp%
    <\arrowtypea`\arrowtypeb`\arrowtypec`\arrowtyped`\arrowtypee;\height>}}
\def\Vtrianglepairp<#1>[#2;#3`#4`#5`#6`#7]{{%  %  #2a---->#2b---->#2c
\settripairparms[#1]%                          %   \      |      /
\settokens[#2]%                                %    \     |     /
\diagram%                                      %   #5\   #6    /#7
\putVtrianglepairp                             %      \   |   /
<\arrowtypea`\arrowtypeb`\arrowtypec`%         %       \  |  /
\arrowtyped`\arrowtypee;\height>%              %        v v v
(0,0)[{#2};#3`#4`#5`#6`{#7}]%                  %         #2d
\enddiagram%%
}}

\def\putCtrianglep<#1>(#2,#3)[#4`#5`#6;#7`#8`#9]{{%
\settriparms[#1]%
\xpos=#2 \ypos=#3
\advance\ypos by\height
\putmorphism(\xpos,\ypos)(1,-1)[``{#9}]{\height}{\arrowtypec}l%
\advance\xpos by\height
\advance\ypos by\height
\putmorphism(\xpos,\ypos)(-1,-1)[#4`#5`{#7}]{\height}{\arrowtypea}l%
{\multiply\height by 2
\putvmorphism(\xpos,\ypos)[`#6`{#8}]{\height}{\arrowtypeb}r}%
}}

\def\putCtriangle{\@ifnextchar <{\putCtrianglep}{\putCtrianglep
    <\arrowtypea`\arrowtypeb`\arrowtypec;\height>}}
\def\Ctriangle{\@ifnextchar <{\Ctrianglep}{\Ctrianglep
    <\arrowtypea`\arrowtypeb`\arrowtypec;\height>}}
\def\Ctrianglep<#1>[#2`#3`#4;#5`#6`#7]{{%%   %                / |
\settriparms[#1]%                            %             #5/  |
\width=\height                               %              /   |
\diagram%                                    %             v    |
\putCtrianglep<\arrowtypea`\arrowtypeb`%     %           #3     |#6
\arrowtypec;\height>%                        %             \    |
(0,0)[#2`#3`#4;#5`#6`{#7}]%                  %            #7\   |
\enddiagram%%                                %               \  |
}}                                           %                v v
\def\putDtrianglep<#1>(#2,#3)[#4`#5`#6;#7`#8`#9]{{%
\settriparms[#1]%
\xpos=#2 \ypos=#3
\advance\xpos by\height \advance\ypos by\height
\putmorphism(\xpos,\ypos)(-1,-1)[``{#9}]{\height}{\arrowtypec}r%
\advance\xpos by-\height \advance\ypos by\height
\putmorphism(\xpos,\ypos)(1,-1)[`#5`{#8}]{\height}{\arrowtypeb}r%
{\multiply\height by 2
\putvmorphism(\xpos,\ypos)[#4`#6`{#7}]{\height}{\arrowtypea}l}%
}}

\def\putDtriangle{\@ifnextchar <{\putDtrianglep}{\putDtrianglep
    <\arrowtypea`\arrowtypeb`\arrowtypec;\height>}}
\def\Dtriangle{\@ifnextchar <{\Dtrianglep}{\Dtrianglep
   <\arrowtypea`\arrowtypeb`\arrowtypec;\height>}}
\def\Dtrianglep<#1>[#2`#3`#4;#5`#6`#7]{{%%  %          | \
\settriparms[#1]%                           %          |  \#6
\width=\height                              %          |   \
\diagram%                                   %          |    v
\putDtrianglep<\arrowtypea`\arrowtypeb`%    %        #5|    #3
\arrowtypec;\height>%                       %          |    /
(0,0)[#2`#3`#4;#5`#6`{#7}]%                 %          |   /#7
\enddiagram%%                               %          |  /
}}                                          %          v v
\def\setrecparms[#1`#2]{\width=#1 \height=#2}%
\def\recursep<#1`#2>[#3;#4`#5`#6`#7`#8]{{%
\width=#1 \height=#2
\settokens[#3]
\settowidth{\tempdimen}{$\tokena$}
\ifdim\tempdimen=0pt
  \savebox{\tempboxa}{\hbox{$\tokenb$}}%
  \savebox{\tempboxb}{\hbox{$\tokend$}}%
  \savebox{\tempboxc}{\hbox{$#6$}}%
\else
  \savebox{\tempboxa}{\hbox{$\hbox{$\tokena$}\times\hbox{$\tokenb$}$}}%
  \savebox{\tempboxb}{\hbox{$\hbox{$\tokena$}\times\hbox{$\tokend$}$}}%
  \savebox{\tempboxc}{\hbox{$\hbox{$\tokena$}\times\hbox{$#6$}$}}%
\fi
\ypos=\height
\divide\ypos by 2
\xpos=\ypos
\advance\xpos by \width
\bfig
\putCtrianglep<-1`1`1;\ypos>(0,0)[`\tokenc`;#5`#6`{#7}]%
\puthmorphism(\ypos,0)[\tokend`\usebox{\tempboxb}`{#8}]{\width}{-1}b%
\puthmorphism(\ypos,\height)[\tokenb`\usebox{\tempboxa}`{#4}]{\width}{-1}a%
\advance\ypos by \width
\putvmorphism(\ypos,\height)[``\usebox{\tempboxc}]{\height}1r%
\efig
}}

\def\recurse{\@ifnextchar <{\recursep}{\recursep<\width`\height>}}

\def\puttwohmorphisms(#1,#2)[#3`#4;#5`#6]#7#8#9{{%
\puthmorphism(#1,#2)[#3`#4`]{#7}0a
\ypos=#2
\advance\ypos by 20
\puthmorphism(#1,\ypos)[\phantom{#3}`\phantom{#4}`#5]{#7}{#8}a
\advance\ypos by -40
\puthmorphism(#1,\ypos)[\phantom{#3}`\phantom{#4}`#6]{#7}{#9}b
}}

\def\puttwovmorphisms(#1,#2)[#3`#4;#5`#6]#7#8#9{{%
\putvmorphism(#1,#2)[#3`#4`]{#7}0a
\xpos=#1
\advance\xpos by -20
\putvmorphism(\xpos,#2)[\phantom{#3}`\phantom{#4}`#5]{#7}{#8}l
\advance\xpos by 40
\putvmorphism(\xpos,#2)[\phantom{#3}`\phantom{#4}`#6]{#7}{#9}r
}}

\def\puthcoequalizer(#1)[#2`#3`#4;#5`#6`#7]#8#9{{%
\setpos(#1)%
\puttwohmorphisms(\xpos,\ypos)[#2`#3;#5`#6]{#8}11%
\advance\xpos by #8
\puthmorphism(\xpos,\ypos)[\phantom{#3}`#4`#7]{#8}1{#9}
}}

\def\putvcoequalizer(#1)[#2`#3`#4;#5`#6`#7]#8#9{{%
\setpos(#1)%
\puttwovmorphisms(\xpos,\ypos)[#2`#3;#5`#6]{#8}11%
\advance\ypos by -#8
\putvmorphism(\xpos,\ypos)[\phantom{#3}`#4`#7]{#8}1{#9}
}}

\def\putthreehmorphisms(#1)[#2`#3;#4`#5`#6]#7(#8)#9{{%
\setpos(#1) \settypes(#8)
\if a#9 %
     \vertsize{\tempcounta}{#5}%
     \vertsize{\tempcountb}{#6}%
     \ifnum \tempcounta<\tempcountb \tempcounta=\tempcountb \fi
\else
     \vertsize{\tempcounta}{#4}%
     \vertsize{\tempcountb}{#5}%
     \ifnum \tempcounta<\tempcountb \tempcounta=\tempcountb \fi
\fi
\advance \tempcounta by 60
\puthmorphism(\xpos,\ypos)[#2`#3`#5]{#7}{\arrowtypeb}{#9}
\advance\ypos by \tempcounta
\puthmorphism(\xpos,\ypos)[\phantom{#2}`\phantom{#3}`#4]{#7}{\arrowtypea}{#9}
\advance\ypos by -\tempcounta \advance\ypos by -\tempcounta
\puthmorphism(\xpos,\ypos)[\phantom{#2}`\phantom{#3}`#6]{#7}{\arrowtypec}{#9}
}}

\def\setarrowtoks[#1`#2`#3`#4`#5`#6]{%
\def\toka{#1}
\def\tokb{#2}
\def\tokc{#3}
\def\tokd{#4}
\def\toke{#5}
\def\tokf{#6}
}
\def\hex{\@ifnextchar <{\hexp}{\hexp<1000`400>}}
\def\hexp<#1`#2>[#3`#4`#5`#6`#7`#8;#9]{%
\setarrowtoks[#9]
\yext=#2 \advance \yext by #2
\xext=#1 \advance\xext by \yext
\bfig
\putCtriangle<-1`0`1;#2>(0,0)[`#5`;\tokb``\tokd]
\xext=#1 \yext=#2 \advance \yext by #2
\putsquare<1`0`0`1;\xext`\yext>(#2,0)[#3`#4`#7`#8;\toka```\tokf]
\advance \xext by #2
\putDtriangle<0`1`-1;#2>(\xext,0)[`#6`;`\tokc`\toke]
\efig
}